\documentclass[12pt]{amsart}
\usepackage{amssymb}
\usepackage{latexsym}
\usepackage{epsfig}
\usepackage{float}
\usepackage{epstopdf}
\usepackage{graphicx}
\usepackage{tikz}
\usepackage[linktocpage=true]{hyperref}
\usepackage{ulem}
\usepackage[left=3.8cm, right=3.8cm, top=3cm, bottom=3cm]{geometry}
\usepackage{color}
\newcommand{\re}[1]{(\ref{#1})}

\newcommand{\rl}[1]{Lemma~\ref{#1}}
\newcommand{\nrc}[1]{Corollary~\ref{#1}}
\newcommand{\rp}[1]{Proposition~\ref{#1}}
\newcommand{\rt}[1]{Theorem~\ref{#1}}

\hypersetup{ colorlinks   = true, 
urlcolor  = blue, 
linkcolor    = blue, 
citecolor   = blue 
}
\definecolor{indianyellow}{rgb}{0.89, 0.66, 0.34}
\definecolor{goldenpoppy}{rgb}{0.99, 0.76, 0.0}
\definecolor{gamboge}{rgb}{0.89, 0.61, 0.06}
\definecolor{frenchbeige}{rgb}{0.65, 0.48, 0.36}
\definecolor{burntorange}{rgb}{0.8, 0.33, 0.0}

\date{\today}

\def\beq{\begin{equation}}
\def\eeq{\end{equation}}

\newcommand{\cL}{\mathcal}

\newcommand{\Z}{{\mathbb Z}}
\newcommand{\R}{{\mathbb R}}
\newcommand{\Q}{{\mathbb Q}}
\newcommand{\C}{{\mathbb C}}

\newcommand{\N}{{\mathbb N}}

\newcommand{\CA}{{\mathcal A}}

\newcommand{\CH}{{\mathcal H}}

\newcommand{\CM}{{\mathcal M}}
\newcommand{\CN}{{\mathcal N}}

\newcommand{\CP}{{\mathcal P}}
\newcommand{\CS}{{\mathcal S}}

\newcommand{\CU}{{\mathcal U}}

\newcommand{\CO}{{\mathcal O}}

\newcommand{\CV}{{\mathcal V}}

\newtheorem{thm}{Theorem}[section]
\newtheorem{cor}[thm]{Corollary}
\newtheorem{prop}[thm]{Proposition}
\newtheorem{lemma}[thm]{Lemma}
\newtheorem{remark}[thm]{Remark}

\theoremstyle{definition}

\newcommand{\ov}{\overline}

\newcommand{\RE}{\operatorname{Re}}

\newcommand{\gaa}{\gamma}

\newcommand{\var}{\varphi}

\newcommand{\om}{\omega}

\usepackage{enumitem}
\DeclareMathOperator{\imag}{Im}

\begin{document}

\title[]{Geometry of hyperbolic Cauchy-Riemann singularities and KAM-like theory for holomorphic involutions}

\author{Laurent Stolovitch}
\address{CNRS and Laboratoire J.-A. Dieudonn\'e
	U.M.R. 7351, Universit\'e C\^ote d'Azur, Parc Valrose
	06108 Nice Cedex 02, France}
\email{stolo@unice.fr}
\thanks{This work has been supported by the French government throught the ANR grant ANR-15-CE40-0001-03 for the project Bekam and through the UCAJEDI Investments in the Future project managed by the National Research Agency (ANR) with the reference number ANR-15-IDEX-01}

\author{Zhiyan Zhao}
\address{Laboratoire J.-A. Dieudonn\'e
	U.M.R. 7351, Universit\'e C\^ote d'Azur, Parc Valrose
	06108 Nice Cedex 02, France}
\email{zhiyan.zhao@univ-cotedazur.fr}

\begin{abstract}
This article is concerned with the geometry of germs of real analytic surfaces in $(\mathbb{C}^2,0)$ having an isolated Cauchy-Riemann (CR) singularity at the origin. These are perturbations of {\it Bishop quadrics}. There are two kinds of CR singularities stable under perturbation~: {\it elliptic} and {\it hyperbolic}. Elliptic case was studied by Moser-Webster \cite{moser-webster} who showed that such a surface is locally, near the CR singularity, holomorphically equivalent to {\it normal form} from which lots of geometric features can be read off.

In this article we focus on perturbations of {\it hyperbolic} quadrics. As was shown by Moser-Webster \cite{moser-webster}, such a surface can be transformed to a formal {\it normal form} by a formal change of coordinates that may not be holomorphic in any neighborhood of the origin.
	
Given a {\it non-degenerate} real analytic surface $M$ in $(\mathbb{C}^2,0)$ having a {\it hyperbolic} CR singularity at the origin, we prove the existence of a non-constant Whitney smooth family of connected holomorphic curves intersecting $M$ along holomorphic hyperbolas.
This is the very first result concerning hyperbolic CR singularity not equivalent to quadrics.
	
 This is a consequence of a non-standard KAM-like theorem for pair of germs of holomorphic involutions $\{\tau_1,\tau_2\}$ at the origin, a common fixed point. We show that such a pair has large amount of invariant analytic sets biholomorphic to $\{z_1z_2=const\}$ (which is not a torus) in a neighborhood of the origin, and that they are conjugate to restrictions of linear maps on such invariant sets.

\



\end{abstract}

\maketitle


\section{Introduction}
In this article, we are concerned with the local holomorphic invariants of a real analytic submanifold $M$ in $\C^n$.
If the tangent space of $M$ at a point $p_0$ contains a maximal {\it complex} subspace, the dimension $d$ of which does not depend on $p_0$, then we say that $M$ is a {\it Cauchy-Riemann (CR) submanifold}. Since the work of Cartan \cite{cartan} in the 30's, lots of studies were devoted to this geometry (see for instance \cite{chernRealHypersurfacesComplex1974, ber, lamelConvergenceDivergenceFormal2018, kossovskiyConvergentNormalForm2019}).
As ``baby" example, one can consider an open neighborhood $U$ of a point $p_0$ in $\R^n$ in $\C^n$. The local hull of holomorphy of $U$ is the largest open set in $\C^n$ containing $U$ and over which all holomorphic functions defined on $U$ can be holomorphically extended to. It can be shown that, in that case, the hull of holomorphy of $U$ is nothing but $U$. This situation is quite different when considering a neighborhood of a CR singularity, that is a point $p_0$ in the real submanifold $M$ in $\C^n$ such that the maximal complex tangent spaces do not have a constant dimension in any neighborhood of $p_0$.
A real submanifold with a CR singularity must have codimension at least $2$.

The study of real submanifolds with CR singularities was initiated by Bishop \cite{bishop} in his pioneering work, and followed by Moser-Webster \cite{moser-webster}. They considered higher-order analytic perturbations of the elementary models called Bishop quadrics ${\cL Q}_{\gaa}\subset\C^2 $, depending on the
Bishop invariant $0\leq\gaa\leq\infty$~:
 \begin{itemize}
   \item for $0\leq\gaa<\infty$,
   ${\cL Q}_\gaa : z_2=Q_{\gamma}(z_1,\bar z_1):=|z_1|^2+\gaa(z_1^2+\ov z_1^2)$,
   \item for $\gamma=\infty$,
   ${\cL Q}_\infty\colon z_2=z_1^2+\ov z_1^2$.
 \end{itemize}
When $\gaa\neq \frac{1}{2}$, such a surface has an isolated CR singularity at the origin as it is totally real (i.e., $d=0$) everywhere
 but at the origin at which the tangent space is the complex line $\{z_2=0\}$ (i.e., $d=1$).
When $0<\gaa<\frac{1}{2}$, one says that this singularity is {\it elliptic}. In their seminal work, Moser-Webster \cite{moser-webster} considered higher-order analytic perturbations of elliptic quadrics ${\cL Q}_{\gaa}$.
 They proved that such a submanifold is holomorphically equivalent to a {\it normal form}, $z_2=|z_1|^2+(\gamma+\epsilon\RE(z_2)^s)(z_1^2+\bar z_1^2)$, for some $\epsilon\in \{-1,0,1\}$ and $s\in \mathbb{N^*\cup\{\infty\}}$.
 Lots of geometric features can be read off from such a normal form. They also considered $n-$dimensional submanifolds in $\mathbb C^n$ which have a complex tangent at the origin of minimal (positive) dimension. This has been recently extended to CR singularity with {\it maximal complex tangent} by Gong and the first author \cite{stolo-gong1, stolo-gong2}. When $\gaa=0$ (degenerate elliptic case), Moser \cite{moser-zero} constructed a formal power series normal form. Although it is still not known whether such a normal form can be obtained through a convergent transformation, Huang-Yin \cite{huang-yin-zero} did the achievement of obtaining the holomorphic classification of analytic perturbations of ${\cL Q}_0$. Relatively recently, related problems such as flattening \cite{huang-yin-mathann, huang-yin-advmath, huang-fang-gafa} or quadric rigidity \cite{huang-yin-imrn} have been successfully considered by Huang and co-authors. Some results on CR singularities of $k-$dimensional submanifolds in $\mathbb{C}^n$, $k\neq n$ have been obtained by A. Coffman \cite{coffman-cross, coffman-illinois}.

In the so-called {\it hyperbolic} case, i.e., higher-order analytic perturbation of ${\cL Q}_{\gamma}$ with $\gamma>\frac{1}{2}$, not much is known. Moser-Webster \cite{moser-webster} showed that some analytic perturbations of ${\cL Q}_{\gaa}$ may not be holomorphically equivalent to a normal form as in the elliptic case. Forstneri\v{c}-Stout \cite{fs-hyperbolic-polhull} proved that such a perturbation is always polynomially convex near such a hyperbolic CR singularity. Gong \cite{gong-hyperbolic} showed that if the higher-order analytic perturbation of ${\cL Q}_{\gaa}$ is {\it formally} equivalent to ${\cL Q}_{\gaa}$ (i.e., by the mean of formal power series transformation) and if a {\it Diophantine condition} associated to $\gaa$ is satisfied, then the perturbation is actually holomorphically equivalent to the quadric. He also proved the existence of higher-order analytic perturbations of a hyperbolic quadric which are formally equivalent to the hyperbolic quadric but not holomorphically equivalent to it \cite{gong-hyper-div}. On the other hand, Klingenberg \cite{klingenberg} showed that under a similar Diophantine condition, for a given higher-order analytic perturbation $M$ of the quadric, there always exists a holomorphic curve that intersects $M$ along two transverse totally real curves. Both results have been extended in higher dimension in the case of maximal complex tangent \cite{stolo-gong1}.

In both elliptic and hyperbolic cases, the CR singularity is stable under perturbation and is not removable.

The aim of this work is to prove that {\it non-degenerate analytic perturbations} of hyperbolic quadrics, i.e., perturbations which are not formally equivalent to quadrics, contain a large number of analytic hyperbolas.
By this, we mean that there exists a compact set $\cL K\subset \mathbb{R}$ of positive measure such that for all $\om\in \cL K$, there exists a connected holomorphic curve $\cL S_{\om}$ that intersects the non-degenerate analytic perturbation $M$ along two distinguished real analytic curves that are simultaneously holomorphically mapped to the two branches of the real hyperbola $\{\xi\eta=\om\}$ (in a neighborhood of the origin).  We remark that it is elementary that a real analytic curve in the real analytic surface is contained in a holomorphic curve. Having a connected holomorphic curve that intersects $M$ in two distinct real analytic curves is, however, one of main conclusions of this paper.

To do so, we shall develop a new KAM theory (named after Kolmogorov-Arnold-Moser \cite{Kolmogorov, arnold-kol, moser-anneau}) for a pair of (germs of) holomorphic involutions in a neighborhood of a fixed point (say $0$) in $\mathbb{C}^2$, which is swapped by conjugacy with some anti-holomorphic involution. Initially, KAM theory was conceived as an answer to the fundamental problem arising in Dynamical Systems and in particular in Celestial Mechanics \cite{fejoz-herman,chierchia-nbody}. It can be formulated as follows~: Given a {\it completely integrable} Hamiltonian dynamical system written in {\it action-angle coordinates} $(\theta, I)\in\mathbb{T}^n\times \mathbb{R}^n$ of the form  $\dot\theta=\omega(I),\; \dot I=0$, where $\omega$ denotes an analytic function. For each $I_0$, the manifold $\mathbb{T}^n\times\{I_0\}$ is invariant and the motion on it is a constant rotation of angle $\omega(I_0)$. In the nature, these systems are rather rare but one encounters small perturbations of them under the form $(*)\; \dot\theta=\omega(I)+ \epsilon f(I,\theta),\; \dot I= \epsilon g(I,\theta)$ with $f,g$ analytic functions and $\epsilon$ a small number. Essentially, KAM theorem states that if the system is {\it non-degenerate} in some sense then there exists a large (in measure) compact set $\mathcal{K}$ such that for all $I\in \mathcal{K}$, the system $(*)$ has an invariant manifold which is diffeomorphic to a torus the dynamical system on which is conjugated to the rotation of angle $\tilde \omega(I)$ on that torus. In some sense, a lot of the invariant tori $\mathbb{T}^n\times\{I\}$ ``survive" under a small non-degenerate perturbation.

As mentioned earlier, to obtain a connected holomorphic curve intersecting the real surface is a main result. In our context, we shall go one step further by proving, for a sufficiently small $R>0$, the existence of a compact set $\cL{O}_{\infty}(R)\subset]- R^2, R^2[$ of positive measure such that for each $\om\in\cL{O}_{\infty}(R)$, there exists an invariant connected complex submanifold $\tilde {\cL S}_\om$ in $\Delta_2(0,R^{\frac12}):=\{(\xi,\eta)\in\C^2: |\xi|,|\eta|<R^{\frac12}\}$, which is the image of the connected holomorphic manifold ${\cL C}_\omega^R:=\{(\xi,\eta)\in\C^2:  \xi\eta=\om, \, |\xi|, |\eta|<R\}$ by a biholomorphism $\Psi_\om :{\cL C}_\omega^R\rightarrow \tilde{\cL S}_\om$. The inverse of the latter, $\Psi_\om^{-1}$, conjugates the restrictions of nonlinear involutions to $\tilde{\mathcal{S}}_{\om}$ to the restrictions to ${\cL C}_\omega^R$ of linear ones. We emphasize that ${\cL C}_\omega^R \cap \Delta_2(0, a\sqrt{|\om|})$ contains the graph $\zeta\mapsto (\zeta, \frac{\om}{\zeta})$ over the annulus $\frac{\sqrt{|\om|}}{a}<|\zeta|<\sqrt{|\om|}a$. This KAM-like result is non-standard as one does not expect to obtain {\it invariant tori} as in \cite{bost-bourb, eliasson-lindstedt, russmann-weak,BBP13} but different kind of invariant manifolds of the form $\{z_1z_2=\om\}$ (in a neighborhood of the origin; when $\om\neq 0$). The role played by the rotation is played by linear involutions.
In a similar spirit, but in a different context, a KAM-like theory was obtained by the first author for germs of holomorphic vector fields at a fixed point \cite{Stolo-kam}.
We emphasize that the KAM-like statement in this paper is different from an apparently similar {\it real} problem for which one obtains a lot of invariant tori. The main achievements in this direction are due to Sevryuk \cite{sevryuk-lnm, sevryuk-2a, sevryuk-2b} near an elliptic fixed point. The non-standard hyperbola character of our KAM-like result near an elliptic fixed point of reversible holomorphic mappings unveils new unexpected difficulties.

{\it Hard implicit function theorem, Nash-Moser theorem, Newton Scheme or KAM process} are various names in the literature that stand for ``rapid iteration scheme" usually needed to solve functional equations in Fr\'echet spaces \cite{hamilton}. This appears in particular in conjugacy problem to normal forms of vector fields at a fixed point \cite{bruno,stolo-ihes}, of  interval exchange maps \cite{mmy-nonstandardKAM} or in reducibility problems of quasi-periodic cocycles \cite{eliasson-floquet,hou-you,AFK11}, the latter being related to spectral theory as well.\\

\noindent

{\bf Acknowledgment:} The authors thank X. Gong for his interest, stimulating discussions that help them to substantially improve the exposition of the main results. The authors would also like to thank L. Lempert and M. Procesi for their interests and sharp comments which help them to improve the text.
The authors thank also the anonymous referees for the careful review and helpful suggestions.

\section{Main results}\label{sec_main}

We shall here summarize some statements of \cite{moser-webster}.
Let us consider Bishop's {\it hyperbolic} quadric, a real quadratic surface in $\C^2$ given by
$${\cL Q}_{\gaa}:z_2=Q_{\gaa}(z_1,\bar z_1)=|z_1|^2+\gaa(z_1^2+\ov z_1^2),\quad \gaa>\frac{1}{2}.$$
Let $M$ be a higher-order analytic perturbation of ${\cL Q}_{\gaa}$ given by
\beq\label{defM}
M: z_2=Q_{\gaa}(z_1,\bar z_1)+ f(z_1,\bar z_1) ,\quad f(z_1,\bar z_1) = O^3(z_1,\bar z_1).
\eeq
To such a real surface, one associates a local dynamical system in $(\Bbb C^2,0)$, $\{\tau_1^o,\tau_2^o,\rho\}$, where $\tau_1^o,\tau_2^o$ are local holomorphic involutions fixing $0$, $\rho$ is an anti-holomorphic involution. They satisfy $\tau_j^o\circ\tau_j^o={\rm Id} $ and $\tau_2^o=\rho\circ\tau_1^o\circ\rho$. Moser-Webster's construction goes as follow. We first "complexify" $M$ as a complex surface $\cL M$ in $\C^4$ by considering two new complex independent variables $w_1$, $w_2$, playing the role of $\bar z_1$, $\bar z_2$ respectively. In these coordinates,
$$
\cL M:\left\{\begin{array}{ccc}
z_2 & = & Q_{\gaa}(z_1,w_1)+f(z_1,w_1)\\
w_2 & = & Q_{\gaa}(z_1,w_1)+\bar f(w_1,z_1).\\
\end{array}\right.
$$
There are two natural holomorphic mappings $\pi_i: (\C^4,0)\cap \cL M\rightarrow (\C^2,0)$, $i=1,2$, defined as $\pi_1(z,w)=z$ and $\pi_2(z,w)=w$. It happens that these are $2-1$ branched coverings. Each mapping $\tau_i^o$ is defined to be the deck transformation (different from identity) of $\pi_i$, that is $$\pi_1(\tau_1^o(z,w))=z, \quad \pi_2(\tau_2^o(z,w))=w.$$
For instance, $\tau_1^o$ can be regarded as a mapping defined as $\tau_1^o(z_1,w_1)=(z_1, \phi_1(z_1,w_1))$ such that $$Q_{\gaa}(z_1,\phi_1(z_1,w_1))+f(z_1,\phi_1(z_1,w_1))=Q_{\gaa}(z_1,w_1)+f(z_1,w_1).$$
The linear part $T_1$ of the mapping $\tau_1^o$ at the fixed point $0$, is obtained by solving the equation $Q_{\gaa}(z_1,T_1(z_1,w_1))=Q_{\gaa}(z_1,w_1)$. An immediate computation shows that $T_1(z_1,w_1)=(z_1,-\gamma^{-1}z_1-w_1)$. In good local holomorphic coordinates $(\xi,\eta)$, $T_1$ is rewritten as $T_1(\xi,\eta)=(\delta\eta,\delta^{-1}\xi)$ for some complex number $\delta$.

Such a triple $\{\tau_1^o,\tau_2^o,\rho\}$ completely characterizes the holomorphic equivalent class of the real surface $M$ (cf. \cite{moser-webster}[Proposition 1.1] or \cite{stolo-gong1}[Proposition 2.8]).
It is also useful to consider the germ of biholomorphism $\sigma_o:=\tau_1^o\circ\tau_2^o$. In good local holomorphic coordinates $(\xi,\eta)$, we have $\rho(\xi,\eta)=(\bar \xi,\bar\eta)$,
\begin{align}
	\tau^o_1(\xi,\eta)&=\left(\begin{array}{l}
		e^{\frac{\rm i}2\lambda}\eta+p^o(\xi,\eta) \\
		e^{-\frac{\rm i}2\lambda}\xi+q^o(\xi,\eta)
	\end{array}
	\right),\label{tau1}\\
	\tau^o_2(\xi,\eta)&=\left(\rho\circ\tau^o_1\circ\rho\right)(\xi,\eta)=\left(\begin{array}{l}
		e^{-\frac{\rm i}2\lambda}\eta+\bar p^o(\xi,\eta) \\
		e^{\frac{\rm i}2\lambda}\xi+\bar q^o(\xi,\eta)\label{tau2}
	\end{array}
	\right),
\end{align}
where $\bar h$ denotes $\bar h(\xi,\eta):=\sum_{k,l\geq 0} \bar{\breve h}_{k,l} \xi^k\eta^l $ if $h(\xi,\eta):=\sum_{k,l\geq 0} \breve h_{k,l} \xi^k\eta^l$. We also have
\begin{equation}\label{sigma_o}
\sigma_o(\xi,\eta)=\left(\begin{array}{l}
                           \mu\xi+f^o(\xi,\eta) \\
                           \mu^{-1}\eta+g^o(\xi,\eta)
                          \end{array}
\right),\quad \mu=e^{\rm i\lambda},\quad |e^{\rm i\lambda}|=1.
\end{equation}

Here, $e^{\frac{\rm i}2\lambda}$, $e^{-\frac{\rm i}2\lambda}$ are the roots of the quadratic equation $\gamma X^2-X+\gamma=0$ and $p^o$, $q^o$, $f^o$, $g^o$ are germs of holomorphic functions of order $\geq 2$ at the origin (i.e., the functions and their first-order derivatives vanish at $0$).
In the case $M={\cL Q}_{\gamma}$, $\tau^o_1$, $\tau^o_2$ are the linear involutions
$$\tau_1^o(\xi,\eta) =\left(\begin{array}{l}
e^{\frac{\rm i}2\lambda}\eta \\
e^{-\frac{\rm i}2\lambda}\xi
\end{array}
\right),\quad \tau_2^o(\xi,\eta) =\left(\begin{array}{l}
e^{-\frac{\rm i}2\lambda}\eta \\
e^{\frac{\rm i}2\lambda}\xi
\end{array}
\right).$$

In the sequel, we shall assume that the submanifold $M$ (or their associated involutions $\tau^o_1$, $\tau^o_2$) is {\it non-exceptional}, meaning that $e^{\frac{\rm i}2\lambda}$ is not a root of unity. In this case, Moser-Webster showed (cf. \cite{moser-webster}[Lemma 3.2, Theorem 3.4]) that there exists a formal transformation $\hat \Psi$ satisfying $\hat\Psi\circ \rho= \rho\circ\hat \Psi$ such that
\begin{eqnarray}
	\hat\tau_{1}:=(\hat\Psi^{-1}\circ \tau^o_1\circ\hat\Psi)(\xi,\eta)&=& \left(\begin{array}{l}
		\Lambda(\xi\eta)\eta \\
		\Lambda^{-1}(\xi\eta)\xi
	\end{array}
	\right),\label{formal_tau1}\\
	\hat\tau_{2}:=(\hat\Psi^{-1}\circ \tau^o_2\circ\hat\Psi)(\xi,\eta)&=& \left(\begin{array}{l}
		\Lambda(\xi\eta)^{-1}\eta \\
		\Lambda(\xi\eta)\xi
	\end{array}
	\right),\label{formal_tau2}\\
	\hat\sigma:=(\hat\Psi^{-1}\circ \sigma_o\circ\hat\Psi)(\xi,\eta) &=&\left(\begin{array}{l}
		\hat M(\xi\eta)\xi \\
		\hat M(\xi\eta)^{-1}\eta
	\end{array}
	\right). \label{formal_sigma}
\end{eqnarray}
Here, $\Lambda(z)$ and $\hat M(z)$ are formal power series of the one-dimensional variable $z$ and satisfy:
$$
\Lambda(z)\bar\Lambda(z)=1,\quad \hat M(z)= \Lambda(z)^2,\quad \Lambda(0)=e^{\frac{\rm i}2\lambda},\quad \hat M(0)=\mu.
$$
The maps $\hat \tau_j$ and $\hat \sigma$ are called formal {\it normal form}. Furthermore, the pair $\{\hat\tau_1,\hat \tau_2\}$ is said to be {\it formally integrable}. It would have been called {\it integrable} over a domain in $\mathbb{C}^2$ if $\Lambda$ was holomorphic in that domain. The map $\hat\Psi$ is called the {\it normalizing transformation}.
Contrary to the elliptic case, one cannot expect the normalizing transformation to converge in a neighborhood of the origin.
This is due to the presence of {\it small divisors} (we recall that $\{|\mu^{k}-1|\}_{k\in\N^*}$ accumulate at the origin when $|\mu|=1$) as emphasized in {\cite{moser-webster}[Section 6 (b)].

If $\Lambda(z)=\Lambda(0)$, then $\{\tau^o_1,\tau^o_2\}$ is formally linearizable by a formal transformation that commutes with $\rho$. Hence, the submanifold is formally equivalent to the quadric ${\cL Q}_{\gaa}$. Gong's theorem \cite{gong-hyperbolic} asserts that, if a {\it Diophantine condition} is satisfied, i.e., there exist $r,c>0$, such that for $k\in \mathbb N^*$,  $|\mu^{k}-1|\geq \frac{c}{k^{r}}$, then the submanifold is actually holomorphically equivalent to the quadric ${\cL Q}_{\gaa}$ near the origin.

In what follows, we shall focus on the {\it non-degenerate} case, i.e., we assume that $ \Lambda(z)\neq\Lambda(0)$ and assume that $s$ is the smallest positive integer $l$ such that $\Lambda^{(l)}(0)\neq 0$. We can normalize $\frac{\Lambda^{(s)}(0)}{s!}=1$.

\subsection{KAM-like theorem for reversible holomorphic maps}
We 
assume that $\tau^o_1, \tau^o_2, \sigma_o$ are defined
in $\{|\xi|,|\eta|<r\}$ for some $0<r<\frac1{4}$ as in (\ref{tau1}) -- (\ref{sigma_o}), and
\begin{enumerate}[label=(\Alph*)]
	\item \label{nondeg}
$\lambda\in[0,4\pi[$ with $\frac{\lambda}{\pi}\in\R\setminus\Q$,
	\item \label{nd2} $p^o$ and $q^o$ are convergent power series on $\{|\xi|,|\eta|<r\}$ of order $\geq 2$, i.e.,
$$p^o(\xi,\eta)=\sum_{l+j\geq2\atop{l,j\geq 0}}\breve p^o_{l,j}\xi^l\eta^j,\quad q^o(\xi,\eta)=\sum_{l+j\geq2\atop{l,j\geq 0}}\breve q^o_{l,j}\xi^l\eta^j,$$
with coefficients $\breve p^o_{l,j}, \ \breve q^o_{l,j}\in\C$.
\end{enumerate}
It is easy to verify that $\sigma_o$ is {\it reversible} w.r.t. the involution $\rho$, i.e.,
$\sigma_o^{-1}=\rho\circ \sigma_o \circ \rho$.

As above, let $\hat\Psi$ be the unique normalized formal transformation together with the formal power series $\Lambda=\Lambda(z)$.
We assume that $\Lambda(z)$ is not constant. Let $s\in\N^*$ be the smallest positive integer such that $\Lambda^{(s)}(0)\neq 0$.
More precisely, we assume that
\beq\label{nd3}
\Lambda(z)=e^{\frac{\rm i}2\lambda}+\sum_{j\geq s}\tilde C_j z^s,\quad \tilde C_s\neq 0.
\eeq

For $r>0$, let $\Delta_2(0,r):=\{(\xi,\eta)\in\C^2: |\xi|,|\eta|<r\}$
and for $\omega\in \R$, let ${\cL C}_{\omega}^r:=\left\{(\xi,\eta)\in \C^2: \xi\eta =\omega, \  |\xi|, |\eta|<r \right\}$. Obviously, ${\cL C}_{\omega}^r$ is empty if $|\omega|\geq r^2$.
The following theorem shows that there is a family of invariant closed curves for the involutions $\tau_j^o$ and the reversible map $\sigma_o$ in any neighborhood of the origin.

\begin{thm}\label{thm_inv}
With the notations above and under assumption \re{nd3}, there exists a small enough $R=R(\lambda,r,s)>0$ such that
there is a compact set ${\cL O}_{\infty}(R)\subset]- R^2, R^2[$ satisfying
\begin{equation}\label{mesure_asymp}
\frac{|{\cL O}_{\infty}(R)|}{2R^2}\to 1, \footnote{Through the paper, for any ${\cL S}\subset \R$, $|{\cL S}|$ denotes its Lebesgue measure.} \quad R\to 0,
\end{equation}
such that for any $\omega\in{\cL O}_{\infty}(R)$, one can find $\mu_\om\in \R$ and a holomorphic transformation $\Psi_{\om}:{\cL C}_{\omega}^{R}\to \Delta_2(0, R^{\frac12})$ with $\Psi_{\om}\circ\rho=\rho\circ\Psi_{\om}$, such that, on ${\cL C}_{\omega}^{R}$~:
\begin{gather*}
	(\Psi_{\om}^{-1}\circ \tau^o_1\circ\Psi_{\om})(\xi,\eta)=\left(\begin{array}{l}
		e^{\frac{\rm i}{2}\mu_\om}\eta \\
		e^{-\frac{\rm i}{2}\mu_\om}\xi
	\end{array}\right),\quad (\Psi_{\om}^{-1}\circ \tau^o_2\circ\Psi_{\om})(\xi,\eta)=\left(\begin{array}{l}
		e^{-\frac{\rm i}{2}\mu_\om}\eta \\
		e^{\frac{\rm i}{2}\mu_\om}\xi
	\end{array}\right),\\
	(\Psi_{\om}^{-1}\circ\sigma_o\circ\Psi_{\om})(\xi,\eta)=
	\left(
	\begin{array}{l}
		e^{{\rm i}\mu_\om}\xi\\
		e^{-{\rm i}\mu_\om}\eta
	\end{array}
	\right),\quad (\xi,\eta)\in {\cL C}_{\omega}^{R}.
\end{gather*}
In other words, $\tau^o_1,\tau^o_2$ and $\sigma_o$ have $\Psi_{\om}({\cL C}_{\omega}^{R})$ as holomorphic invariant set and their restrictions to it are conjugate to the restrictions to ${\cL C}_{\omega}^{R}$ of linear maps defined above.
Moreover, $\mu_\om\in ]\lambda-\frac{\pi}{4}, \lambda+\frac{\pi}{4}[$ depends on $\om$ smoothly in the sense of Whitney, and $\Psi_{\om}=\check\Psi\circ({\rm Id}+\phi_\omega)$, with $\check\Psi$ biholomorphic on the neighborhood $\Delta_2(0,R)$, fixing the origin, and $\phi_\omega$ is smooth with respect to $\om$ and sufficiently small in the sense of Whitney.
\end{thm}


\begin{remark} If the surface $M$ can be holomorphically flattened, that is, if it can be holomorphically mapped into $\imag(z_2)=0$, then the situation is much simpler. Indeed, in that case, the associated dynamical system has an extra holomorphic first integral \cite{gong-integral}. It implies that automatically, in good holomorphic coordinates near the origin, all curves $\{\xi\eta= constant\}$ are left invariant by the original dynamics. One thus needs to prove that for suitable values of these constants (i.e. $\omega$'s), one has a conjugacy to linear maps on the associated $\{\xi\eta= constant\}$ as mentioned by Sevryuk \cite{sevryuk-ams} in his Mathematical review of Gong's article \cite{gong-integral}.
\end{remark}

\begin{remark} $\mu_\omega$ obtained in Theorem \ref{thm_inv} is such that $\frac{\mu_{\om}}{\pi}$ is irrational and $\{|e^{{\rm i}n\mu_\omega}-1|\}_{n\in\N^*}$ does not accumulate at the origin too quickly (see (\ref{tilde_O})), which guarantees that the restriction $\sigma|_{{\CS}_\omega}$ is an irrational rotation on ${\CS}_\omega$. This contrasts with the example given Section 6(b) of \cite{moser-webster}. Indeed,  the divergence of the formal normalizing transformation $\hat\Psi$ in (\ref{formal_tau1}) -- (\ref{formal_sigma}) cannot be avoided because of the periodic orbits of $\sigma$.
	Such periodic orbits, as well as its invariant curves, do not have any immediate geometrical significance since they do not lie on $M$ but only on its complexification ${\mathcal M}\subset \C^4$.
\end{remark}
\

\noindent
{\it Sketch of proof of Theorem \ref{thm_inv}.} For $\tau^o_1$ and $\tau^o_2$ given as in (\ref{tau1}) and (\ref{tau2}), our aim is to eliminate the perturbation $p^o$ and $q^o$ (hence $\bar p^o$ and $\bar q^o$) by a sequence of holomorphic transformations which commute with $\rho$.

After finitely many steps of normalization in the sense of Poincar\'e-Dulac in the neighborhood of origin, we obtain a pair of involutions of the form
$$\check\tau_1(\xi,\eta)=\left(\begin{array}{l}
		e^{\frac{\rm i}2\check\alpha(\xi\eta)}\eta+\check p(\xi,\eta) \\
		e^{-\frac{\rm i}2\check\alpha(\xi\eta)}\xi+\check q(\xi,\eta)\end{array}\right),\quad \check\tau_2=\rho\circ\check\tau_1\circ\rho,$$
with a non-degenerate $\check\alpha=\check\alpha(z)$ (as in \re{check_alpha}) and higher-order perturbations $\check p$, $\check q$ (as in \re{check_p_q}). Hence, we can make the norm of the perturbation small enough by choosing a small enough neighborhood of origin $\{|\xi|, |\eta|<r_*\}$.
By a possible normalization on the ``crown", $\{|\xi\eta-\omega|<\beta\}$, around $\{\xi\eta=\omega\}$ with $\omega$ well chosen from a compact positive-measure subset of $]-r_*^2,r_*^2[$, the system enters into a general iteration scheme (a KAM-like process, see Proposition \ref{prop_KAM}), under an additional assumption on the perturbation.

By the iteration process, we build a sequence of involutions $\tau_\nu^{(1)}$, $\nu\in\N$, (hence $\tau_\nu^{(2)}=\rho\circ\tau_\nu^{(1)}\circ\rho$ and $\sigma_\nu=\tau_\nu^{(1)}\circ\tau_\nu^{(2)}$) of the form
$$\tau_\nu^{(1)}(\xi,\eta)=\left(\begin{array}{l}
		e^{\frac{\rm i}2\alpha_\nu(\xi\eta)}\eta+p_\nu(\xi,\eta) \\
		e^{-\frac{\rm i}2\alpha_\nu(\xi\eta)}\xi+q_\nu(\xi,\eta)\end{array}\right)$$
on crowns around $\{\xi\eta=\omega\}$, that shrink to the connected holomorphic curve $\{\xi\eta=\omega\}$ when $\nu$ tends to infinite. On the other hand, when restricted to $\{\xi\eta=\omega\}$, $\alpha_\nu$ tends to a real number $\alpha_{\infty}(\omega)$ and the perturbation $(p_\nu, q_\nu)$ tends to zero, as $\nu$ tends to infinite. In order to control this process, one has to exclude some parameters $\omega$ from the previous set and to show that, the set of admissible parameters $\omega$ for full process is non-void.

The required supplementary condition mentioned above is that the ``crossing term", called {\it skew term} below, $e^{\frac{\rm i}2\alpha_\nu(\xi\eta)}\eta q_\nu+e^{-\frac{\rm i}2\alpha_\nu(\xi\eta)}\xi p_\nu$ of $\tau_\nu^{(1)}$ is much smaller than $p_\nu$ and $q_\nu$ (see (\ref{esti_Lpq})) on the crown. With this condition, we are able to construct a suitable holomorphic transformation of the form
$$\psi_\nu(\xi,\eta)=\left(\begin{array}{l}
\xi+ u_\nu(\xi,\eta)\\
\eta+ v_\nu(\xi,\eta)
\end{array}\right),$$
which conjugates $\tau_\nu^{(1)}$ into $\tau_{\nu+1}^{(1)}$ with perturbation of much smaller size on a smaller crown around $\{\xi\eta=\omega\}$. Here $\omega$ is chosen from a suitable real parameter set, which is related to the small-divisor conditions and guarantees the convergence for the product of sequence of transformations $\{\psi_\nu\}$ on $\{\xi\eta=\omega\}$.

Indeed, the supplementary condition on the skew term of $\tau_\nu^{(1)}$ implies the skew term $\eta u_\nu +\xi v_\nu$ of transformation $\psi_\nu$ is much smaller (similar to \re{esti_Luv} in Lemma \ref{lem_cohomo}). As a consequence, the error term coming from the {\it non-degeneracy} of the ``eigenvalues" $e^{\frac{\rm i}2\alpha_\nu(\xi\eta)}$:
$$e^{\pm\frac{\rm i}2\alpha_\nu(\xi\eta+\eta u_\nu+\xi v_\nu+ u_\nu v_\nu)}- e^{\pm\frac{\rm i}2\alpha_\nu(\xi\eta)}$$
is so small that it can be directly put into the new perturbation.
We emphasize that this supplementary hypothesis on the skew terms of $\tau_\nu^{(1)}$ has to be assumed only at the initial KAM step (i.e.,  for $\nu=0$), as  after each $(\nu+1)-{\rm th}$ KAM step, the new skew term of $\tau_{\nu+1}^{(1)}$ is automatically much smaller than $p_{\nu+1}$ and $q_{\nu+1}$.
This is due to a subtle cancellation of main parts (see \re{esti_p_+q_+} and \re{esti_Lp_+q_+} in Theorem \ref{thm_trans_cohomo} and its proof). As mentioned above, by an initial preparation of the involutions, we can make them to satisfy this supplementary condition required in the iteration process.

\subsection{Geometry of hyperbolic CR singularity}
We recall that $M$ is {\it non-exceptional}, since $\frac{\lambda}{\pi}\in\R\setminus\Q$ in the associated involutions $\tau_1^o$ and $\tau_2^o$ given in \re{tau1} and \re{tau2}. Hence, \re{formal_tau1} and \re{formal_tau2} hold. Let us show that Theorem \ref{thm_inv} enables us to obtain the result on the geometry of real analytic surfaces with a hyperbolic CR singularity.

As mentioned above, the triple $\{\tau_1^o,\tau_2^o,\rho\}$ given in \re{tau1}, \re{tau2} completely characterizes the holomorphic equivalent class of the submanifold $M$ given in \re{defM}. Indeed, following Moser-Webster \cite{moser-webster}, we can reconstruct a submanifold from a pair of involutions.
 Let us define two holomorphic mappings $\varphi_1$, $\Phi$ fixing the origin of $\mathbb{C}^2$ as follows~:
 \begin{equation}\label{varphi_1_2}
 \varphi_1(\xi,\eta):=\xi+\xi\circ\tau^o_1, \quad
	\varphi_2:=\ov{\varphi_1\circ \rho},\quad \rho(\xi,\eta)=(\bar\xi,\bar \eta).
 \end{equation}
The latter implies that the biholomorphic mapping, fixing the origin of $\mathbb{C}^2$,
\begin{equation}\label{varphi_}
\varphi(\xi,\eta)= (\varphi_1(\xi,\eta),\varphi_2(\xi,\eta))=:(z',w')
 \end{equation}
	 transforms $\rho$ into the standard
	complex conjugation $(z',w')\to(\ov w',\ov z')$. Define
\begin{equation}\label{Phi}
	\Phi(\xi,\eta):=(\xi\circ\tau^o_1(\xi,\eta))\cdot \xi,
 \end{equation}
where $\xi\circ\tau^o_1(\xi,\eta)$ denotes the $\xi-$coordinate of $\tau^o_1(\xi,\eta)$.
	We verify that $\varphi_1$ and $\Phi$ are invariant by $\tau^o_1$. Then the local analytic submanifold defined by the local equation
	\beq\label{invol2manifld}
	z_2=(\Phi\circ \var^{-1})(z_1,\bar z_1),
	\eeq
	has $\{\tau^o_1, \tau^o_2,\rho\}$  as associated Moser-Webster involutions. Assume that the pair of involutions $\tau_1^o$ and $\tau_2^o$ associated to the real surface $M$ satisfies the assumption of Theorem \ref{thm_inv}.
	 Then, for each $\om\in \cL{O}_{\infty}(R)$, there is a holomorphic map $\Psi_{\om}$ defined on ${\cL C}_{\omega}^R=\left\{(\xi,\eta)\in \C^2: \xi\eta =\omega, \  |\xi|, |\eta|<R \right\}$. It is a small perturbation of the the identity. Let $H_{\om}^R$ be the real hyperbola of ${\cL C}_{\omega}^R$, i.e., ${H}_{\omega}^R:=\left\{(\xi,\eta)\in \R^2: \xi\eta =\omega, \  |\xi|, |\eta|<R \right\}$. Let us consider the connected holomorphic curve defined as the image of ${\cL C}_{\omega}^R$ (which contains a graph over an annulus)~: $$
	 \mathcal{S}_{\omega}: \left\{\begin{array}{ccc}z_1&=&\varphi_1\circ \Psi_{\om}(\xi,\eta)\\
	 	z_2&=&\Phi\circ \Psi_{\om}(\xi,\eta)\end{array}\right.  ,\quad (\xi,\eta)\in \mathcal{C}_{\omega}^{R}.
	 $$ as well as $\mathcal{H}_{\om}^R$, the image of the restriction, $\tilde\Psi_{\om}:=(\varphi_1,\Phi)\circ \Psi_{\om}|_{{H}_{\om}^R}$, to ${H}_{\om}^R$. Hence, $\mathcal{H}_{\om}^R\subset \mathcal{S}_{\omega}$ shrinks to zero with $R$.

\begin{thm}
Under the assumption of Theorem \ref{thm_inv} and the notation above, the family $\{\mathcal{S}_{\om}\}_{\om\in \cL{O}_{\infty}(R)}$ is a non-constant Whitney smooth family of connected holomorphic curves. Each of them intersects M, in a neighborhood of the origin, along the holomorphic hyperbola $\mathcal{H}^R_{\om}$.
\end{thm}
\begin{remark}
Assumptions of the previous theorem, through \re{nd3}, implies that the real analytic surface $M$ given in (\ref{defM}) is not formally equivalent to ${\cL Q}_\gamma$.
\end{remark}

\begin{remark} The conclusion of the previous theorem contrasts with that of the elliptic case treated by Moser-Webster. Indeed, in the holomorphic normalizing coordinates, there is a real analytic family of holomorphic curves $\CS_c: z_2=c$ for $c$ in a real neighborhood of the origin, and for every $c$, $\CS_c$ intersects $M$ along the ellipse $c=|z_1|^2+(\gamma+\epsilon c^s)(z_1^2+\bar z_1^2)$.
\end{remark}

\begin{figure*}[hbtp]
    \begin{center}
        \leavevmode
        \includegraphics[scale=.5]{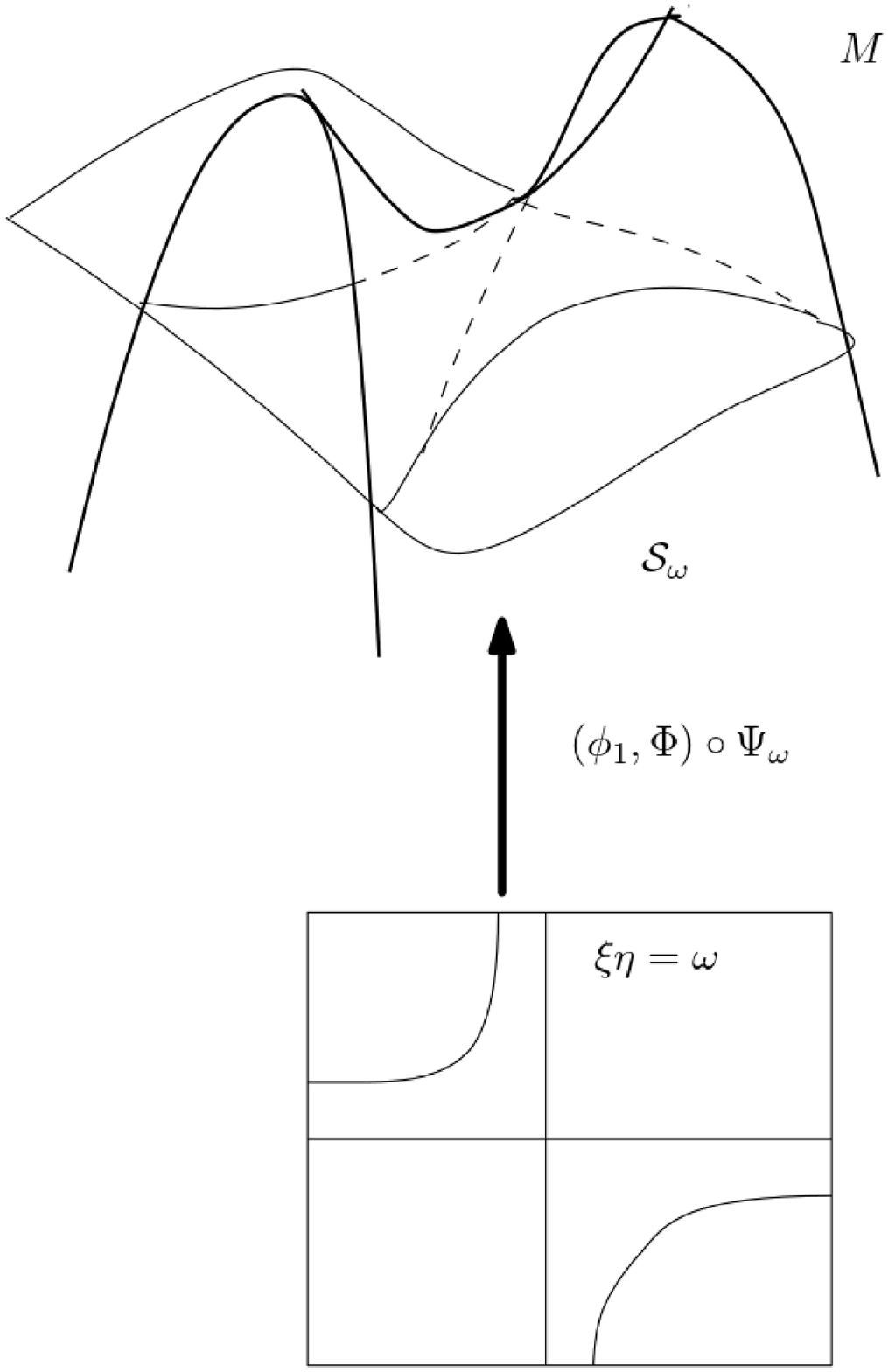}
        \caption{Holomorphic hyperbola~: Intersection of $M$ by a holomorphic curve.}
    \end{center}
\end{figure*}


\proof According to \rt{thm_inv}, for any good parameter $\om\in {\cL O}_{\infty}(R)$, there exists a connected holomorphic curve $\cL S_{\om}$ invariant by the dynamics and $\rho$, such that $\tau^o_{j}|_{\cL S_{\om}}$ is conjugated to the restriction to ${\cL C}_{\omega}^{R}$ of the linear involutions~:
$$
T_j:(\xi,\eta)\mapsto\left(e^{\frac{\rm i}2(-1)^{j-1}\mu_{\om}}\eta, \,  e^{-\frac{\rm i}2(-1)^{j-1}\mu_{\om}}\xi
\right).
$$
Indeed, with the definitions given in (\ref{varphi_1_2}) and (\ref{Phi}), we have, for all $(\xi,\eta)\in {\cL C}_{\omega}^{R}$,
\begin{align*}
(\varphi_1\circ \Psi_{\om})(\xi,\eta)&= (\xi\circ \Psi_{\om})(\xi,\eta)+\left(\xi\circ \Psi_{\om}\circ \left(\Psi_{\om}^{-1}\circ\tau^o_1\circ \Psi_{\om}\right)\right)(\xi,\eta) \\
&=(\xi\circ \Psi_{\om})(\xi,\eta)+(\xi\circ \Psi_{\om})\left(e^{\frac{\rm i}2\mu_{\om}}\eta,e^{-\frac{\rm i}2\mu_{\om}}\xi\right),\\
(\Phi\circ \Psi_{\om})(\xi,\eta)&=(\xi\circ \Psi_{\om})(\xi,\eta)\cdot \left((\xi\circ \Psi_{\om})\left(e^{\frac{\rm i}2\mu_{\om}}\eta,e^{-\frac{\rm i}2\mu_{\om}}\xi\right)\right).
\end{align*}
We define the connected holomorphic curve $\mathcal{S}_{\omega}$ as the image of the holomorphic curve $\mathcal{C}_{\omega}^{R}$ by the holomorphic map~:
$$
\mathcal{S}_{\omega}: \left\{\begin{array}{ccc}z_1&=&(\varphi_1\circ \Psi_{\om})(\xi,\eta)\\
z_2&=&(\Phi\circ \Psi_{\om})(\xi,\eta)\end{array}\right.  ,\quad (\xi,\eta)\in \mathcal{C}_{\omega}^{R}.
$$

On the one hand, since $\Psi_{\om}$ and $\rho$ commutes, we have
$$
\var\circ \Psi_{\om}=\left(\varphi_1\circ \Psi_{\om}, \overline{\varphi_1\circ \Psi_{\om}\circ\rho}\right).
$$
On the other hand, ${H}^R_{\omega}:={\cL C}_{\omega}^{R}\cap{\rm Fix}(\rho)$ is the union of two branches of real hyperbola $\{\xi\eta=\omega, \ \xi=\bar\xi, \ \eta=\bar \eta, \ |\xi|,|\eta|<R \}$. Hence, we have
$$
\var\circ \Psi_{\om}|_{{H}^R_{\omega}}=(z_1,\bar z_1)|_{{H}^R_{\omega}}.
$$
As a consequence, the complex curve $\mathcal{S}_{\om}$ intersects $M$ given in \re{invol2manifld} along the image of ${H}^R_{\om}$ by $(\varphi_1,\Phi)\circ\Psi_{\om}$,
$$
z_2=(\Phi\circ \Psi_{\om})\circ \left(\Psi_{\om}^{-1}\circ\var^{-1}\right)(z_1,\bar z_1)|_{{H}^R_{\omega}}=(\Phi\circ \Psi_{\om})\circ (\var\circ \Psi_{\om})^{-1}(z_1,\bar z_1)|_{{H}^R_{\omega}}.
$$
We recall that $\Psi_{\om}=\check\Psi\circ({\rm Id} +\phi_{\om}) $.
We can assume that $\check\Psi={\rm Id}$ for convenience. Then, with $\Lambda_\om:=e^{\frac{\rm i}2\mu_{\om}}$, $\CS_{\om}$ is defined as
$$
 \left\{\begin{array}{ccc}z_1&=&\xi+\Lambda_{\om}\eta+(\xi\circ\phi_{\om})(\xi,\eta)+(\xi\circ\phi_{\om})(\Lambda_{\om}\eta,\Lambda_{\om}^{-1}\xi)\\
	z_2&=&\Lambda_{\om}\om+ \xi\cdot((\xi\circ\phi_{\om})(\Lambda_{\om}\eta,\Lambda_{\om}^{-1}\xi))+ \Lambda_{\om}\eta\cdot((\xi\circ\phi_{\om})(\xi,\eta))\\
	& & +((\xi\circ\phi_{\om})(\xi,\eta))\cdot((\xi\circ\phi_{\om})(\Lambda_{\om}\eta,\Lambda_{\om}^{-1}\xi)) \end{array}\right.,
$$
where $(\xi,\eta)\in {\cL C}_{\omega}^{R}$, that is $\eta=\frac{\om}{\xi}$, $\frac{|\om|}{R}<|\xi|<R$. Since $\phi_{\om}$ is sufficiently small and smooth with respect to $\om$ in the sense of Whitney, $\CS_{\om}$ is a small perturbation of the main part
$$
\left\{\begin{array}{ccc}z_1&=&\xi+\Lambda_{\om}\frac{\om}{\xi}\\
	z_2&=&\Lambda_{\om}\om \end{array}\right.,
$$ which varies with $\om$ as $\Lambda_{\om}\om$ does. Indeed, assume that there exist $\omega,\omega'\in{\cL O}_{\infty}(R)\subset\R$ with $\omega\neq \omega'$ such that $\Lambda_{\om}\omega= \Lambda_{\om'}\omega'$. Since $\Lambda_\om=e^{\frac{\rm i}2\mu_{\om}}$ and $\Lambda_{\om'}=e^{\frac{\rm i}2\mu_{\om'}}$, we have
	$$e^{\frac{\rm i}{2}(\mu_{\om'}-\mu_\om)}=\frac{\om}{\om'}\in\R,$$
	which implies that $\mu_{\om'}-\mu_\om = (4k+2)\pi$ for some $k\in\Z$. This contradicts with the fact that, $\mu_\om\in ]\lambda-\frac{\pi}{4}, \lambda+\frac{\pi}{4}[$ for every $\omega\in{\cL O}_{\infty}(R)$.
Hence, $\CS_{\om}$ varies with $\om$.

The Whitney smoothness of $\CS_\om$ follows immediately from that of $\Psi_\om$ and $\mu_\om$.
\qed

\vspace{.5cm}

The rest of paper will be organized as follows. In Section \ref{sec_pre}, the precise definition of crowns around the curve $\{\xi\eta=constant\}\subset\C^2$ and the norm of holomorphic functions on them are introduced, and basic properties associated with reversible map are given.
In Section \ref{sec_KAM}, we give an abstract KAM-like theorem, which is used to prove Theorem \ref{thm_inv} in Section \ref{sec_proof}. A preliminary normalization (which is required to start the KAM-like process, also known as Newton method), as well as the Whitney smoothness of the family of invariant curves, is also given in Section \ref{sec_proof}. In Section \ref{sec_rev}, we describe properties of the pair of holomorphic involutions $\{\tau_1,\tau_2=\rho\circ\tau_1\circ\rho\}$ and in particular of their non-degenerate principal parts, as well as of their reversible associated composition $\sigma=\tau_1\circ\tau_2$.
In Section \ref{sec_trans}, two types of holomorphic transformations of $\{\tau_1,\tau_2\}$, commuting with $\rho$, are introduced, which is used to complete the proof of the KAM-like theorem.

\section{Preliminaries and notations}\label{sec_pre}

\subsection{Basic property of reversible map}

Let us define the involution $\rho(\xi,\eta)=(\bar\xi,\bar\eta)$. An invertible map $\sigma:\C^2\to \C^2$ is called {\it reversible} with respect to $\rho$ if $\sigma^{-1}=\rho\circ\sigma\circ\rho$.
\begin{lemma}\label{lem_reve_psi}
	Given $\psi:\C^2\to\C^2$ with $\psi=(u,v)$. 
	Then $\psi\circ\rho=\rho\circ\psi$ if and only if
	$u=\bar u$, $v=\bar v$.
\end{lemma}
\proof Since $\rho\circ\rho={\rm Id}$, we have $\rho^{-1}\circ\psi\circ\rho=\rho\circ\psi\circ\rho=(\bar u,\bar v)$ which equals to $\psi$ if and only if $u=\bar u$, $v=\bar v$.\qed

\begin{lemma}\label{lem_reve}
Let $\sigma:\C^2\to \C^2$ be a reversible map w.r.t. $\rho$, and let $\psi:\C^2\to \C^2$ be an invertible map commuting with $\rho$, 
Then $\psi^{-1}\circ\sigma\circ\psi$ is also reversible w.r.t. $\rho$.
\end{lemma}
\proof The conclusion follows from
\begin{eqnarray*}
(\psi^{-1}\circ\sigma\circ\psi)^{-1}&=&\psi^{-1}\circ\sigma^{-1}\circ\psi\\
&=&(\rho\circ \psi^{-1}\circ\rho)\circ (\rho\circ\sigma\circ \rho) \circ(\rho\circ \psi\circ\rho) \, = \, \rho\circ(\psi^{-1}\circ\sigma\circ\psi)\circ\rho.\qed
\end{eqnarray*}

\subsection{Function space and norms}

Given $0<r<\frac14$ and $0\leq \beta< r^2$, for $\omega\in ]-r^2+\beta,r^2-\beta[$, we define
\begin{align}
{\cL C}^r_{\omega}:=&\left\{(\xi,\eta)\in \C^2:  \xi\eta =\omega, \ |\xi|,|\eta|<r \right\},\quad {\cL C}_{\omega}:=\bigcup_{r> 0}{\cL C}^r_{\omega},\label{cw}\\
{\cL C}^r_{\omega,\beta}:=&\left\{(\xi,\eta)\in \C^2: |\xi\eta -\omega|\leq\beta, \ |\xi|,|\eta|<r \right\},\quad {\cL C}_{\omega,\beta}:=\bigcup_{r> 0}{\cL C}^r_{\omega,\beta}.\label{cwb}
\end{align}

For a power series
$$
f(\xi,\eta)=\sum_{l,j\geq 0}\breve f_{l,j}\xi^l\eta^j,\quad \breve f_{l,j}\in\C,
$$
we have the unique decomposition
\begin{equation}\label{decomp}
f(\xi,\eta)=f_{0,0}(\xi\eta)+\sum_{l\geq 1} f_{l,0}(\xi\eta)\,  \xi^l+\sum_{j\geq 1} f_{0,j}(\xi\eta) \, \eta^j=\sum_{l,j\geq 0\atop{l\cdot j=0}}f_{l,j}(\xi\eta)\xi^l\eta^j,
\end{equation}
with the {\it coefficients} of $f$, depending on the product $\xi\eta$, given by
$$
f_{l,j}(\xi\eta)   = \sum_{k\geq 0}\breve f_{k+l,k+j}\cdot(\xi\eta)^k,\quad l\cdot j=0.
$$
Sometimes, by defining $f_{l,j}=0$ for $lj\neq 0$, we rewrite $f$ as
$$
f(\xi,\eta)=\sum_{l,j\geq0}f_{l,j}(\xi\eta)\,  \xi^l\eta^j.
$$
Let us consider the anti-holomorphic involution $\rho:(\xi,\eta)\mapsto (\bar\xi,\bar\eta)$. We define the conjugate of $f$ to be $\bar f$, whose Taylor expansion coefficients at the origin are the complex conjugates of those of $f$.
Obviously, $f=\bar f$ if and only if $\breve f_{l,j}\in\R$ for all $l,j\geq 0$.

Let $h=h(\xi\eta)$ be a function of the product $\xi\eta$. We define
$$
|h|_{\omega,\beta}:=\sup_{(\xi,\eta)\in {\cL C}_{\omega,\beta}}|h(\xi\eta)|= \sup_{|z-\omega|<\beta}|h(z)|,\quad |h|_{\omega}:=\sup_{(\xi,\eta)\in  {\cL C}_{\omega}}|h(\xi\eta)|=\sup_{z=\omega}|h(z)|.
$$
Given a power series $f$, we define the norms
\beq\label{norm-0}
|f|_r:=\sum_{l,j\geq 0}|\breve f_{l,j} \, |r^{l+j},\quad \|f\|_{\omega,\beta,r}:=\sum_{l,j\geq0\atop{lj=0}}|f_{l,j}|_{\omega,\beta}\, r^{l+j}.
\eeq
In particular, for a function $f$ of $\xi\eta$, we have $\|f\|_{\omega,\beta,r}=|f|_{\omega,\beta}$.
For $\omega\in ]-r^2+\beta,r^2-\beta[$, it is easy to see that
\begin{equation}\label{comp_norm}
\sup_{(\xi,\eta)\in{\cL C}_{\omega,\beta}^r} |f(\xi,\eta)|\leq \|f\|_{\omega,\beta,r}\leq |f|_r.
\end{equation}
The definition of norm $\|\cdot\|_{\omega,\beta,r}$ implies that
\begin{equation}\label{esi_coeff}
\|f(\xi,\eta)\|_{\omega,\beta,r}=\|f(\eta,\xi)\|_{\omega,\beta,r}=\|\bar f(\xi,\eta)\|_{\omega,\beta,r},\quad|f_{l,j}|_{\omega,\beta}\leq  \|f\|_{\omega,\beta,r} r^{-(l+j)}, \;\ l,j\geq 0.
\end{equation}
For ${\CO}\subset ]-r^2+\beta,r^2-\beta[$, we also define the norm
$$\|f\|_{{\CO},\beta,r}:=\sup_{\omega\in{\CO} \atop{|\omega|<r^2-\beta}}\|f\|_{\omega,\beta,r}.
$$
In particular, for the coefficients of $f$,
$$ \|f_{l,j}\|_{{\CO},\beta,r}=\sup_{\omega\in{\CO} \atop{|\omega|<r^2-\beta}} |f_{l,j}|_{\omega,\beta},\quad l,j\geq 0, \ lj=0.$$
\begin{lemma}\label{lem_norm00} For given power series $f$ and $g$, if $|\omega|< r^2-\beta$, then
$\|fg\|_{\omega,\beta,r}\leq\|f\|_{\omega,\beta,r}\|g\|_{\omega,\beta,r}$.
\end{lemma}
\proof Decompose $fg$ as in (\ref{decomp}), we have
$$(fg)(\xi,\eta)=(fg)_{0,0}+\sum_{l\geq 1}(fg)_{l,0}\xi^l+\sum_{j\geq 1}(fg)_{0,j}\eta^j$$
with the coefficients given by
\begin{eqnarray}
(fg)_{0,0}&=&f_{0,0}  g_{0,0}+\sum_{k\geq 1}(f_{k,0}g_{0,k}+f_{0,k}g_{k,0})(\xi\eta)^k,\label{fg00}\\
(fg)_{l,0}&=&\sum_{k=0}^{l-1}(f_{l-k,0}g_{k,0}+f_{k,0}g_{l-k,0})+ \sum_{k\geq1}(f_{l+k,0}g_{0,k}+f_{0,k}g_{l+k,0})(\xi\eta)^k,\label{fgl0}
\end{eqnarray}
and similar expression for $(fg)_{0,j}$. If $|\omega|< r^2-\beta$ and $(\xi,\eta)\in {\cL C}_{\omega,\beta}$, then $|\xi\eta|<r^2$. Therefore, we have
\begin{eqnarray*}
|(fg)_{0,0}|_{\omega,\beta}&\leq&|f_{0,0}|_{\omega,\beta}|g_{0,0}|_{\omega,\beta}+\sum_{k\geq 1}\left(|f_{k,0}|_{\omega,\beta}|g_{0,k}|_{\omega,\beta}+|f_{0,k}|_{\omega,\beta}|g_{k,0}|_{\omega,\beta}\right)r^{2k},\\
|(fg)_{l,0}|_{\omega,\beta}r^l&\leq&\sum_{k=0}^{l-1}\left(|f_{l-k,0}|_{\omega,\beta}|g_{k,0}|_{\omega,\beta}+|f_{k,0}|_{\omega,\beta}|g_{l-k,0}|_{\omega,\beta}\right)r^l\\
& &+ \sum_{k\geq1}\left(|f_{l+k,0}|_{\omega,\beta} |g_{0,k}|_{\omega,\beta}+|f_{0,k}|_{\omega,\beta}|g_{l+k,0}|_{\omega,\beta}\right)r^{2k+l},
\end{eqnarray*}
and a similar estimate for $|(fg)_{0,j}|_{\omega,\beta}r^j$.
Since $\|fg\|_{\omega,\beta,r}=\sum_{l,j\geq 0}|(fg)_{l,j}|_{\omega,\beta}r^{l+j}$ with $(fg)_{l,j}=0$ whenever $lj\neq 0$,
we obtain
$$
\|fg\|_{\omega,\beta,r}\leq\left(\sum_{l,j\geq 0} |f_{l,j}|_{\omega,\beta}r^{l+j}\right)\left(\sum_{l,j\geq 0} |g_{l,j}|_{\omega,\beta}r^{l+j}\right)=\|f\|_{\omega,\beta,r}\|g\|_{\omega,\beta,r}.\qed
$$

Let $0\leq\beta<r^2$ and ${\CO}\subset ]-r^2, \, r^2 [$. We define $\CH_{\beta,r}(\CU)$ to be the set of holomorphic functions in a complex neighborhood $\CU$ of
\beq\label{Orbeta}{\CO}(r,\beta):={\CO} \, \cap \, ]-r^2+\beta,r^2-\beta[,\eeq
and the collections of power series
\begin{eqnarray*}
{\CA}_{\beta, r}({\CO})&:=&\left\{f=\sum_{l,j}f_{l,j}(\xi\eta)\xi^l\eta^j: \begin{array}{l}
f \ {\rm is \ holomorphic \ on} \ \bigcup_{\omega\in{\CO}(r,\beta)} {\cL C}_{\omega,\beta}^r,\\
f_{l,j}\in  \CH_{\beta,r}(\CU),\; \CU {\rm \ complex \ neighborhood \ of\ } {\CO}(r,\beta) 
\end{array}  \right\}, \\
{\CA}^\R_{\beta, r}({\CO})&:=&\left\{f\in {\CA}_{\beta, r}({\CO}): f=\bar f\right\}.
\end{eqnarray*}
In the above definition, $\CU$ denotes an unprecised neighborhood over which all $f_{l,j}$'s are holomorphic.
Finally for any $\tilde\CO\subset \R$ and any function $h$ defined on $\tilde \CO$, we define the norm
$|h|_{\tilde\CO}:=\sup_{\omega\in\tilde\CO}|h(\omega)|$.
By (\ref{esi_coeff}), we see that
\begin{equation}\label{esi_coeff_uniform}
 |f_{l,j}|_{{\CO}(r,\beta)} = \|f_{l,j}\|_{{\CO},\beta,r}\leq \|f\|_{{\CO},\beta,r}r^{-(l+j)}.
\end{equation}
It is easy to see that:
\begin{itemize}
\item {\it (linear structure)} for $f,g\in {\CA}_{\beta, r}(\CO)$ (or ${\CA}^\R_{\beta, r}(\CO)$), we have $a_1f+a_2g\in {\CA}_{\beta, r}(\CO)$ (or ${\CA}^\R_{\beta, r}(\CO)$), $a_1,a_2\in\C$(or $\R$), with
\begin{equation}\label{esti_linear}
\|a_1f+a_2g\|_{{\CO},\beta,r}\leq |a_1| \|f\|_{{\CO},\beta,r}+|a_2|\|g\|_{{\CO},\beta,r}.
\end{equation}
\item {\it (monotonicity)} If $\CO'\subset \CO$ and $r'\leq r$, $\beta'\leq \beta$, $r'^2-\beta'\leq r^2-\beta$, then ${\CA}_{\beta, r}(\CO)\subset {\CA}_{\beta', r'}(\CO')$ and ${\CA}^\R_{\beta, r}(\CO)\subset {\CA}^\R_{\beta', r'}(\CO')$ with
\begin{equation}\label{esti_mono}
\|f\|_{\CO',\beta', r'}\leq \|f\|_{\CO,\beta, r},\quad \forall \ f\in {\CA}_{\beta, r}(\CO) \ {\rm or } \  {\CA}^\R_{\beta, r}(\CO).
\end{equation}
\end{itemize}

%
%
%
%
%
%

\begin{lemma}\label{lemma_analytic_prod}
For $f,g$ with $\|f\|_{{\CO},\beta,r}$, $\|g\|_{{\CO},\beta,r}<\infty$,
we have that
\begin{equation}\label{prod_O}
\|fg\|_{{\CO},\beta,r}\leq \|f\|_{{\CO},\beta,r}\|g\|_{{\CO},\beta,r}.
\end{equation}
Moreover, if $f,g\in {\CA}_{\beta,r}(\CO)$, then $fg\in {\CA}_{\beta,r}(\CO)$.
\end{lemma}
\proof In view of Lemma \ref{lem_norm00}, we obtain \re{prod_O}. Provided that $f,g\in {\CA}_{\beta,r}(\CO)$, let us prove the analyticity of coefficients $(fg)_{l,j}(\cdot)$, $l,j\geq 0$ with $lj=0$, on a same neighborhood of ${\CO}(r,\beta)$. By (\ref{fgl0}), we have, for $l\geq 1$,
\begin{eqnarray*}
(fg)_{l,0}(\omega)&=&\sum_{k=0}^{l-1}\left(f_{l-k,0}(\omega)g_{k,0}(\omega)+f_{k,0}(\omega)g_{l-k,0}(\omega)\right)\\
& & + \,  \sum_{k\geq1}\left(f_{l+k,0}(\omega)g_{0,k}(\omega)+f_{0,k}(\omega)g_{l+k,0}(\omega)\right)\omega^k.
\end{eqnarray*}
Let $\CU$ be a complex neighborhood of ${\CO}(r,\beta)$ over which the $f_{k,j}, g_{k,j}$'s are holomorphic. Since $|\omega|<r^2-\beta$ and, in view of (\ref{esi_coeff_uniform}), we have
$$
|f_{l+k,0}|_{{\CO}(r,\beta)}|g_{0,k}|_{{\CO}(r,\beta)}+|f_{0,k}|_{{\CO}(r,\beta)}|g_{l+k,0}|_{{\CO}(r,\beta)} \leq\frac{2 \|f\|_{{\CO}, \beta, r}\|g\|_{{\CO,\beta, r}}}{r^{l+2k}}.
$$
The analyticity of the $(fg)_{l,0}$'s on a same neighborhood of ${\CO}(r,\beta)$ follows.
The proof for $(fg)_{0,j}$, $j\geq 0$, is similar by the proof of \rl{lem_norm00}.
\qed

\begin{lemma}\label{lemma_analytic_exp} Given $f\in {\CA}_{\beta, r}(\CO)$ with $\|f\|_{{\CO}, \beta, r}<\infty$ and $a\in\C$, we have $e^{a f}\in {\CA}_{\beta, r}(\CO)$.
\end{lemma}
\proof By Lemma \ref{lemma_analytic_prod}, we see that $f^k\in {\CA}_{\beta, r}(\CO)$ with $\|f^k\|_{{\CO}, \beta, r}\leq \|f\|^k_{{\CO}, \beta, r}$ for every $k\in\N$.
Then, according to (\ref{esi_coeff_uniform}), we have that, for $l,j\geq 0$ with $lj=0$,
$$|(a^kf^k)_{l,j}|_{{\CO}(r,\beta)}\leq \frac{a^k \|f^k\|_{{\CO}, \beta, r}}{r^{l+j}} \leq\frac{a^k }{r^{l+j}}\|f\|^k_{{\CO}, \beta, r}.$$
Developing the exponential function around $0$, we have, for $\omega\in {\CO}(r,\beta)$,
$$(e^{a f})_{l,j}(\omega)=1+\sum_{k\geq 1} \frac{a^k(f^k)_{l,j}(\omega)}{k!}.$$
Hence we obtain the analyticity of $(e^{a f})_{l,j}$, since
$$\left|\sum_{k\geq 1} \frac{(a^kf^k)_{l,j}(\omega)}{k!}\right|\leq \frac{1}{r^{l+j}}\sum_{k\geq 1}\frac{a^k \|f\|^k_{{\CO}, \beta, r}}{k!}.\qed$$

\begin{lemma}\label{lem_norm-pre}
Given $0<r''<r'<\frac{1}{4}$, $\CO\subset]-r'^2,r'^2[$ and $0<2\beta''\leq\beta'$, if $\beta'$ is sufficiently small such that
\begin{equation}\label{smallness_beta}
e^{\frac98\beta'}\frac{r''}{r'}<1-\frac{\beta'^2}{16},\quad 8\beta'^{\frac12}<(r'-r'')r'',
\end{equation}
then for
$h \in {\CA}_{\beta',r'}\left({\CO}\right)$ with $\|h\|_{\CO,\beta',r'}< +\infty$, for $f_1,f_2,g_1,g_2$ satisfying
$$
\|f_m\|_{\CO,\beta'',r''},  \|g_m\|_{\CO,\beta'',r''} < \frac{\beta'^2}{16} , \quad m=1,2,
$$
we have that
$$\|h(\xi+f_1,\eta+g_1)-h(\xi+f_2,\eta+g_2)\|_{{\CO},\beta'',r''}\nonumber\\
<\frac{3r'\|h\|_{\CO,\beta',r'}}{(r'-r'')\beta'}\max\left\{ \|f_1-f_2\|_{\CO,\beta'',r''}, \|g_1-g_2\|_{\CO,\beta'',r''}\right\}.$$
Moreover, if $f_1,f_2,g_1,g_2\in {\CA}_{\beta'',r''}\left({\CO}\right)$, then
$$h(\xi+f_1,\eta+g_1)-h(\xi+f_2,\eta+g_2)
 \in   {\CA}_{\beta'',r''}\left({\CO}\right).$$
\end{lemma}
\begin{remark} Note that the second inequality in (\ref{smallness_beta}) implies that $r''^2-\beta''<r'^2-\beta'$. Hence, by the monotonicity, ${\CA}_{\beta',r'}\left({\CO}\right)\subset {\CA}_{\beta'',r''}\left({\CO}\right)$.
\end{remark}
A more general version of Lemma \ref{lem_norm-pre} will be given in Section \ref{sec_rev} (see Lemma \ref{lem_norm}), and will be shown in Appendix \ref{proof_lemma_norm}.

\smallskip

Given $(f,g)\in ({\CA}_{\beta, r}(\CO))^2$, we define, for $\omega\in {\CO}(r,\beta)$,
$$
\|(f,g)\|_{\omega,\beta,r}:= \|f\|_{\omega,\beta,r}+\|g\|_{\omega,\beta,r},\quad \|(f,g)\|_{\CO,\beta,r}:= \|f\|_{\CO,\beta,r}+\|g\|_{\CO,\beta,r}.$$

\begin{lemma}\label{lem_tildeU}
Given $0<r''<r'<\frac{1}{4}$, $\CO\subset ]-r'^2,r'^2 [$ and $\beta'>0$, consider the transformation $\phi={\rm Id} +{\cL U}$ on ${\cL C}_{\omega,\beta'}^{r'}$ with $\CU\in ({\CA}_{\beta', r'}(\CO))^2$. If $\beta'$ is sufficiently small such that (\ref{smallness_beta}) is satisfied, and
\begin{equation}\label{small_uv}
\|\CU\|_{\CO, \beta', r'}<\frac{\beta'(r'-r'')}{30r'},
\end{equation}
then $\phi$ is invertible on ${\cL C}_{\omega,\beta'}^{r'}$, with $\phi^{-1}-{\rm Id}\in (\CA_{\frac{\beta'}2,r''}(\CO))^2$ and
$$\|\phi^{-1}-{\rm Id}+\CU\|_{\CO,\frac{\beta'}2,r''}\leq \frac{8r'\|{\cL U}\|^2_{\CO,\beta', r'}}{(r'-r'')\beta'}.$$
\end{lemma}
\proof In view of (\ref{small_uv}), we see that $\phi$ is close to identity, hence it is biholomorphic on ${\cL C}_{\omega,\beta'}^{r'}$.
Let us write $\phi^{-1}=:{\rm Id}+{\cL V}$. The identity $\phi\circ\phi^{-1} ={\rm Id}$ means that
\beq\label{formual-inv}
{\cL V}=- \,  {\cL U}\circ({\rm Id}+{\cL V})= - \,  {\cL U} -({\cL U}\circ({\rm Id}+{\cL V})-{\cL U}).
\eeq
By Lemma \ref{lem_norm-pre} (since (\ref{smallness_beta}) is satisfied) and \re{formual-inv}, we have
$$\|\cL V\|_{\CO,\frac{\beta'}{2}, r''}  <  \|{\cL U}\|_{\CO,\beta', r'}+\frac{6r'}{(r'-r'')\cdot \beta'}\|\cL V\|_{\CO,\frac{\beta'}{2}, r''}\|{\cL U}\|_{\CO,\beta', r'}
 < \|{\cL U}\|_{\CO,\beta', r'}+\frac{\|\cL V\|_{\CO,\frac{\beta'}{2}, r''}}{5},$$
which implies that
$
\|\cL V\|_{\CO,\frac{\beta'}{2}, r''} < \frac54\|{\cL U}\|_{\CO,\beta', r'}
$.
Let us set
$\CU_1=- \, \CU$, $\CV_1={\cL U}\circ({\rm Id}+{\cL V})-   {\cL U}$.
Hence, $\CV=\CU_1-  \CV_1$, and, by Lemma \ref{lem_norm-pre},  $\CU_1\in \CA_{\frac{\beta'}2, r''}(\CO)$,
$$\| \CV_1\|_{\CO,\frac{\beta'}2, r''}<\frac{6r'}{(r'-r'')\cdot \beta'}\|\cL V\|_{\CO,\frac{\beta'}{2}, r''}\|{\cL U}\|_{\CO,\beta', r'}
<\frac{\|\cL V\|_{\CO,\frac{\beta'}{2}, r''}}{5}<\frac{\|{\cL U}\|_{\CO,\beta', r'}}4.$$
Then, we have
\begin{eqnarray*}
\CV&=& - \,  {\cL U}\circ({\rm Id}+\CU_1-  \CV_1)\\
&=&- \, {\cL U}\circ({\rm Id}+\CU_1)-\left({\cL U}\circ({\rm Id}+\CU_1-  \CV_1)-{\cL U}\circ({\rm Id}+\CU_1)\right)  \, =: \, \CU_2-\CV_2,
\end{eqnarray*}
and, by Lemma \ref{lem_norm-pre}, $\CU_2\in\CA_{\frac{\beta'}{2}, r''}(\CO)$,
$$\| \CV_2\|_{\CO,\frac{\beta'}2, r''}<\frac{6r'}{(r'-r'')\cdot \beta'}\|\cL V_1\|_{\CO,\frac{\beta'}{2}, r''}\|{\cL U}\|_{\CO,\beta', r'}<\frac{\|\cL V_1\|_{\CO,\frac{\beta'}{2}, r''}}{4}<\frac{\|{\cL U}\|_{\CO,\beta', r'}}{4^2}.$$
Assume that, for some $n\in\N$, we have $\CV=\CU_{n}-\CV_{n}$ with
$$\CU_{n}\in \CA_{\frac{\beta}2, r''}(\CO),\quad \| \CV_n\|_{\CO,\frac{\beta'}2, r''}<4^{-n}\|{\cL U}\|_{\CO,\beta', r'}.$$
Then we have $\CV=\CU_{n+1}-\CV_{n+1}$ with
$$
\CU_{n+1}:= - \, \CU\circ({\rm Id}+\CU_{n}),\quad \CV_{n+1}:=\CU\circ({\rm Id}+\CU_{n}-\CV_{n})-\CU\circ({\rm Id}+\CU_{n}).
$$
By Lemma \ref{lem_norm-pre}, we have $\CU_{n+1}\in \CA_{\frac{\beta}2, r''}(\CO)$ and
$$
\|\CV_{n+1}\|_{\CO,\frac{\beta'}2, r''}<\frac{6r'}{(r'-r'')\cdot \beta'}\|\cL V_n\|_{\CO,\frac{\beta'}{2}, r''}\|{\cL U}\|_{\CO,\beta', r'}<\frac{\|\cL V_n\|_{\CO,\frac{\beta'}{2}, r''}}{4}<\frac{\|{\cL U}\|_{\CO,\beta', r'}}{4^{n+1}}.
$$
As $n\to\infty$, we see that $\CV\in \CA_{\frac{\beta}2, r''}(\CO)$.

Let $\tilde{\cL U}:=\phi^{-1}-{\rm Id} +{\cL U}$.
By \re{formual-inv}, we have
$\tilde{\cL U}={\cL U}-{\cL U}\circ({\rm Id}+\cL V)$. Hence, by Lemma \ref{lem_norm-pre},
$$
\|\tilde{\cL U}\|_{\CO,\frac{\beta'}{2}, r''}=\|{\cL U}-{\cL U}\circ({\rm Id}+\cL V)\|_{\CO,\frac{\beta'}{2}, r''}
<\frac{6r'\|{\cL U}\|_{\CO,\beta', r'}}{(r'-r'')\cdot \beta'}\|\cL V\|_{\CO,\frac{\beta'}{2}, r''}
<\frac{8r'\|{\cL U}\|^2_{\CO,\beta', r'}}{(r'-r'')\cdot \beta'}.\qed
$$


\section{An abstract KAM-like theorem}\label{sec_KAM}

In this section, we give an abstract KAM-like theorem for pairs of holomorphic involutions near a fixed point, which are pairwise conjugate by an anti-holomorphic involution. From this, we obtain the existence of a lot of analytic invariant sets in a neighborhood of the fixed point. This is the core of the proof of Theorem \ref{thm_inv}.

\subsection{Sequences of quantities}\label{subsec_seq}

With fixed $s\in \N^*$, $0<r_0<\frac14$, $0<\varepsilon_0<r_0^2$, $\zeta_0:=\varepsilon_0^{\frac13}$,
define the sequences, with $\nu\in\mathbb N$, $\{\varepsilon_\nu\}$, $\{\beta_\nu\}$, $\{\tilde\beta_\nu\}$, $\{\zeta_\nu\}$, $\{r_\nu\}$ and $\{K_\nu\}$ by:
\begin{equation}\label{sequences}
\begin{array}{lll}
\displaystyle \varepsilon_{\nu+1}:=\varepsilon_\nu^{\frac54}, & \displaystyle \beta_\nu:=\varepsilon_\nu^{\frac1{40s}}, & \displaystyle \tilde\beta_\nu:=16\beta_{\nu+1}=16\varepsilon_\nu^{\frac1{32s}},  \\[2mm]
\displaystyle \zeta_{\nu+1}:=\zeta_\nu+\varepsilon_\nu^{\frac1{3}},& \displaystyle r_{\nu+1}:= r_\nu-\frac{r_0}{2^{\nu+2}}, &\displaystyle  K_\nu:=\frac{|\ln\varepsilon_{\nu}|}{\left|\ln\left(\frac{7r_{\nu}+r_{\nu+1}}{8r_{\nu}}\right)\right|}.
\end{array}
\end{equation}
Between $r_{\nu}$ and $r_{\nu+1}$, we define
\begin{equation*}
r_{\nu}^{(m)}:=r_{\nu+1}+\frac{m}{8}(r_{\nu}-r_{\nu+1}),\quad m=0,1,\cdots,8,\qquad \tilde r_\nu:=r_{\nu}^{(4)}=\frac{r_\nu+r_{\nu+1}}{2}.
\end{equation*}
We assume that $\varepsilon_0$ is small enough such that
\begin{equation}\label{varepsilon_0_small}
\left(\frac{\left|\ln\varepsilon_0\right|}{\left|\ln\left(\frac{7}{8}+\frac{r_1}{8r_0}\right)\right|}+2\right)\frac{(16s+1)^{16s}  \varepsilon_0^{\frac{1}{2400 s^2}}}{(r_0-r_1)r_1}<1.
\end{equation}

\begin{lemma}\label{lemma_varepsilon_nu_small}
 Under the assumption (\ref{varepsilon_0_small}), we have
\begin{equation}\label{varepsilon_nu_small}
\left(\frac{\left|\ln\varepsilon_\nu\right|}{\left|\ln\left(\frac{7}{8}+\frac{r_{\nu+1}}{8r_\nu}\right)\right|}+2\right)\frac{(16s+1)^{16s}  \varepsilon_\nu^{\frac{1}{2400 s^2}}}{(r_\nu-r_{\nu+1})r_{\nu+1}}<1,\quad \nu\in\N.
\end{equation}
\end{lemma}
\proof The above inequality holds for $\nu=0$ under the assumption (\ref{varepsilon_0_small}). Now, assume that, for some $\nu_*\in\N$, we have
\begin{equation}\label{assump_nu_star}
\left(\frac{\left|\ln\varepsilon_{\nu_*}\right|}{\left|\ln\left(\frac{7}{8}+\frac{r_{\nu_*+1}}{8r_{\nu_*}}\right)\right|}+2\right)\frac{(16s+1)^{16s}  \varepsilon_{\nu_*}^{\frac{1}{2400 s^2}}}{(r_{\nu_*}-r_{\nu_*+1})r_{\nu_*+1}}<1.
\end{equation}

The definition of sequence $\{\varepsilon_\nu\}$ implies that, for $\nu\in\N$,
\begin{equation}\label{ele1}
\left|\ln\varepsilon_{\nu+1}\right|=\frac54\left|\ln\varepsilon_{\nu}\right|,\quad \varepsilon_{\nu+1}^{\frac{1}{2400 s^2}}= \varepsilon_{\nu}^{\frac54\cdot \frac{1}{2400 s^2}},
\end{equation}
and the definition of $\{r_\nu\}$ implies that, for $\nu\in\N$,
\begin{equation}\label{ele2}
r_{\nu+1}=r_0\left(1-\sum_{j=0}^{\nu}\frac{1}{2^{j+2}}\right),\quad r_{\nu+1}-r_{\nu+2}=\frac{r_{\nu}-r_{\nu+1}}{2}.
\end{equation}
Hence, for $\nu\in\N$, we have
\begin{equation}\label{ele3}
\frac{r_{\nu}}{r_{\nu+1}}\leq \frac43,\quad 1-\frac{r_{\nu+1}}{r_{\nu}}=\frac{1}{2^{\nu+1}+2}.
\end{equation}
Indeed, it is true for $\nu=0$, and for $\nu\in\N^*$,
$$
\frac{r_{\nu}}{r_{\nu+1}}= \frac{1-\sum_{j=0}^{\nu-1}\frac{1}{2^{j+2}}}{1-\sum_{j=0}^{\nu}\frac{1}{2^{j+2}}}
=1+\frac{1}{2^{\nu+1}+1}
\leq\frac{4}{3},\quad 1-\frac{r_{\nu+1}}{r_{\nu}} =\frac{1}{2^{\nu+1}+2}.
$$
Then, we obtain
\begin{equation}\label{ele4}
\left|\frac{\ln\left(1-\frac1{8}(1-\frac{r_{\nu+1}}{r_{\nu}})\right)}{\ln\left(1-\frac1{8}(1-\frac{r_{\nu+2}}{r_{\nu+1}})\right)}\right| \leq\frac{\frac32\cdot\frac1{8}(1-\frac{r_{\nu+1}}{r_{\nu}})}{\frac34\cdot\frac1{8}(1-\frac{r_{\nu+2}}{r_{\nu+1}})}=4-\frac{2}{2^{\nu}+1}\leq 4.
\end{equation}
By (\ref{ele1}) --- (\ref{ele4}), combining with the assumption (\ref{assump_nu_star}), we have
\begin{eqnarray*}
 & & \left(\frac{\left|\ln\varepsilon_{\nu_*+1}\right|}{\left|\ln\left(\frac{7}{8}+\frac{r_{\nu_*+2}}{8r_{\nu_*+1}}\right)\right|}+2\right)\frac{(16s+1)^{16s}  \varepsilon_{\nu_*+1}^{\frac{1}{2400 s^2}}}{(r_{\nu_*+1}-r_{\nu_*+2})r_{\nu_*+2}} \\
&=&  \left(\frac54\frac{\left|\ln\varepsilon_{\nu_*}\right|}{\left|\ln\left(\frac{7}{8}+\frac{r_{\nu_*+1}}{8r_{\nu_*}}\right)\right|}\cdot\left|\frac{\ln\left(1-\frac1{8}(1-\frac{r_{\nu_*+1}}{r_{\nu_*}})\right)}{\ln\left(1-\frac1{8}(1-\frac{r_{\nu_*+2}}{r_{\nu_*+1}})\right)}\right|+2\right) \cdot \, \frac{2r_{\nu_*+1}}{r_{\nu_*+2}}\frac{(16s+1)^{16s}  \varepsilon_{\nu_*}^{\frac54\cdot \frac{1}{2400 s^2}}}{(r_{\nu_*}-r_{\nu_*+1})r_{\nu_*+1}}\\
&\leq& \frac54\cdot4\cdot 2\cdot\frac{4}{3} \cdot \varepsilon_{\nu_*}^{\frac14\cdot \frac{1}{2400 s^2}} \left(\frac{\left|\ln\varepsilon_{\nu_*}\right|}{\left|\ln\left(\frac{7}{8}+\frac{r_{\nu_*+1}}{8r_{\nu_*}}\right)\right|} +2\right)\frac{(16s+1)^{16s}  \varepsilon_{\nu_*}^{\frac{1}{2400 s^2}}}{(r_{\nu_*}-r_{\nu_*+1})r_{\nu_*+1}}\\
&<&\frac{40}{3} \varepsilon_{\nu_*}^{\frac14\cdot \frac{1}{2400 s^2}} <1.
\end{eqnarray*}
The last inequality follows from (\ref{assump_nu_star}) since $\varepsilon_{\nu_*}^{\frac14\cdot \frac{1}{2400 s^2}}<\varepsilon_{0}^{\frac14\cdot \frac{1}{2400 s^2}}<(16s+1)^{-4s}<\frac{3}{40}$. \qed

\subsection{Iteration argument}

Let the sequences of quantities be given as in \re{sequences}.
For $\omega\in{\cL O}_0\subset \, ]-r_0^2,r_0^2[$ with $|{\cL O}_0|>r_0^2$, we consider the pair of germs of holomorphic involutions $\tau_{0}^{(1)},\tau_{0}^{(2)}:(\C^2,0)\to (\C^2,0)$,
i.e., they satisfy $\tau_0^{(k)}\circ\tau_0^{(k)}={\rm Id}$, $k=1,2$. Recalling the notation in \re{Orbeta}, we set ${\CO}_{0}(r_0,\beta_0):={\CO}_0 \, \cap \, ]-r_0^2+\beta_0, r_0^2-\beta_0[$.
We assume that they are of the form
\begin{eqnarray}
\tau_{0}^{(1)}(\xi,\eta)&=&\left(\begin{array}{l}
                           e^{\frac{\rm i}2\alpha_{0}(\xi\eta)}\eta+p_{0}(\xi,\eta) \\
                           e^{-\frac{\rm i}2\alpha_{0}(\xi\eta)}\xi+q_{0}(\xi,\eta)
                          \end{array}
\right),\label{tau_1}\\
\tau_{0}^{(2)}(\xi,\eta)&=& \left(\rho\circ\tau_{0}^{(1)}\circ\rho\right)(\xi,\eta)=\left(\begin{array}{l}
e^{-\frac{\rm i}2\alpha(\xi\eta)}\eta+\bar p_{0}(\xi,\eta) \\
e^{\frac{\rm i}2\alpha(\xi\eta)}\xi+\bar q_{0}(\xi,\eta)
\end{array}
\right),\label{tau_2}
\end{eqnarray}
where, for the fixed $s\in\N^*$,
\begin{itemize}
  \item {\it (the principal part)} $\alpha_0=\alpha_0(\xi\eta)\in{\CA}^\R_{\beta_0,r_0}(\CO_0)$ with
  \begin{equation}\label{alpha_bdd-0}
  \alpha_0(\omega)\in \, \left]-\frac18,4\pi+\frac18\right[, \quad  \omega\in{\CO}_0(r_0,\beta_0),
  \end{equation}
\begin{eqnarray}
\|\alpha_0(\xi\eta)\|_{\CO_0,\beta_0,r_0}&<& 4\pi+\frac14,\label{alpha_bdd-c-0}\\
\left|\alpha_0^{(s)}-s!\right|_{{\CO}_{0}(r_0,\beta_0)}&<&\frac{s!}{16}, \label{Ds0}\\
\left|\alpha_0^{(k)}\right|_{{\CO}_{0}(r_0,\beta_0)}&<& \frac{1}{16},\qquad 1\leq k\leq  s-1,  \  {\rm if} \   s\geq 2, \label{Ds-0}\\
\left|\alpha_0^{(k)}\right|_{{\CO}_{0}(r_0,\beta_0)}&<& \frac{r_0^{-1}}{4},\qquad s+1\leq k\leq 16s,\label{Ds+0}
  \end{eqnarray}
 \item {\it (the perturbation)} $p_0=p_0(\xi,\eta)$, $q_0=q_0(\xi,\eta)\in{\CA}_{\beta_0,r_0}(\CO_0)$ with
 \begin{equation}\label{esti_pq_sec4}
 \|p_0\|_{\CO_0,\beta_0,r_0}, \ \|q_0\|_{\CO_0,\beta_0,r_0}\leq \frac{\varepsilon_0}{10},\quad \|e^{\frac{\rm i}2\alpha_{0}(\xi\eta)}\eta q_{0}+e^{-\frac{\rm i}2\alpha_{0}(\xi\eta)}\xi p_{0}\|_{\CO_0,\beta_0,r_0}<\frac{\varepsilon_0^{\frac32}}{3}.
 \end{equation}
\end{itemize}
\begin{remark}
For instance, $\alpha_0(z)=\lambda+z^s+\sum_{j=s+1}^m c_j z^j$, with $\lambda\in[0,4\pi[$, arbitrary $m\geq s+1$, $c_j\in\R$, and $r_0$ sufficiently small, is an example of such a function satisfying (\ref{alpha_bdd-0}) -- (\ref{Ds+0}).
\end{remark}

\begin{remark} In (\ref{esti_pq_sec4}), besides the smallness of perturbation $(p_0, q_0)$, the smallness of $e^{\frac{\rm i}2\alpha_{0}(\xi\eta)}\eta q_{0}+e^{-\frac{\rm i}2\alpha_{0}(\xi\eta)}\xi p_{0}$, called the ``skew term" of $\tau_0^{(1)}$, is crucial in the iteration.
\end{remark}

We also consider the germ of map
$\sigma_0=\tau_0^{(1)}\circ\tau_0^{(2)}$. It is reversible with respect to both $ \rho$ and $\tau_0^{(1)}$ since
$$
\rho\circ\sigma_0\circ\rho= \rho\circ\tau_0^{(1)}\circ\rho\circ\rho\circ\tau_0^{(2)}\circ\rho=\tau_0^{(2)}\circ\tau_0^{(1)}=\sigma_0^{-1}=\tau_0^{(1)}\circ\sigma_0\circ\tau_0^{(1)}.
$$
We can write $\sigma_{0}=\tau_{0}^{(1)}\circ \tau_{0}^{(2)}$ as
$$\sigma_{0}(\xi,\eta)=\left(\begin{array}{l}
                           e^{{\rm i}\alpha_{0}(\xi\eta)}\xi+f_{0}(\xi,\eta) \\
                           e^{-{\rm i}\alpha_{0}(\xi\eta)}\eta+g_{0}(\xi,\eta)
                          \end{array}
\right).$$
It will be shown in Section \ref{sec_rev} (see Lemma \ref{lem_appro} and Corollary \ref{cor_fg}) that, if $\varepsilon_0$ satisfies (\ref{varepsilon_0_small}), then $f_0,g_0\in {\CA}_{\tilde\beta_0,r_0^{(7)}}(\CO_0)$ with
$$\|f_0\|_{\CO_0,\tilde\beta_0,r_0^{(7)}}, \ \|g_0\|_{\CO_0,\tilde\beta_0,r_0^{(7)}}\leq \frac{\varepsilon_0}{4}.$$

\begin{prop}[Iteration scheme]\label{prop_KAM}  Assume that $\varepsilon_0>0$ satisfies (\ref{varepsilon_0_small}).
There exist a sequence of sets $\{{\cL O}_{\nu}\}$ with ${\cL O}_{\nu}\subset \, ]-r_\nu^2,r_\nu^2[$ satisfying that
\begin{equation}\label{mes}
 {\cL O}_{\nu+1}\subset{\cL O}_{\nu} \, \cap \, ]-r_{\nu+1}^2,r_{\nu+1}^2[, \quad \left|({\cL O}_{\nu}\setminus{\cL O}_{\nu+1}) \, \cap \, \left]-r_{\nu+1}^2+\beta_{\nu+1}, \, r_{\nu+1}^2-\beta_{\nu+1}\right[\right|<\varepsilon_\nu^{\frac{1}{100s^2}},
\end{equation}
and a sequence of maps $\{\sigma_{\nu}\}$ given by $\sigma_{\nu}=\tau_{\nu}^{(1)}\circ \tau_{\nu}^{(2)}$ with
\begin{equation}\label{tau_nu}
\tau_{\nu}^{(1)}(\xi,\eta)=\left(\begin{array}{l}
                           e^{\frac{\rm i}2\alpha_{\nu}(\xi\eta)}\eta+p_{\nu}(\xi,\eta) \\
                           e^{-\frac{\rm i}2\alpha_{\nu}(\xi\eta)}\xi+q_{\nu}(\xi,\eta)
                          \end{array}
\right), \quad \tau_{\nu}^{(2)}=\rho\circ\tau_{\nu}^{(1)}\circ\rho,
\end{equation}
satisfying $\tau_{\nu}^{(k)}\circ\tau_{\nu}^{(k)}={\rm Id}$, $k=1,2$, such that the following holds.
\begin{enumerate}
\item $\alpha_\nu=\alpha_\nu(\xi\eta)\in {\CA}^\R_{\beta_\nu,r_\nu}(\CO_\nu)$ satisfies that
\begin{eqnarray}
\left|\alpha_\nu^{(s)}-s!\right|_{{\CO}_{\nu}(r_\nu,\beta_\nu)}&<&  \left(\frac{1}{16}+\zeta_\nu\right) s!,\label{Dsnu}\\
\left|\alpha_\nu^{(k)}\right|_{{\CO}_{\nu}(r_\nu,\beta_\nu)}&<& \frac{1}{16}+\zeta_\nu,\qquad 1 \leq k\leq s-1, \  if  \ s\geq 2,\label{Ds-nu}\\
\left|\alpha_\nu^{(k)}\right|_{{\CO}_{\nu}(r_\nu,\beta_\nu)}&<& \left(\frac{1}{4}+\zeta_\nu\right)r_\nu^{-1},\qquad s+1 \leq k\leq 16s,\label{Ds+nu}
\end{eqnarray}
and for $0\leq k \leq 16 s$,
\begin{equation}\label{error_alpha-nu}
\left\|(\alpha_{\nu+1}-\alpha_\nu)^{(k)}\right\|_{\CO_{\nu+1},\beta_{\nu+1},r_{\nu+1}}< \varepsilon_\nu^{\frac{1}{3}}.
\end{equation}

\smallskip

\item $p_\nu$, $q_\nu\in {\CA}_{\beta_\nu,r_\nu}(\CO_\nu)$ satisfy that
\begin{equation}\label{esti_Lpq}
\|p_{\nu}\|_{\CO_\nu,\beta_\nu,r_\nu}, \|q_{\nu}\|_{\CO_\nu,\beta_\nu,r_\nu}< \frac{\varepsilon_{\nu}}{10},\quad
\|e^{\frac{\rm i}2\alpha_{\nu}(\xi\eta)}\eta q_{\nu}+e^{-\frac{\rm i}2\alpha_{\nu}(\xi\eta)}\xi p_{\nu}\|_{\CO_\nu,\beta_\nu, r_\nu}<\frac{\varepsilon_\nu^{\frac32}}{3}.
\end{equation}

\item The reversible map (w.r.t. $\rho$)
$\sigma_{\nu}=\tau_{\nu}^{(1)}\circ \tau_{\nu}^{(2)}$ has the form
$$\sigma_{\nu}(\xi,\eta)=\left(\begin{array}{l}
                           e^{{\rm i}\alpha_{\nu}(\xi\eta)}\xi+f_{\nu}(\xi,\eta) \\
                           e^{-{\rm i}\alpha_{\nu}(\xi\eta)}\eta+g_{\nu}(\xi,\eta)
                          \end{array}
\right),$$
with
 $\|f_{\nu}\|_{\CO_\nu,\tilde\beta_{\nu},r_\nu^{(7)}}$, $\displaystyle \|g_{\nu}\|_{\CO_\nu,\tilde\beta_{\nu},r_\nu^{(7)}}<\frac{\varepsilon_{\nu}}4$.
\smallskip

\item There is a sequence of transformations $\{\psi_{\nu}\}$ of the form
\begin{equation}\label{trans_Prop}
\psi_{\nu}(\xi,\eta)=({\rm Id}+{\cL U}_{\nu})(\xi,\eta)=\left(\begin{array}{l}
\xi+u_{\nu}(\xi,\eta)\\
\eta+v_{\nu}(\xi,\eta)
\end{array}\right),
\end{equation}
with ${\cL U}_{\nu}\in ({\CA}^\R_{\beta_{\nu+1},r_{\nu+1}}({\cL O}_{\nu+1}))^2$ satisfying
$$\|u_{\nu}\|_{\CO_{\nu+1},\beta_{\nu+1},r_{\nu+1}}, \ \|v_{\nu}\|_{\CO_{\nu+1},\beta_{\nu+1},r_{\nu+1}}<\frac{\varepsilon_{\nu}^{\frac{49}{50}}}2$$
such that, for every $\omega\in{\cL O}_{\nu+1}(r_{\nu+1},\beta_{\nu+1})$, $\psi_{\nu}:{\cL C}^{r_{\nu+1}}_{\omega,\beta_{\nu+1}}\to {\cL C}^{r_{\nu}}_{\omega,\beta_{\nu}}$ and, on ${\cL C}^{r_{\nu+1}}_{\omega,\beta_{\nu+1}}$,
$$\sigma_{\nu+1}=\psi_{\nu}^{-1}\circ\sigma_{\nu}\circ\psi_{\nu},\quad \tau^{(k)}_{\nu+1}=\psi_{\nu}^{-1}\circ\tau^{(k)}_{\nu}\circ\psi_{\nu}, \;\ k=1,2.$$
\end{enumerate}
\end{prop}

\begin{remark}\label{rem_nu} According to the definition of sequence $\varepsilon_{\nu}=\varepsilon_0^{\left(\frac54\right)^\nu}$ and the fact that $(\frac54)^x-1>\frac{x}{8}$, $\forall \ x>0$, we have
\beq\label{sumeps}
\sum_{\nu\geq 0}\varepsilon_\nu^{\varsigma}= \varepsilon_0^{\varsigma}\sum_{\nu\geq 0}\varepsilon_0^{\varsigma\left[\left(\frac54\right)^\nu-1\right]}< \varepsilon_0^{\varsigma}\sum_{\nu\geq 0}\varepsilon_0^{\frac{\varsigma\nu}{8}}=\frac{\varepsilon_0^{\varsigma}}{1-\varepsilon_0^{\frac{\varsigma}{8}}},\quad \forall \ 0<\varsigma<1.
\eeq
In particular, $\sum_{\nu\geq 0}\varepsilon_\nu^{\frac13}<\frac{\varepsilon_0^{\frac13}}{1-\varepsilon_0^{\frac{1}{24}}}<\frac1{240}.$
Since $\left\|(\alpha_{\nu+1}-\alpha_\nu)^{(k)}\right\|_{\CO_{\nu+1},\beta_{\nu+1},r_{\nu+1}}< \varepsilon_\nu^{\frac13}$, $0\leq k \leq 16s$, 
we obtain, according to \re{alpha_bdd-0} and \re{alpha_bdd-c-0}, for $\forall \ \nu\in \N$,
\begin{equation}\label{al_nu_sum}
\alpha_\nu(\omega)\in\, \left]-\frac14,4\pi+\frac14\right[, \quad \forall \ \omega\in {\CO}_{\nu}(r_\nu,\beta_\nu), \qquad  \|\alpha_\nu\|_{\CO_\nu,\beta_{\nu}, r_\nu}<4\pi+\frac12.
\end{equation}
Moreover, in view of the definition of $\{\zeta_\nu\}$, we see that $\zeta_\nu<\frac1{240}$ for every $\nu\in\N$,
which implies that $\frac{1}{16}+\zeta_\nu<\frac1{15}$, $\frac{1}{4}+\zeta_\nu<\frac12$.
Then, by \re{Dsnu} -- \re{Ds+nu},
\begin{eqnarray}
\left|\alpha_\nu^{(s)}-s!\right|_{{\CO}_{\nu}(r_\nu,\beta_\nu)}&<&  \frac{s!}{15} ,\label{Dsnu_sum}\\
\left|\alpha_\nu^{(k)}\right|_{{\CO}_{\nu}(r_\nu,\beta_\nu)}&<& \frac{1}{15},\qquad 1 \leq k\leq s-1, \  if  \ s\geq 2,\label{Ds-nu_sum}\\
\left|\alpha_\nu^{(k)}\right|_{{\CO}_{\nu}(r_\nu,\beta_\nu)}&<& \frac{r_\nu^{-1}}{2},\qquad  s+1 \leq k\leq 16s.\label{Ds+nu_sum}
\end{eqnarray}
%
\end{remark}

\subsection{Proof of Proposition \ref{prop_KAM}}

Suppose that, at the $(\nu+1)^{\rm th}$ step, $\nu\geq 0$, we have
$$\tau_{\nu}^{(1)}  =\left(\begin{array}{l}
e^{\frac{\rm i}2\alpha_{\nu}(\xi\eta)}\xi+p_{\nu}(\xi,\eta)\\
e^{-\frac{\rm i}2\alpha_{\nu}(\xi\eta)}\eta+q_{\nu}(\xi,\eta)
\end{array}\right),\quad \tau_{\nu}^{(2)}=\rho\circ\tau_{\nu}^{(1)}\circ\rho,\quad \sigma_\nu=  \tau_{\nu}^{(1)}\circ \tau_{\nu}^{(2)},
$$
as described in Proposition \ref{prop_KAM}.
Our aim is to construct the transformation $\psi_{\nu}$ as in \re{trans_Prop},
such that $\sigma_{\nu+1}:=\psi_{\nu}^{-1}\circ\sigma_{\nu}\circ\psi_{\nu}$, $\tau^{(k)}_{\nu+1}:=\psi_{\nu}^{-1}\circ\tau^{(k)}_{\nu}\circ\psi_{\nu}$, $k=1,2$, possess similar properties as those of $\sigma_{\nu}$, $\tau^{(k)}_{\nu}$. This will describe an iteration step, hence will give the proof of Proposition \ref{prop_KAM} .

Before starting the construction of $\psi_{\nu}$, we first introduce another type of transformations that conjugate the pairs of involutions (\ref{tau_nu}) to a perturbation of a new integrable pair. The main feature of the new involutions is that the new perturbation part is much smaller than the initial one, provided that the initial skew term smallness condition is satisfied.



For given $0<r_+<r<\frac14$ and fixed $s\in\N^*$,
 assume that $\varepsilon>0$ satisfies
 $$\left(\frac{\left|\ln\varepsilon\right|}{\left|\ln\left(\frac{7}{8}+\frac{r_{+}}{8r}\right)\right|}+2\right)\frac{(16s+1)^{16s}  \varepsilon^{\frac{1}{2400 s^2}}}{(r-r_{+})r_{+}}<1.$$
 Let us set \begin{equation}\label{new_def_beta}
 	\beta\in [\varepsilon^{\frac{1}{40s}},  \varepsilon^{\frac{1}{60s}}], \quad
 	\beta_+=\beta^{\frac54}\in [\varepsilon^{\frac{1}{32s}}, \varepsilon^{\frac{1}{48s}}],
 \end{equation}
\begin{equation}\label{r_m}
r^{(m)}:=r_++\frac{m}{8}(r-r_+),\quad m=0,1,\cdots,8, \qquad  \tilde r:=r^{(4)}=\frac{r+r_+}{2}.
\end{equation}
Given $\CO \subset  ]-r^2, r^2[$, we consider the involutions  $\tau_1,\tau_2=\rho\circ\tau_1\circ\rho$~:
\begin{equation}\label{gen-invol}
\tau_1(\xi,\eta)=\left(\begin{array}{l}
e^{\frac{\rm i}2\alpha(\xi\eta)}\eta+ p(\xi,\eta) \\
e^{-\frac{\rm i}2\alpha(\xi\eta)}\xi+ q(\xi,\eta)
\end{array}\right), \quad \tau_2(\xi,\eta)=\left(\begin{array}{l}
	e^{-\frac{\rm i}2\alpha(\xi\eta)}\eta+\bar p(\xi,\eta) \\
	e^{\frac{\rm i}2\alpha(\xi\eta)}\xi+\bar q(\xi,\eta)
\end{array}
\right)
\end{equation}
with $\alpha=\alpha(\xi\eta)\in{\CA}^\R_{\beta,r}(\CO)$ satisfying \re{al_nu_sum} -- \re{Ds+nu_sum} as $\alpha_\nu$,
together with $p,q\in{\CA}_{\beta,r}(\CO)$ satisfying
\begin{equation}\label{esti_pq}\|p\|_{\CO,\beta,r}, \ \|q\|_{\CO,\beta,r}<\frac{\varepsilon}{10}.
\end{equation}
\begin{remark}\label{rem_quantities_nu}
Here, all assumptions of \rp{prop_KAM} are satisfied but 
the smallness condition of the skew term $e^{\frac{\rm i}2\alpha(\xi\eta)}\eta q+e^{-\frac{\rm i}2\alpha(\xi\eta)}\xi p$ in (\ref{esti_Lpq}).
\end{remark}

\begin{thm}(Main step)\label{thm_trans_cohomo}
Given $\delta\in ] \, 80\varepsilon^{\frac1{60s}}, 1  \, [$, let
\begin{equation}\label{SD_general}
\CO_{\delta}:=\left\{\omega \in \CO  :  \  |e^{{\rm i}n\alpha(\omega)}-1|\geq\delta, \quad  \forall \ 0<|n|\leq K+1, \quad K:=\frac{|\ln\varepsilon|}{\left|\ln(r^{(7)}/r)\right|}\right\}.
\end{equation}
There exists a transformation $\psi$ of the form
$$\psi(\xi,\eta)=({\rm Id}+{\CU})(\xi,\eta)=\left(\begin{array}{l}
\xi+ u(\xi,\eta) \\
\eta+ v(\xi,\eta)
\end{array}
\right)$$
with $u,v\in {\CA}^\R_{\beta_+,r_+}(\CO_\delta)$ satisfying
$$\|u\|_{\CO_\delta,\beta_+,r_+}, \ \|v\|_{\CO_\delta,\beta_+,r_+} <\frac{\varepsilon^{\frac{49}{50}}}2,$$
such that for every $\omega\in{\cL O}_\delta(r_+, \beta_+)= {\cL O}_\delta \, \cap \, ]-r_+^2+\beta_+, r_+^2-\beta_+[$,
$\psi$ is biholomorphic on ${\cL C}^{r_+}_{\omega,\beta_+}$ with $\psi\left({\cL C}^{r_+}_{\omega,\beta_+}\right)\subset{\cL C}^{r}_{\omega,\beta}$,
and, on ${\cL C}^{r_+}_{\omega,\beta_+}$,
$$(\psi^{-1}\circ\tau_1\circ\psi)(\xi,\eta)=\left(\begin{array}{l}
e^{\frac{\rm i}2\alpha_+(\xi\eta)}\eta+ p_+(\xi,\eta) \\
e^{-\frac{\rm i}2\alpha_+(\xi\eta)}\xi+ q_+(\xi,\eta)
\end{array}
\right),\quad (\xi,\eta)\in {\cL C}^{r_+}_{\omega,\beta_+},$$
where $\alpha_+=\alpha_+(\xi\eta)\in {\CA}^\R_{\beta_+,r_+}(\CO_\delta)$, with
\begin{equation}\label{error_alpha}
\left\|(\alpha_+-\alpha)^{(k)}\right\|_{\CO_{\delta},\beta_+,r_+}<\frac{\varepsilon^{\frac13}}{10},\quad 0\leq k\leq 16s,
\end{equation}
and $p_+$, $q_+\in {\CA}_{\beta_+,r_+}(\CO_\delta)$, with
\begin{equation}\label{esti_p_+q_+}
\|p_+\|_{\CO_{\delta},\beta_+,r_+}, \|q_+\|_{\CO_{\delta},\beta_+,r_+}< \frac{\varepsilon^{\frac{61}{32}}}{2}+24(K+1)\delta^{-1}\|e^{\frac{\rm i}2\alpha}\eta q+e^{-\frac{\rm i}2\alpha}\xi p\|_{ \CO,\beta, r},
\end{equation}
\begin{equation}\label{esti_Lp_+q_+}
\|e^{\frac{\rm i}2\alpha_+}\eta q_++e^{-\frac{\rm i}2\alpha_+}\xi p_+\|_{\CO_{\delta},\beta_+, r_+}<\varepsilon^{\frac{61}{32}}.
\end{equation}
\end{thm}

\begin{remark} Theorem \ref{thm_trans_cohomo} can be applied in two ways, which are indeed two cases described in Section \ref{sec_app_KAM}.
\begin{itemize}
\item If the skew term $e^{\frac{\rm i}2\alpha}\eta q+e^{-\frac{\rm i}2\alpha}\xi p$ of $\tau_1$ satisfies
$$
\delta^{-1}\|e^{\frac{\rm i}2\alpha}\eta q+e^{-\frac{\rm i}2\alpha}\xi p\|_{\CO,\beta, r}<\varepsilon^{1+\varsigma},\quad 0<\varsigma<1,
$$
then, according to \re{esti_p_+q_+}, Theorem \ref{thm_trans_cohomo} describes an iteration step (as in Proposition \ref{prop_KAM}) with $r=r_\nu$, $r_+=r_{\nu+1}$, $\varepsilon=\varepsilon_\nu$, and we can simply take $\beta=\varepsilon^{\frac{1}{40s}}$.
\item If it is not the case, we cannot apply our Iteration scheme. However, noting that \re{esti_p_+q_+} implies $\|p_+\|_{\CO_{\delta},\beta_+,r_+}$, $\displaystyle \|q_+\|_{\CO_{\delta},\beta_+,r_+}<\frac{\varepsilon^{\frac{49}{50}}}{10}$,
we see that the Iteration scheme is applicable to $\psi^{-1}\circ\tau_1\circ\psi$ in view of \re{esti_Lp_+q_+}. Hence, Theorem \ref{thm_trans_cohomo} describes a preliminary step of iteration. In this case, we need to take $\beta>\varepsilon^{\frac{1}{40s}}$ since the new perturbation may be of size $\varepsilon^{\frac{49}{50}}$. This is why $\beta$ is defined in an interval as in (\ref{new_def_beta}).
\end{itemize}
\end{remark}

We postpone the proof of Theorem \ref{thm_trans_cohomo} to Section \ref{sec_trans}. The rest of the section is devoted to the proof of \rp{prop_KAM} from \rt{thm_trans_cohomo}.

\smallskip

We want to conjugate the involution $\tau_\nu^{(1)}$ to a new one $\tau_{\nu+1}^{(1)}$. To do so, we need to exclude some parameters so we define the new parameter set as follow~:
\begin{equation}\label{tilde_O}
{\cL O}_{\nu+1}:=\left\{\omega\in{\cL O}_\nu \, \cap \, ]-r_{\nu+1}^2,r_{\nu+1}^2[ \, :
|e^{{\rm i}n\alpha_\nu(\omega)}-1|> \varepsilon_\nu^{\frac{1}{64s}}, \;
\forall \ 0<|n| \leq K_\nu+1\right\}.
\end{equation}
In order to measure its size, let us first recall Pyartli's lemma~:
\begin{lemma}\label{pyartli}\cite{pyartli-dioph, russmann-weak}
Let $f: [a,b] \mapsto \mathbb R$ with $a<b$ be a $q$-times continuously differentiable function satisfying
$$|f^{(q)}(t)|\geq\delta,\quad \forall \  t\in [a,b]$$
for some $q\in \N^*$ and $\delta> 0$. Then, for any $ A>0$,
$$\left|\left\{t\in [a,b] : |f(t)|\leq  A\right\}\right|\leq 4\left(q!\frac{ A}{2\delta}\right)^{\frac1q}.$$
\end{lemma}

We then have
\begin{lemma}\label{lemma_meas_esti}
$\left|({\cL O}_{\nu}\setminus{\cL O}_{\nu+1}) \, \cap \, ] \, -r_{\nu+1}^2+\beta_{\nu+1},r_{\nu+1}^2-\beta_{\nu+1} \, [ \right|<\varepsilon_\nu^{\frac1{80s^2}}$.
\end{lemma}
\proof Since $|\alpha_{\nu}(\omega)|<4\pi+1$, it is sufficient to bound from above the measure of the parameter set $\bigcup_{0<|n| \leq K_\nu+1}{\cL R}_{\nu,n}$
, where
$${\cL R}_{\nu, n}:=\left\{\omega\in{\cL O}_\nu \, \cap \, ]-r_{\nu+1}^2,r_{\nu+1}^2[ \,  : \left|n\alpha_\nu(\omega)-2k\pi\right|\leq\frac{3\varepsilon_\nu^{\frac{1}{64s}}}2,\;\  k\in\Z , \ |k|\leq3|n| \right\}.$$

In view of (\ref{Dsnu}) -- (\ref{Ds+nu}), we see $\inf_{\omega\in{\cL O}_\nu}\left|(n\alpha_\nu)^{(s)}(\omega)\right|\geq \frac34 |n| s!$.
Applying \rl{pyartli} with $q=s$, we have
$$|{\cL R}_{\nu, n}|\leq 3|n|\cdot 4\left(\frac{\varepsilon_\nu^{\frac{1}{64s}}}{2 |n| }\right)^{\frac1s}\leq 20(K_\nu+1)^{1-\frac1s}\varepsilon_\nu^{\frac{1}{64s^2}}.$$
Therefore, we obtain
$$\bigcup_{0<|n| \leq K_\nu+1}\left|{\cL R}_{\nu, n}\right|\leq 40(K_\nu+1)^{2-\frac{1}{s}}\varepsilon_\nu^{\frac{1}{64s^2}}<\varepsilon_\nu^{\frac{1}{80s^2}},$$
noting that (\ref{varepsilon_nu_small}) implies that
$$ (K_\nu+1)^2\varepsilon_\nu^{\frac1{64s^2}-\frac{1}{80s^2}} =  \left(\frac{|\ln\varepsilon_\nu|}{\left|\ln\left(\frac{7r_\nu+r_{\nu+1}}{8r_\nu}\right)\right|}+1\right)^2\varepsilon_\nu^{\frac1{320s^2}}
< (16s+1)^{-16s}
< \frac{1}{40}.\qed$$

Applying Theorem \ref{thm_trans_cohomo} to $\tau_1=\tau_\nu^{(1)}$ with $\delta=\varepsilon_\nu^{\frac{1}{64s}}>80 \varepsilon_\nu^{\frac{1}{60s}}$, we obtain a transformation $\psi_\nu$ of the form
$$\psi_\nu(\xi,\eta)=({\rm Id}+{\CU}_\nu)(\xi,\eta)=\left(\begin{array}{l}
\xi+ u_\nu(\xi,\eta) \\
\eta+ v_\nu(\xi,\eta)
\end{array}
\right)$$
with $u_\nu$, $v_\nu\in \CA^\R_{\beta_{\nu+1}, r_{\nu+1}}(\CO_{\nu+1})$ satisfying
$$\|u_\nu\|_{\CO_{\nu+1},\beta_{\nu+1},r_{\nu+1}},\ \|v\|_{\CO_{\nu+1},\beta_{\nu+1},r_{\nu+1}}<\frac{\varepsilon_\nu^{\frac{49}{50}}}2,$$ such that, on ${\cL C}_{\omega,\beta_{\nu+1}}^{r_{\nu+1}}$, $\omega\in {\cL O}_{\nu+1}(r_{\nu+1},\beta_{\nu+1})$,
 $\psi_\nu:{\cL C}_{\omega,\beta_{\nu+1}}^{r_{\nu+1}}\to{\cL C}_{\omega,\beta_\nu}^{r_\nu}$ is injective and holomorphic, and
 $$(\psi_\nu^{-1}\circ\tau_\nu^{(1)}\circ\psi_\nu)(\xi,\eta)=\left(\begin{array}{l}
e^{\frac{\rm i}2\alpha_{\nu+1}(\xi\eta)}\eta+ p_{\nu+1}(\xi,\eta) \\
e^{-\frac{\rm i}2\alpha_{\nu+1}(\xi\eta)}\xi+ q_{\nu+1}(\xi,\eta)
\end{array}
\right),
$$
where $\alpha_{\nu+1}=\alpha_{\nu+1}(\xi\eta)\in \CA^\R_{\beta_{\nu+1}, r_{\nu+1}}(\CO_{\nu+1})$ with
$$
\|(\alpha_{\nu+1}-\alpha_\nu)^{(k)}\|_{\CO_{\nu+1},\beta_{\nu+1},r_{\nu+1}}<\varepsilon_\nu^{\frac13},\quad 0\leq k\leq 16s,
$$
which, combining with (\ref{Dsnu}) -- (\ref{Ds+nu}), implies that
 \begin{eqnarray*}
|\alpha_{\nu+1}^{(s)}-s!|_{{\cL O}_{\nu+1}(r_{\nu+1},\beta_{\nu+1})}&<&  \left(\frac{1}{16}+\zeta_{\nu+1}\right) s!,\\
|\alpha_{\nu+1}^{(k)}|_{{\cL O}_{\nu+1}(r_{\nu+1},\beta_{\nu+1})}&<& \frac{1}{16}+\zeta_{\nu+1},\qquad  1 \leq k\leq s-1, \ {\rm if} \  s\geq 2,\\
|\alpha_{\nu+1}^{(k)}|_{{\cL O}_{\nu+1}(r_{\nu+1},\beta_{\nu+1})}&<& \left(\frac{1}{4}+\zeta_{\nu+1}\right)r^{-1},\qquad  s+1 \leq k\leq 16s.
\end{eqnarray*}
In view of \re{esti_Lpq} and (\ref{esti_p_+q_+}), (\ref{esti_Lp_+q_+}), we have the new perturbation $p_{\nu+1}$ and $q_{\nu+1}$ satisfy that
\begin{align*}
\|p_{\nu+1}\|_{\CO_{\nu+1},\beta_{\nu+1},r_{\nu+1}}, \|q_{\nu+1}\|_{\CO_{\nu+1},\beta_{\nu+1},r_{\nu+1}}
< &\frac{\varepsilon_\nu^{\frac{61}{32}}}{2}+24(K_\nu+1)\varepsilon_\nu^{-\frac{1}{64s}}\|e^{\frac{\rm i}2\alpha_\nu}\eta q_\nu+e^{-\frac{\rm i}2\alpha_\nu}\xi p_\nu\|_{\CO_{\nu},\beta_\nu, r_\nu},\\
< &\frac{\varepsilon_\nu^{\frac54}}{10}= \frac{\varepsilon_{\nu+1}}{10},
\end{align*}
and the new skew term satisfies
$$\|e^{\frac{\rm i}2\alpha_{\nu+1}}\eta q_{\nu+1}+e^{-\frac{\rm i}2\alpha_{\nu+1}}\xi p_{\nu+1}
\|_{\CO_{\nu+1},\beta_{\nu+1}, r_{\nu+1}}
<\varepsilon_\nu^{\frac{61}{32}}<\frac{\varepsilon_\nu^{\frac{15}{8}}}{3}=\frac{\varepsilon_{\nu+1}^{\frac32}}{3}.$$
According to Lemma \ref{lem_reve_psi}, $u_\nu$, $v_\nu\in \CA^\R_{\beta_{\nu+1}, r_{\nu+1}}(\CO_{\nu+1})$ implies that $\rho\circ \psi_\nu=\psi_\nu\circ\rho$. Then, for
$$\tau_{\nu+1}^{(1)}:=\psi_\nu^{-1}\circ\tau_\nu^{(1)}\circ\psi_\nu,\quad  \tau_{\nu+1}^{(2)}:=\psi_\nu^{-1}\circ\tau_\nu^{(2)}\circ\psi_\nu=\rho\circ \tau_{\nu+1}^{(1)} \circ\rho,$$
we still have $\tau_{\nu+1}^{(k)}\circ\tau_{\nu+1}^{(k)}={\rm Id}$, $k=1,2$. Moreover, in view of Lemma \ref{lem_reve},
 for $\sigma_{\nu+1}=\tau_{\nu+1}^{(1)}\circ\tau_{\nu+1}^{(2)}$, it is still reversible w.r.t. $\rho$.

\section{Proof of Theorem \ref{thm_inv}}\label{sec_proof}

This section is dedicated to the proof of Theorem \ref{thm_inv} by applying the Iteration scheme \rp{prop_KAM}.

\subsection{Preliminary normalization}

Let us consider the pair of involutions $\tau^o_1$ and $\tau^o_2$ given in \re{tau1} and \re{tau2} (and the reversible map $\sigma_o=\tau^o_1\circ\tau^o_2$). In order to start the Iteration scheme \rp{prop_KAM}, we need the involutions to be in a well prepared form as in (\ref{tau_1}) and (\ref{tau_2}).

First of all, for any $N>s$, there exists a holomorphic transformation $\check\Psi$ in the neighborhood of origin, tangent to identity up to order $2$, with $\check\Psi\circ\rho=\rho\circ\check\Psi$, such that
\begin{equation}\label{check_tau}
\check\tau_1(\xi,\eta):=\left(\check\Psi^{-1}\circ \tau^o_1\circ \check\Psi\right)(\xi,\eta)=\left(\begin{array}{l}
		e^{\frac{\rm i}2\check\alpha(\xi\eta)}\eta+\check p(\xi,\eta) \\
		e^{-\frac{\rm i}2\check\alpha(\xi\eta)}\xi+\check q(\xi,\eta)\end{array}\right),
\end{equation}
where we have
\begin{equation}\label{check_alpha}
\check\alpha(\xi\eta):=\lambda+(\xi\eta)^s+\sum_{n=s+1}^{N} c_n(\xi\eta)^n,\quad c_n\in \R,
\end{equation}
and convergent power series at the origin
\begin{equation}\label{check_p_q}
\check p(\xi,\eta)=\sum_{l+j\geq 2N+2\atop{l,j\geq 0}}\breve{\check p}_{l,j}\xi^l\eta^j,\quad \check q(\xi,\eta)=\sum_{l+j\geq 2N+2\atop{l,j\geq 0}}\breve{\check q}_{l,j}\xi^l\eta^j.
\end{equation}
Indeed, since $\frac{\lambda}{\pi}\in\R\setminus\Q$,
by classical normal form theory \cite{Arn2, ilyashenko-yakovenko-book, moser-webster}, combining with the fact that $\tau^o_1\circ\tau^o_1={\rm Id}$,
there is a polynomial transformation $\Psi_{PD}$, tangent to identity up to order $2$ at the origin (composed by finitely many steps of normalization in the sense of Poincar\'e-Dulac) satisfying $\Psi_{PD}\circ\rho=\rho\circ\Psi_{PD}$, such that
$$\left(\Psi_{PD}^{-1}\circ \tau^o_1\circ\Psi_{PD}\right)(\xi,\eta)=\left(\begin{array}{l}
		(e^{\frac{\rm i}2\lambda}+ \tilde C(\xi\eta))\eta+\tilde p(\xi,\eta) \\
		(e^{\frac{\rm i}2\lambda}+\tilde C(\xi\eta))^{-1}\xi+\tilde q(\xi,\eta)\end{array}\right),\quad \tilde C(z)=\sum_{j=s}^{N}\tilde c_j z^j.$$
Here $\tilde p$, $\tilde q$ are holomorphic at the origin and of order $\geq2N+2$ there.
We note that $e^{\frac{\rm i}2\lambda}+ \tilde C(z)$ is actually the truncation of $\Lambda(z)$ in (\ref{formal_tau1}) -- (\ref{formal_sigma}). Recalling the non-degeneracy assumption of Theorem \ref{thm_inv}, we see that $\tilde c_s\neq 0$.
According to the proof of Theorem 3.4 of \cite{moser-webster}, we can
change $e^{\frac{\rm i}2\lambda}+ \tilde C$ to $(e^{\frac{\rm i}2\lambda}+ \tilde C)\mu^{-2}$ by applying the transformation
\begin{equation}\label{scaling_MW}
(\xi,\eta)\mapsto(\mu(\xi\eta) \xi, \mu^{-1}(\xi\eta)\eta),
\end{equation}
where $\mu=\mu(\xi\eta)$ is the fourth root $\mu(\xi\eta):=\left((e^{\frac{\rm i}2\lambda}+ \tilde C(\xi\eta))(e^{-\frac{\rm i}2\lambda}+ \bar{\tilde C}(\xi\eta))\right)^{\frac14}$. We see that $\mu(\xi\eta)$ is sufficiently close to $1$ and hence well-defined since $(\xi,\eta)$ belongs to a sufficiently small neighborhood of the origin. A direct computation shows that
$$(e^{\frac{\rm i}2\lambda}+ \tilde C)\mu^{-2}=(e^{\frac{\rm i}2\lambda}+ \tilde C)^{\frac12}(e^{-\frac{\rm i}2\lambda}+ \bar{\tilde C})^{-\frac12}=e^{\frac{\rm i}2\lambda}(1+e^{-\frac{\rm i}2\lambda}\tilde C)^{\frac12}(1+e^{\frac{\rm i}2\lambda}\bar{\tilde C})^{-\frac12}.$$
We consider the principal determination of the logarithm to be defined on $\C\setminus \R^{-}$.  Let $A(t)$ be holomorphic germ vanishing at the origin. If $t$ is small enough, there are real numbers $\tilde b_k$ such that
	$$
	\ln(1+A(t))-\ln (1+\bar A(t))=\sum_{k\geq 1}\frac{A^k(t)-\bar A^{k}(t)}{k}={\rm i}\sum_{k\geq 1} \tilde b_kt^k
	$$
	and it converges at the origin. Applying this, we have
$$\ln\left((1+e^{-\frac{\rm i}2\lambda}\tilde C)^{\frac12}(1+e^{\frac{\rm i}2\lambda}\bar{\tilde C})^{-\frac12} \right)=\frac12\left( e^{-\frac{\rm i}2\lambda} \tilde c_s - e^{\frac{\rm i}2\lambda}\bar{\tilde c}_s\right)(\xi\eta)^s+\sum_{n\geq s+1} b_n (\xi\eta)^n,$$
with $\{b_n\}\subset {\rm i}\R$ coefficients of a convergent power series.
Define
$$\check\alpha(\xi\eta):=\lambda-{\rm i}\left(e^{-\frac{\rm i}2\lambda}\tilde c_s-e^{\frac{\rm i}2\lambda}\bar{\tilde c}_s\right)(\xi\eta)^s-2{\rm i}\sum_{n= s+1}^N b_n (\xi\eta)^n,$$
which is of the form (\ref{check_alpha}) up to a scaling on the neighborhood of origin.
By rewriting $(1+e^{-\frac{\rm i}2\lambda}\tilde C)^{\frac12}(1+e^{\frac{\rm i}2\lambda}\bar{\tilde C})^{-\frac12}$ as
$$(1+e^{-\frac{\rm i}2\lambda}\tilde C)^{\frac12}(1+e^{\frac{\rm i}2\lambda}\bar{\tilde C})^{-\frac12} =e^{\frac{\rm i}2(\check\alpha-\lambda)}\exp\left\{\sum_{n\geq N+1} b_n (\xi\eta)^n\right\},$$
we see that
$$
(e^{\frac{\rm i}2\lambda}+ \tilde C)\mu^{-2}-e^{\frac{\rm i}2\check\alpha}= e^{\frac{\rm i}2\lambda}\cdot e^{\frac{\rm i}2(\check\alpha-\lambda)} \cdot \left(\exp\left\{\sum_{n\geq N+1} b_n (\xi\eta)^n\right\}-1\right),
$$
which contains terms of $(\xi,\eta)$ of order $\geq 2N+2$. It is similar for $(e^{\frac{\rm i}2\lambda}+ \tilde C)^{-1}\mu^{2}-e^{-\frac{\rm i}2\check\alpha}$.
Hence, we obtain $\check\tau_1$ in (\ref{check_tau}), up to a scaling on the neighborhood of origin.

\

The following lemma shows that, for $N\geq 16s$ large enough, there exists $0<r_*<\frac14$ sufficiently small, depending on the coefficients $c_{n}$, $n=s+1,\cdots,N$, such that (\ref{varepsilon_0_small}) can be satisfied in the two cases (where $ A:=10 \max\{|\check p|_{r_*}, \ |\check q|_{r_*}\}$):
\begin{itemize}
	\item $\varepsilon_0= A^{\frac{49}{50}}$, $r_0=\frac{3}{4}r_*$ and $r_1=\frac{9}{16}r_*$,
	\item $\varepsilon_0= A$, $r_0=r_*$ and $r_1=\frac{3}{4}r_*$.
\end{itemize}



\begin{lemma}\label{lemma_r_varepsilon}
If $N$ is large enough, there exists $0<r_*<\frac14$ sufficiently small such that
\begin{enumerate}
	\item [(1)] $\sup_{|z|<r_*^2} |\check\alpha(z)-\lambda|<\frac18$, $\sup_{|z|<r_*^2}|\check\alpha^{(s)}(z)-s!|<\frac{s!}{20}$,
	\item [(2)] $\sup_{|z|<r_*^2}|\check\alpha^{(k)}(z)|<\frac{1}{20}$, $1\leq k\leq s-1$, if $s\geq 2$,
	\item [(3)] $\sup_{|z|<r_*^2}|\check\alpha^{(k)}(z)|<\frac{r_*^{-1}}{4}$, $ s+1\leq k\leq 16s$,
	\item [(4)] $A=10 \max\{|\check p|_{r_*}, \ |\check q|_{r_*}\}$ satisfies that
\begin{equation}\label{smallnessA}
	\left(\frac{\left|\ln A\right|}{\left|\ln\left(\frac{7}{8}+\frac18\cdot\frac34\right)\right|}+2\right)\frac{(16s+1)^{16s}   A^{\frac{49}{50}\cdot\frac{1}{2400 s^2}}}{(\frac{3}{4}r_*-\frac{9}{16}r_*)\cdot\frac{9}{16}r_*}<1.
\end{equation}
\end{enumerate}
\end{lemma}
\proof In view of (\ref{check_alpha}), it is easy to see that (1) -- (3) are satisfied for $\check\alpha$ for any $N\geq 16s$ and $r_*$ sufficiently small. Let us choose $N= 4900s^2$.
We recall that all terms of $\check p$ and $\check q$ are of order $\geq 2N+2$. According to the definitions of norms in  \re{norm-0}, if we replace $r_*$ by $r_k':=2^{-k}r_*$, then $ A_k':=10 \max\{|\check p|_{r'_k}, \ |\check q|_{r'_k}\}$ satisfies that $ A'_k\leq A\cdot \left(\frac{r_k'}{r_*}\right)^{2N+2}=2^{-(2N+2)k} A$.
Since $0<r_*<\frac14$, we have that
$\left((\frac{3}{4}r_*-\frac{9}{16}r_*)\cdot\frac{9}{16}r_*\right)^{-1}=\frac{256}{27}r_*^{-2}<10r_*^{-2}$.
To show (4), it is sufficient to show that
\begin{eqnarray}
(16s+1)^{16s} \left(\frac{\left|\ln A\right|}{\left|\ln\left(\frac{31}{32}\right)\right|}+2\right)  A^{\frac{49}{50}\cdot\frac{1}{4800 s^2}} &<& 1, \label{r_varepsilon_1}\\
10r_*^{-2} A^{\frac{49}{50}\cdot\frac{1}{4800 s^2}}  &<& 1.\label{r_varepsilon_2}
\end{eqnarray}
Replacing $r_*$ by $r_k'$ in (\ref{r_varepsilon_1}) and (\ref{r_varepsilon_2}), we see that, as $k\to \infty$,
\begin{align*}
 & (16s+1)^{16s} \left(\frac{\left|\ln A'_k\right|}{\left|\ln\left(\frac{31}{32}\right)\right|}+2\right) ( A_k')^{\frac{49}{50}\cdot\frac{1}{4800 s^2}}\\
< \, &(16s+1)^{16s}\left(\frac{\ln(2)\cdot (9800s^2+2)k+\left|\ln A\right|}{\left|\ln\left(\frac{31}{32}\right)\right|}+2\right)  2^{-\frac{49}{50}\cdot\frac{9800s^2+2}{4800s^2}k}  A^{\frac{49}{50}\cdot\frac{1}{4800 s^2}}\to 0,
\end{align*}

$$ 10r_k'^{-2}\cdot ( A_k')^{\frac{49}{50}\cdot\frac{1}{4800 s^2}}<10r_*^{-2} 2^{2k} \cdot 2^{-\frac{49}{50}\cdot\frac{9800s^2+2}{4800s^2}k}  A^{\frac{49}{50}\cdot\frac{1}{4800 s^2}}
 =10r_*^{-2} A^{\frac{1}{4800 s^2}} \cdot 2^{-(\frac{49}{50}\cdot\frac{9800s^2+2}{4800s^2}-2)k}
 \to 0.$$
Hence, there exists a $k_*\in\N^*$ such that if we replace $r_*$ by $r'_{k_*}$, then \re{r_varepsilon_1} and \re{r_varepsilon_2} are both satisfied.\qed

\subsection{Application of KAM-like Theorem}\label{sec_app_KAM}
Take $r_*$ as in Lemma \ref{lemma_r_varepsilon}. Since $\check p$ and $\check q$ are convergent power series, in view of (\ref{comp_norm}), we have that, for any $\beta_*\in [  A^{\frac{1}{40s}},   A^{\frac{1}{60s}}]$,
\begin{equation}\label{esti_check_p_q}
\|\check p\|_{]-r_*^2+\beta_*,r_*^2-\beta_*[, \, \beta_*, \, r_*}, \  \|\check q\|_{]-r_*^2+\beta_*,r_*^2-\beta_*[, \, \beta_*, \, r_*}\leq \max\{|\check p|_{r_*}, \ |\check q|_{r_*}\}=\frac{ A}{10}.
\end{equation}

\begin{itemize}
  \item {\bf Case 1. Small skew term}
\end{itemize}
If we have
\begin{equation}\label{crossing_check}
\|e^{-\frac{\rm i}2\check\alpha}\xi\check p+e^{\frac{\rm i}2\check\alpha}\eta\check q\|_{]-r_*^2+\beta_*,r_*^2-\beta_*[, \, \beta_*, \, r_*}<\frac{ A^{\frac32}}{3},
\end{equation}
then, from \rl{lemma_r_varepsilon}
and \re{esti_check_p_q}, with $\alpha_0(\xi\eta):=\check\alpha(\xi\eta)$, $p_0:=\check p$, $q_0:=\check q$, $r_0:=r_*$, $\varepsilon_0:= A$, $\beta_0:=\varepsilon_0^{\frac{1}{40s}}= A^{\frac{1}{40s}}$ and ${\CO}_0:=]  -r_0^2, r_0^2 \, [$, we have that
$\alpha_0(\xi\eta)\in {\CA}^\R_{\beta_0, r_0}(\CO_0)$, $p_0$, $q_0\in {\CA}_{\beta_0, r_0}(\CO_0)$, and (\ref{alpha_bdd-0}) -- \re{esti_pq_sec4} are satisfied. Since, for $r_1:=\frac34 r_0=\frac34 r_*$,
\re{varepsilon_0_small} is satisfied from Lemma \ref{lemma_r_varepsilon},
we can apply Proposition \ref{prop_KAM} to $\tau_0^{(1)}:=\check\tau_1$ (hence $\tau_0^{(2)}=\rho\circ\tau_0^{(1)}\circ\rho$ and $\sigma=\tau_0^{(1)}\circ\tau_0^{(2)}$) on ${\cL C}_{\omega, \beta_0}^{r_0}$,  $\omega\in\CO_0(r_0,\beta_0)= {\CO}_0 \, \cap \, ]  -r_0^2+\beta_0, r_0^2-\beta_0 \, [$.

\smallskip

\begin{itemize}
  \item {\bf Case 2. Non-small skew term}
\end{itemize}
Now assume that \re{crossing_check} is not satisfied.
In view of (\ref{check_alpha}) and \rl{lemma_r_varepsilon}, we see that
$$|e^{\pm\frac{\rm i}2\check\alpha}|_{\omega,\beta_*}=|e^{\pm\frac{\rm i}2(\check\alpha-\lambda)}|_{\omega,\beta_*}\leq\sum_{k\geq 0}\frac{1}{k!}\left(\frac{|\check\alpha-\lambda|_{\omega,\beta_*}}2\right)^k\leq e^{\frac1{16}},\quad \omega\in \, ]-r_*^2+\beta_*,r_*^2-\beta_*[.$$
Then, by (\ref{esti_check_p_q}),
$$\|e^{-\frac{\rm i}2\check\alpha}\xi\check p+e^{\frac{\rm i}2\check\alpha}\eta\check q\|_{\omega,\beta_*,r_*}<\frac{e^{\frac1{16}}}{2}\cdot\frac{A}{10}<\frac{ A}{10}.$$

With $r=r_*$, $r_+=\frac34 r_*$ and $r_*^{(j)}:=r_++\frac{j}{8}(r_*-r_+)$, $j=0,\cdots,8$, $\delta:=100 A^{\frac{1}{60s}}$, $\beta:= A^{\frac{1}{60s}}$, $\beta_+:=  A^{\frac{49}{50}\cdot\frac{1}{40s}}\in[A^{\frac{1}{32s}}, A^{\frac{1}{48s}}]$, and
$$\CO_\delta := \left\{\omega\in \, ]  -r_*^2, r_*^2 \, [ \, : |e^{{\rm i}n\check\alpha(\omega)}-1|>\delta, \quad  \forall \ 0<|n|\leq K_*+1:=\frac{|\ln A|}{\left|\ln(r_*^{(7)}/r_*)\right|}+1\right\},$$
we apply Theorem \ref{thm_trans_cohomo} to $\check\tau_1$ for $\omega\in\CO_\delta(r_+,\beta_+):=\CO_\delta \, \cap \, ]-r_+^2+\beta_+,r_+^2-\beta_+[$. We obtain, for all $\om\in \CO_\delta(r_+,\beta_+)$, a biholomorphic transformation $\check\psi={\rm Id}+\check{\CU}:{\cL C}^{r_+}_{\omega,\beta_+}\to {\cL C}^{r}_{\omega,\beta}$, with $\check\CU\in ({\CA}^\R_{\beta_+,r_+}(\CO_\delta))^2$ and $\|\check \CU\|_{\omega,\beta_+,r_+}<\frac{A^{\frac{49}{50}}}2$.
Furthermore, there are $\check\alpha_+=\check\alpha_+(\xi\eta)\in {\CA}^\R_{\beta_+,r_+}(\CO_\delta)$,  $\check p_+$, $\check q_+\in {\CA}_{\beta_+,r_+}(\CO_\delta)$ such that
$$(\check\psi^{-1}\circ\check\tau_1\circ\check\psi)(\xi,\eta)=\left(\begin{array}{l}
e^{\frac{\rm i}2\check\alpha_+(\xi\eta)}\eta+ \check p_+(\xi,\eta) \\
e^{-\frac{\rm i}2\check\alpha_+(\xi\eta)}\xi+ \check q_+(\xi,\eta)
\end{array}
\right),\quad (\xi,\eta)\in {\cL C}^{r_+}_{\omega,\beta_+}.$$
They satisfy
\begin{equation}\label{error_check_alpha}
\left\|(\check\alpha_+-\check\alpha)^{(k)}\right\|_{\CO_{\delta},\beta_+,r_+}<\frac{ A^{\frac13}}{10},\quad 0\leq k\leq 16s,
\end{equation}
$$\|\check p_+\|_{\CO_{\delta},\beta_+,r_+}, \|\check q_+\|_{\CO_{\delta},\beta_+,r_+}< \frac{ A^{\frac{61}{32}}}{2}+18(K_*+1)\delta^{-1}\|e^{\frac{\rm i}2\check\alpha}\eta \check q+e^{-\frac{\rm i}2\check\alpha}\xi \check p\|_{]-r_*^2+\beta_*,r_*^2-\beta_*[,\beta_*, r_*}< \frac{ A^{\frac{49}{50}}}{10},$$
$$\|e^{-\frac{\rm i}2\check\alpha_+}\xi\check p_++e^{\frac{\rm i}2\check\alpha_+}\eta\check q_+\|_{\CO_{\delta},\beta_+,r_+}< A^{\frac{61}{32}}<\frac{ A^{\frac{49}{50}\cdot\frac32}}{3}.$$
Hence, for $\alpha_0:=\check\alpha_+$, $p_0:=\check p_+$, $q_0:=\check q_+$, $r_0:=r_+=\frac34r_*$, $\varepsilon_0:= A^{\frac{49}{50}}$, and $$\beta_0:=\beta_+=\varepsilon_0^{\frac{1}{40s}}= A^{\frac{49}{50}\cdot\frac{1}{40s}},$$
we have \re{esti_pq_sec4} for ${\CO}_0:=\CO_\delta(r_0,\beta_0)$.
According to \re{error_check_alpha}, we obtain (\ref{alpha_bdd-0}) -- (\ref{Ds+0}) from \rl{lemma_r_varepsilon}. In particular, Lemma \ref{lemma_r_varepsilon} (1) and \re{error_check_alpha} imply that $$
|\alpha_0(\cdot)-\lambda|_{{\CO}_0}<\frac{1}{4}.
$$
For $r_0=\frac34r_*$, $r_1:=\frac34 r_0=\frac9{16} r_*$,
\re{varepsilon_0_small} is verified by Lemma \ref{lemma_r_varepsilon}, then we can apply Proposition \ref{prop_KAM} to $\tau_0^{(1)}:=\check\psi^{-1}\circ\check\tau_1\circ\check\psi$ (hence $\tau_0^{(2)}=\rho\circ\tau_0^{(1)}\circ\rho$ and $\sigma=\tau_0^{(1)}\circ\tau_0^{(2)}$) on ${\cL C}_{\omega, \beta_0}^{r_0}$, $\omega\in{\CO}_0(r_0,\beta_0)={\CO}_0$.
According to \rl{lemma_r_varepsilon}}, we see that $|\check\alpha^{(s)}(\omega)|\geq \frac{19}{20}s!$ for $\omega\in \, ]-r_*^2,r_*^2[$, then we deduce from Pyartli's lemma (\rl{pyartli}) that
\begin{equation}\label{meas_preliminary}
\left| \,  ]-r_0^2+\beta_0,r_0^2-\beta_0 \, [ \, \setminus {\CO}_0\right|< A^{\frac{1}{80s^2}}.\end{equation}
The proof of \re{meas_preliminary} is similar to that of Lemma \ref{lemma_meas_esti}.

\smallskip

Let us define $\check\psi={\rm Id}$ in Case 1. In Case 2, we have, as above, $\check\psi={\rm Id}+\check{\CU}$ with $\check\CU\in ({\CA}^\R_{\beta_+,r_+}(\CO_\delta))^2$. In both cases, we define $\tilde\Phi:=\check\Psi\circ \check\psi$.
To summarize, for the involutions $\tau^o_1$ given in (\ref{tau1}), we have

\begin{prop} There exists $r_0>0$ and there exists $\varepsilon_0>0$ satisfying \re{varepsilon_0_small} with $r_1=\frac{3 r_0}{4}$, and there exists a set ${\cL O}_{0}\subset ]-r_0^2,r_0^2[$ with
$$ \left| \,  ]-r_0^2+\beta_0,r_0^2-\beta_0 \, [ \, \setminus {\CO}_0\right|< \varepsilon_0^{\frac{49}{50}\cdot\frac{1}{80s^2}}$$
for $\beta_0=\varepsilon_0^{\frac{1}{40s}}$, such that the following holds for $\omega\in{\cL O}_{0}$.

There exists a transformation $\tilde\Phi: {\cL C}_{\omega,\beta_0}^{r_0}\to \C^2$, with $\tilde\Phi\circ\rho=\rho\circ\tilde\Phi$, such that the involution
$$\tau_0^{(1)}(\xi,\eta)=(\tilde\Phi^{-1}\circ\tau^o_1\circ\tilde\Phi)(\xi,\eta)=\left(\begin{array}{l}
e^{\frac{\rm i}2\alpha_{0}(\xi\eta)}\eta+p_{0}(\xi,\eta) \\
e^{-\frac{\rm i}2\alpha_{0}(\xi\eta)}\xi+q_{0}(\xi,\eta)
\end{array}
\right),$$
with $\alpha_0(\xi\eta)\in{\CA}_{\beta_0,r_0}^\R(\CO_0), p_0 ,q_0\in{\CA}_{\beta_0,r_0}(\CO_0)$ and (\ref{alpha_bdd-0}) -- \re{esti_pq_sec4} satisfied. In particular,
\begin{equation}\label{alpha0_lambda}
|\alpha_0(\cdot)-\lambda|_{\CO_0}<\frac14.
\end{equation}
\end{prop}

\smallskip

As in \re{tau2}, $\tau_0^{(2)}$ is obtained by
$\tau_0^{(2)}=\rho\circ\tau_0^{(1)}\circ\rho$,
and $\sigma_0=\tau_0^{(1)}\circ\tau_0^{(2)}$. Since $\tilde\Phi$ commutes with $\rho$, we see that $\tau_0^{(2)}$ is still an involution and, in view of Lemma \ref{lem_reve}, $\sigma_0$ is still reversible with respect to $\rho$.

\smallskip

By applying Proposition \ref{prop_KAM} to $\tau_0^{(1)}$ (hence $\tau_0^{(2)}=\rho\circ\tau_0^{(1)}\circ\rho$ and $\sigma_0=\tau_0^{(1)}\circ\tau_0^{(2)}$) on ${\cL C}_{\omega, \beta_0}^{r_0}$,  $\omega\in{\CO}_0(r_0,\beta_0)$, we
get sequences of involutions $\{\tau_\nu^{(1)}\}_{\nu\in\N}$ and $\{\tau_\nu^{(2)}\}_{\nu\in\N}$
(and hence $\{\sigma_{\nu}\}_{\nu\in\N}$) on ${\cL C}_{\omega, \beta_\nu}^{r_\nu}$ and a sequence of holomorphic transformations $\{\psi_{\nu}\}_{\nu\in\N}$
of the form $\psi_{\nu}={\rm Id}+{\cL U}_{\nu}$ with ${\cL U}_{\nu}\in ({\CA}^\R_{\beta_{\nu+1},r_{\nu+1}}({\cL O}_{\nu+1}))^2$ and $\|{\cL U}_{\nu}\|_{\CO_{\nu+1},\beta_{\nu+1},r_{\nu+1}}<\varepsilon_{\nu}^{\frac{49}{50}}$,
such that, for all  $\omega\in{\cL O}_{\nu+1}(r_{\nu+1},\beta_{\nu+1})$, $\psi_{\nu}:{\cL C}^{r_{\nu+1}}_{\omega,\beta_{\nu+1}}\to {\cL C}^{r_{\nu}}_{\omega,\beta_{\nu}}$,
$$\sigma_{\nu+1}=\psi_{\nu}^{-1}\circ\sigma_{\nu}\circ\psi_{\nu},\quad \tau^{(k)}_{\nu+1}=\psi_{\nu}^{-1}\circ\tau^{(k)}_{\nu}\circ\psi_{\nu}, \;\ k=1,2.$$

Recall that $r_{\nu+1}= r_\nu-2^{-(\nu+2)}r_0$, we see that $r_{\nu}\to\frac{r_0}{2}=:R$ as $\nu\to\infty$.
By (\ref{mes}), we can see that ${\cL O}_\nu\to {\cL O}_{\infty}(R)$ for some ${\cL O}_{\infty}(R)\subset \left]-\frac{r_0^2}{4}, \frac{r_0^2}{4}\right[=\left]-R^2, R^2\right[$. Moreover, by \re{mes} and \re{sumeps},
$$
\left|\left]-R^2, R^2\right[\setminus {\cL O}_{\infty}(R)\right|<\sum_{\nu\geq 0}\varepsilon_\nu^{\frac{1}{100s^2}}<\frac{\varepsilon_0^{\frac{1}{100s^2}}}{1-\varepsilon_0^{\frac{1}{800s^2}}}<2\varepsilon_0^{\frac{1}{100s^2}}.$$
Hence, we have (\ref{mesure_asymp}) by noting that the Lebesgue density of ${\cL O}_{\infty}(R)$ in $\left]-R^2, R^2\right[$ satisfies
$$1-\frac{\varepsilon_0^{\frac{1}{100s^2}}}{R^2}<\frac{|{\cL O}_{\infty}(R)|}{2R^2}<1,$$
since $r_0=r_*$ or $\frac34r_*$, and by (\ref{smallnessA}), we have
 $\varepsilon_0^{\frac{1}{100s^2}}\leq A^{\frac{49}{50}\cdot\frac{1}{100s^2}}<\left(\frac{r_0}2\right)^{2400s^2\cdot\frac{1}{100s^2}}= R^{24}$.

For any $\nu\in\N$, let $\Psi_\nu:=\check\psi\circ\psi_{0}\circ \cdots \circ \psi_\nu$, which is well defined and injective on ${\cL C}^{r_{\nu+1}}_{\omega,\beta_{\nu+1}}$ for every $\omega\in{\cL O}_{\nu+1}(r_{\nu+1},\beta_{\nu+1})$.
By Lemma \ref{lem_reve_psi}, since $\check\CU\in ({\CA}^\R_{\beta_+,r_+}(\CO_\delta))^2$ and $u_{j},\, v_{j} \in{\CA}_{\beta_{j+1},r_{j+1}}^\R({\cL O}_{j+1})$, $j=0,1,\cdots,\nu$, then $\check\psi\circ\rho=\rho\circ\check\psi$ and $\psi_j\circ\rho=\rho\circ\psi_j$. Hence $\Psi_\nu\circ\rho=\rho\circ\Psi_\nu$.

\begin{lemma}\label{lem_Psi_nu} For every $\nu\in\N$, $\|\Psi_{\nu+1}-\Psi_\nu\|_{\CO_{\nu+2},\beta_{\nu+2},r_{\nu+2}}<\varepsilon_{\nu+1}^{\frac45}$.
\end{lemma}

\proof  With the smallness of $A$ verified in Lemma \ref{lemma_r_varepsilon} and recalling that $\varepsilon_{0}=A$ or $A^{\frac{49}{50}}$, we have, by Lemma \ref{lem_norm-pre},
\begin{eqnarray*}
\|\Psi_0\|_{\CO_1,\beta_1,r_1}&\leq&\|\check\psi\circ\psi_0-\check\psi\|_{\CO_1,\beta_1,r_1}+\|\check\psi\|_{\CO_1,\beta_1,r_1}\\
&\leq& \frac{3r_0 \|\check\psi\|_{\CO_0,\beta_{0},r_{0}}\varepsilon_{0}^{\frac{49}{50}}}{2(r_0-r_1)\beta_0}+\|\check\psi\|_{\CO_0,\beta_0,r_0}\\
&\leq&(2r_0+A^{\frac{49}{50}})\left(\frac{3r_0 \varepsilon_{0}^{\frac{49}{50}}}{2(r_0-r_1)\beta_0}+1\right) \, \leq \, 3r_0.
\end{eqnarray*}

Let us show the lemma by induction on $\nu$.
When $\nu=0$, $\Psi_0=\check\psi\circ\psi_0$ and $\Psi_1=\Psi_0\circ\psi_1$.
For any $(\xi,\eta)\in{\cL C}^{r_1}_{\omega,\beta_{1}}$, we have $\psi_{1}(\xi,\eta)\in {\cL C}^{r_0}_{\omega,\beta_0}$.
Since (\ref{varepsilon_0_small}) implies that (\ref{smallness_beta}) holds for $\beta'=\beta_{1}$, $\beta''=\beta_{2}$, $r'=r_{1}$, $r''=r_{2}$,
 we have, by Lemma \ref{lem_norm-pre},
\begin{eqnarray*}
\|\Psi_1-\Psi_0\|_{\CO_2,\beta_{2},r_{2}}&=&\|\Psi_0(\xi+u_1,\eta+v_1)-\Psi_0(\xi,\eta)\|_{\CO_2,\beta_{2},r_{2}} \\
&<&\frac{3r_1 \|\Psi_0\|_{\CO_1,\beta_{1},r_{1}}}{2(r_1-r_2)\beta_1}\varepsilon_{1}^{\frac{49}{50}}
\, < \, \frac{3r_0r_1 \varepsilon_{1}^{\frac{49}{50}}}{(r_1-r_2)\beta_1}  \, < \, \varepsilon_{1}^{\frac45}.
\end{eqnarray*}
Given $k\in \N^*$,
assume that
$\|\Psi_{j+1}-\Psi_j\|_{\CO_{j+2},\beta_{j+2},r_{j+2}}<\varepsilon_{j+1}^{\frac45}$ for $0\leq j\leq k$. Then
\begin{equation}\label{bounded_by_r}
\|\Psi_{k+1}\|_{\CO_{k+2},\beta_{k+2},r_{k+2}}\leq\|\Psi_0\|_{\CO_{1},\beta_{1},r_{1}}+\sum_{j=0}^k\|\Psi_{j+1}-\Psi_j\|_{\CO_{j+2},\beta_{j+2},r_{j+2}}
<3r_0+\sum_{j=0}^k\varepsilon_{j+1}^{\frac45}<4r_0.
\end{equation}
Hence, by Lemma \ref{lem_norm-pre},
\begin{eqnarray*}
\|\Psi_{k+2}-\Psi_{k+1}\|_{\CO_{k+3},\beta_{k+3},r_{k+3}}
&=&\|\Psi_{k+1}(\xi+u_{k+2},\eta+v_{k+2})-\Psi_{k+1}(\xi,\eta)\|_{\CO_{k+3},\beta_{k+3},r_{k+3}}\\
&<&\frac{3r_{k+2}\|\Psi_{k+1}\|_{\CO_{k+2},\beta_{k+2},r_{k+2}}}{2(r_{k+2}-r_{k+3})\beta_{k+2}}\varepsilon_{k+2}^{\frac{49}{50}}\\
&<& \frac{6r_0 r_{k+2}\varepsilon_{k+2}^{\frac{49}{50}}}{(r_{k+2}-r_{k+3})\beta_{k+2}}\,< \, \varepsilon_{k+2}^{\frac45},
\end{eqnarray*}
since (\ref{varepsilon_nu_small}) implies that (\ref{smallness_beta}) holds for $\beta'=\beta_{k+2}$, $\beta''=\beta_{k+3}$, $r'=r_{k+2}$, $r''=r_{k+3}$.\qed

\smallskip

The above lemma shows that, with $\check\Psi$ in (\ref{check_tau}), for every $\omega\in {\cL O}_{\infty}(R)$, the sequence $\{\check\Psi\circ\Psi_{\nu}\}$ converges uniformly to an injective holomorphic mapping $\Psi_{\om}:{\cL C}^{R}_{\omega}\to \C^2$ as it is a Cauchy sequence:
\begin{eqnarray*}
\sup_{(\xi,\eta)\in {\cL C}^{R}_{\omega}}\| \check\Psi\circ\Psi_{\nu}(\xi,\eta)- \check\Psi\circ\Psi_{\nu'}(\xi,\eta) \|&\leq& \|\check\Psi\|_{R}  \|\Psi_{\nu}-\Psi_{\nu'}\|_{\omega,0,R} \\
 &\leq& \|\check\Psi\|_{R}\sum_{j= \nu}^{\nu'}\|\Psi_{j}-\Psi_{j+1}\|_{\omega,\beta_{j+2},r_{j+2}}\\
 &<&\|\check\Psi\|_{R} \sum_{j\geq \nu}\varepsilon_j^{\frac45} \, \to \, 0,\quad \nu\to \infty.
\end{eqnarray*}
We shall denote $\Psi_{\om,\nu}$ the restriction of $\check\Psi\circ\Psi_{\nu}$ to ${\cL C}^{R}_{\omega}$.
Moreover, recalling that $\check\Psi$ is tangent to identity and combining with (\ref{bounded_by_r}), we have, for every $\omega\in{\cL O}_\infty(R)$,
\begin{equation}\label{estim-sizeimage}
\sup_{(\xi,\eta)\in {\cL C}^{R}_{\omega}}\| \Psi_{\om}(\xi,\eta) \|\leq\|\check\Psi\|_{R}\left(\| \Psi_{0}\|_{{\cL O}_{1},\beta_{1},r_{1}} +\sum_{\nu=0}^\infty\|\Psi_{\nu+1}-\Psi_{\nu}\|_{{\cL O}_{\nu+2},\beta_{\nu+2},r_{\nu+2}} \right)< R^{\frac12},
\end{equation}
which implies that $\Psi_{\om}({\cL C}^{R}_{\omega})\subset\Delta_2(0,R^{\frac12})$.

\smallskip

For $\omega\in {\cL O}_{\nu+1}(r_{\nu+1},\beta_{\nu+1})$, we have
$$
\Psi_\nu^{-1}\circ\check\Psi^{-1}\circ\tau^{(1)}_0\circ (\check\Psi\circ\Psi_\nu) =\tau^{(1)}_{\nu+1}=\left(
\begin{array}{c}
e^{\frac{\rm i}2\alpha_{\nu+1}(\xi\eta)}\xi+p_{\nu+1}\\
e^{-\frac{\rm i}2\alpha_{\nu+1}(\xi\eta)}\eta+q_{\nu+1}
\end{array}
\right),
$$
with $\|p_{\nu+1}\|_{\omega,\beta_{\nu+1},r_{\nu+1}}$, $\displaystyle \|q_{\nu+1}\|_{\omega,\beta_{\nu+1},r_{\nu+1}}<\frac{\varepsilon_{\nu+1}}{10}$
and $(\check\Psi\circ\Psi_\nu)\circ\rho=\rho\circ(\check\Psi\circ\Psi_\nu)$ implies that
$$\tau^{(2)}_{\nu+1}=\rho\circ\tau^{(1)}_{\nu+1}\circ\rho=\Psi_\nu^{-1}\circ\check\Psi^{-1}\circ\tau^{(2)}_0\circ(\check\Psi\circ\Psi_\nu),\quad \sigma_{\nu+1}=\tau^{(1)}_{\nu+1}\circ\tau^{(2)}_{\nu+1}=\Psi_\nu^{-1}\circ\check\Psi^{-1}\circ \sigma_0\circ (\check\Psi\circ\Psi_\nu).$$
Hence, for every $\omega\in {\cL O}_{\infty}(R)\subset  \left]-R^2, R^2\right[$, the sequence $\{\alpha_{\nu}\}$ restricted to ${\cL C}_{\omega}^{R}$ converges to a real number $\mu_{\om}=\alpha_{\infty}(\omega)$ with $\alpha_{\infty}:=\alpha_0+\sum_{\nu\geq 1}(\alpha_{\nu}-\alpha_{\nu-1})$. Indeed, since for any $\nu',\nu\in\N$ with $\nu'\geq \nu$,
	$$|\alpha_{\nu}-\alpha_{\nu'}|_{\omega,0}\leq\sum_{j= \nu}^{\nu'}|\alpha_{j}-\alpha_{j+1}|_{\omega,\beta_{j+1}}<\sum_{j\geq \nu}\varepsilon_j^{\frac13}\to 0,\quad \nu\to \infty,$$ it is a Cauchy sequence $\subset ]-1,4\pi+1[$. In particular, combining (\ref{error_alpha-nu}) and \re{alpha0_lambda}, we obtain
$$|\alpha_\infty-\lambda|_{{\cL O}_{\infty}(R)}<|\alpha_0-\lambda|_{{\cL O}_0}+\sum_{\nu\geq 1}\|\alpha_{\nu}-\alpha_{\nu-1}\|_{\CO_{\nu},\beta_{\nu},r_{\nu}}<\frac14+ \sum_{\nu\geq 1} \varepsilon_{\nu-1}^{\frac{1}{3}}<\frac{\pi}{4}.$$
 Furthermore, we have
$$ (\Psi^{-1}_{\omega}\circ \tau_0^{(1)}\circ \Psi_{\om})(\xi,\eta)=\left(
\begin{array}{c}
e^{\frac{\rm i}2\mu_{\omega}}\eta\\
e^{-\frac{\rm i}2\mu_{\omega}}\xi
\end{array}
\right),\quad (\Psi^{-1}_{\omega}\circ  \sigma_0\circ \Psi_{\omega})(\xi,\eta)=\left(
\begin{array}{c}
e^{{\rm i}\mu_{\omega}}\xi\\
e^{-{\rm i}\mu_{\omega}}\eta
\end{array}
\right).$$

\subsection{Whitney-smoothness of holomorphic hyperbolas}

Now let us show the smoothness (in the sense of Whitney) of $\mu_\omega$ and $\Psi_{\omega}$ with respect to $\omega\in \CO_\infty=\CO_\infty(R)$. We refer to \cite{bhs-book}[Chapter 6.1.4] for Whitney-smoothness notions.

Given $r\in \R_+ \setminus\N$, let $k_r:=\lfloor r\rfloor=\max\{k\in\Z:k<r\}$. Given a closed set $\CO\subset \mathbb{R}$, given a family $f=(f_l)_{0\leq l\leq k_r}\in (C^0(\CO))^{k_r+1}$, 
we define its {\it Whitney-$C^r$ norm} to be~:
$$
|f|_{{C_{W}^{r}(\mathcal{O})}}:= \sup_{\omega\in \CO\atop{ 0\leq l\leq k_r}}|f_{l}(\om)|+ \sup_{(\om,\om')\in \CO^2,\ \om\neq\om'\atop{ 0\leq l\leq k_r}}\frac{|f_l(\om)-P_l(\om,\om')|}{|\om-\om'|^{r-l}},
$$
where $P_l$ is an analogue of the $(k_r-l)-{\rm th}$ Taylor polynomial for $f_l$, i.e.,
$$
P_l(\om,\om'):=\sum_{j=0}^{k_r-l}\frac{1}{j!}f_{l+j}(y)(\om-\om')^{j}.
$$
Such a family $f=(f_l)_{0\leq l\leq k_r}\in (C^0(\CO))^{k_r+1}$ is said to define a Whitney-$C^r$ function $f$ if its Whitney-$C^{r}$ norm is finite~:
$$
|f|_{{C_{W}^{r}(\mathcal{O})}}<+\infty.
$$
Given an open set $\CU$ satisfying $\CO\subset\CU\subset \mathbb{R}$, for $g\in C^{k_r}(\CU)$, we define
\begin{eqnarray*}
|g|_{{C^{r}(\mathcal{\CU})}}&:=& \sup_{\omega\in \CU\atop{ 0\leq l\leq k_r}}|g^{(l)}(\om)|+ \sup_{(\om,\om')\in \CU^2 \atop{\om\neq\om'}}\frac{|D^{k_r}g(\om)-D^{k_r}g(\om')|}{|\om-\om'|^{r-k_r}},\\
|g|_{{C^{r}(\mathcal{\CO})}}&:=& \sup_{\omega\in \CO\atop{ 0\leq l\leq k_r}}|g^{(l)}(\om)|+ \sup_{(\om,\om')\in \CO^2 \atop{\om\neq\om'}}\frac{|D^{k_r}g(\om)-D^{k_r}g(\om')|}{|\om-\om'|^{r-k_r}}.
\end{eqnarray*}
We have $|g|_{{C^{r}(\mathcal{\CO})}}\leq|g|_{{C^{r}(\mathcal{\CU})}}$.
According to Chapter 6.1.4 of \cite{bhs-book}, the norms $|\cdot|_{{C_{W}^{r}(\CO)}}$ and $|\cdot|_{{C^{r}(\CO)}}$ are equivalent.
Given $A>0$, this last norm can be extended to functions defined in a (complex) $A$-neighborhood $\CO+ A:=\{z\in\C : |z-\om|<A,\;\text{for some $\om\in\CO$}\}$ of $\CO$ as well. Following Zehnder \cite{zehnder1}[(2.5), p.109], Cauchy's estimates can be generalized to ``derivatives of non-integer orders" of holomorphic functions $g$ in a neighborhood $\CU$ of $\CO$ such that $|g|_{{C^{r}(\CU)}}<+\infty$~: if $r'<r$ (not necessarily integers) and $ A'< A$, then there exists some constant $C_{r,r'}>1$ such that
\begin{equation}\label{cauchy_esti}
|g|_{{C^{r}(\mathcal{O}+ A)}}\leq \frac{C_{r,r'}}{(A- A')^{r-r'}}|g|_{C^{r'}(\mathcal{O}+ A')}.
\end{equation}

Let us consider the sum $\alpha_{\infty}=\alpha_0+\sum_{\nu\geq 1}(\alpha_{\nu}-\alpha_{\nu-1})$, which converges in $C^{\tilde s}(\CO_\infty)$ if $\tilde s\in\R_+\setminus \N$ with $\tilde s<16s$. Indeed, according to \re{error_alpha-nu}, we apply (\ref{cauchy_esti}) to $\alpha_{\nu+1}-\alpha_{\nu}$, $\nu\in\N$, with $r=\tilde s$, $r'=k_{\tilde s}$. Since $\beta_{\nu+1}<\beta_{\nu+1}^{\tilde s - k_{\tilde s}}$, there is some constant $C_{\tilde s}$ such that
$$
|\alpha_{\nu+1}-\alpha_{\nu}|_{C^{\tilde s}(\mathcal{O}_\infty)}\leq C_{\tilde s}\cdot \frac{\sup_{0\leq l\leq k_{\tilde s}}\|(\alpha_{\nu+1}-\alpha_\nu)^{(l)}\|_{\mathcal{O}_{\nu+1},\beta_{\nu+1}, r_{\nu+1}}}{\beta_{\nu+1}}< C_{\tilde s}\varepsilon_\nu^{\frac{1}{3}-\frac{1}{32s}}.
$$
Since $\frac{1}{3}-\frac{1}{32s}>0$, according to \re{sumeps}, we have
$$
\sum_{l\geq \nu}|\alpha_{l+1}-\alpha_{l}|_{C^{\tilde s}(\CO_\infty)}\leq C_{\tilde s}   \sum_{l\geq \nu}\varepsilon_l^{\frac{1}{3}-\frac{1}{32s}}\leq \frac{C_{\tilde s}  \varepsilon_\nu^{\frac{1}{3}-\frac{1}{32s}}}{1-\varepsilon_\nu^{\frac{32s-3}{8\cdot 96s}}}\rightarrow 0,\quad \nu\to\infty.
$$
Hence, it is a Cauchy sequence in $C^{\tilde s}(\CO_\infty)$. Furthermore, we have
$$|\alpha_{\infty}-\alpha_0|_{C^{\tilde s}_{W}(\CO_\infty)}\leq\sum_{\nu\geq 1}\left|\alpha_{\nu}-\alpha_{\nu-1}\right|_{C^{\tilde s}_{W}(\CO_\infty)}
\leq  \frac{C_{\tilde s}  \varepsilon_0^{\frac{1}{3}-\frac{1}{32s}}}{1-\varepsilon_0^{\frac{32s-3}{8\cdot 96s}}}.$$

With fixed $\om_0\in\CO_\infty$, let us consider the sets
$$
U_{\om_0}:=\left\{\om\in \CO_\infty: \left|\frac{\om}{\om_0}\right|\leq 2\right\},\quad {\cL C}_{\om_0}^{\frac{R}{2}}=\left\{(\xi,\eta)\in (\mathbb{C}^2,0): \xi\eta=\om_0, \  |\xi|, \, |\eta|<  \frac{R}{2}\right\}.
$$
Let us define, for all $\om \in  U_{\om_0}$, the map
$\kappa_{\om}: {\cL C}_{\om_0}^{\frac{R}{2}}\rightarrow {\cL C}_{\om}^{R}$ to be $\kappa_{\om}(\xi,\eta):=(\frac{\om}{\om_0}\xi,\eta)$.
We then define $\widetilde\Psi_{\om,\nu}:= \Psi_{\om,\nu}\circ \kappa_{\om}$ and $\widetilde\Psi_{\om}:= \Psi_{\om}\circ \kappa_{\om}$ on ${\cL C}_{\om_0}^{\frac{R}{2}}$. With the same argument, using Lemma \ref{lem_Psi_nu}, we have
\begin{eqnarray*}
	\sup_{(\xi,\eta)\in {\cL C}_{\om_0}^{\frac{R}{2}}}|\widetilde\Psi_{\om}(\xi,\eta)-\widetilde\Psi_{\om,0}(\xi,\eta)|_{C_{W}^{\tilde s}(U_{\om_0})}&=&\sup_{(\xi,\eta)\in {\cL C}_{\om_0}^{\frac{R}{2}}}\left|\sum_{\nu\geq 1}(\widetilde\Psi_{\om,\nu}(\xi,\eta)-\widetilde\Psi_{\om,\nu-1}(\xi,\eta))\right|_{C_{W}^{\tilde s}(U_{\om_0})}\\
	&\leq& C(\tilde s,\om_0) \sum_{\nu\geq 1}\frac{\|\Psi_{\nu}-\Psi_{\nu-1}\|_{\CO_{\nu+1},\beta_{\nu+1},r_{\nu+1}}}{\beta_{\nu+1}}\\
	&\leq & C(\tilde s,\om_0) \sum_{\nu\geq 1}\varepsilon_{\nu}^{\frac{4}{5}-\frac{1}{32s}}\leq\frac{ C(\tilde s,\om_0) \varepsilon_0^{\frac{4}{5}-\frac{1}{32s}}}{1-\varepsilon_0^{\frac{128s-5}{8\cdot160s}}}.
\end{eqnarray*}
Here, $C(\tilde s,\om_0)>0$ is some constant that depends only on $\tilde s$ and $\om_0$. Hence, $\tilde\Psi_\om$ is Whitney smooth in $\om$ and holomorphic on ${\cL C}_{\om_0}^{\frac{R}{2}}$.

As a consequence, there is a $C^{\tilde s}$-Whitney smooth family of holomorphic invariant curves of $\tau^o_1$, $\tau^o_2$ and $\sigma_o$, for $16s-1<\tilde s<16s$.

\section{Involutions and reversible maps}\label{sec_rev}

In this section, we describe some properties of a pair of germs of involutions $\tau_1$, $\tau_2$ as in \re{gen-invol} with
$\alpha=\alpha(\xi\eta)\in{\CA}^\R_{\beta,r}(\CO)$ satisfying \re{al_nu_sum} -- \re{Ds+nu_sum} with $\CO_{\nu}=\CO$, together with $p,q\in{\CA}_{\beta,r}(\CO)$ with $\|p\|_{\CO,\beta,r}, \ \|q\|_{\CO,\beta,r}<\frac{\varepsilon}{10}$.
Hence, according to Remark \ref{rem_quantities_nu},  if we take the involutions $(\tau_1,\tau_2)=(\tau_\nu^{(1)},\tau_\nu^{(2)})$ for some $\nu\in\N$, described in Proposition \ref{prop_KAM},  then the above assumptions are satisfied. We also consider the germ of map
$\sigma=\tau_1\circ\tau_2$.

Given $0<r_+<r$, recall that we have defined $r^{(m)}$ and $\tilde r$ in (\ref{r_m}) between $r$ and $r_+$:
\begin{equation}\label{r_m-sec6}
r^{(m)}=r_++\frac{m}{8}(r-r_+),\quad m=0,1,\cdots,8, \qquad  \tilde r=r^{(4)}=\frac{r+r_+}{2}.
\end{equation}
We assume, from now on, that $\varepsilon$ is sufficiently small such that
\begin{equation}\label{varepsilon_small-general}
\left(\frac{\left|\ln\varepsilon\right|}{\left|\ln\left(\frac{7}{8}+\frac{r_+}{8r}\right)\right|}+2\right)\frac{(16s+1)^{16s}  \varepsilon^{\frac{1}{2400 s^2}}}{(r-r_+)r_+}<1.
\end{equation}
It is easy to see that (\ref{varepsilon_nu_small}) holds if we have
$$\varepsilon_\nu=\varepsilon,\quad r_\nu=r,\quad r_{\nu+1}= r_+.$$

\subsection{Properties of $\alpha(\cdot)$}\label{sub_alpha}

Recalling \re{Orbeta} and (\ref{new_def_beta}), we define
$$
\beta\in [\varepsilon^{\frac{1}{40s}},  \varepsilon^{\frac{1}{60s}}], \;\
\beta_+=\beta^{\frac54}\in [\varepsilon^{\frac{1}{32s}}, \varepsilon^{\frac{1}{48s}}],
\;\ \tilde\beta=16\beta^{\frac{5}{4}}, \;\  \CO(r,\beta)={\cL O}\cap ]-r^2+\beta, r^2-\beta[.
$$
The definitions are compatible with  (\ref{sequences}) if we take $$\varepsilon=\varepsilon_\nu,\quad \beta=\beta_\nu, \quad \beta_+=\beta_{\nu+1},\quad \tilde\beta=\tilde\beta_\nu,\qquad \nu\in\N.$$
The smallness of $\varepsilon$ in (\ref{varepsilon_small-general}) implies that of $\beta$, $\beta_+$ and $\tilde\beta$, and we have 
\beq \label{tildebeta}
\beta_+<\tilde\beta<\beta < \beta^{\frac{1}{32}} < 2^{-16}.
\eeq
As it is needed below, we also have $r^{-1}<\beta^{-\frac{1}{32}}$ and
\beq\label{exptbeta}
e^{\frac{9\tilde\beta}{8}}< 1+\frac{9\tilde\beta}{8}\sum_{k\geq 0}\left(\frac{9\tilde\beta}{8}\right)^k= 1+\frac{9\tilde\beta}{8}\frac{1}{1-\frac{9\tilde\beta}{8}}<1+\frac{7\tilde\beta}{6}.
\eeq
Then, according to \re{al_nu_sum} -- \re{Ds+nu_sum}, we have
\begin{lemma}\label{lem_deriv_a} $\left|\alpha'\right|_{\CO(r,\beta)}<\frac{11}{10}$ and
 $\left|\alpha^{(k)}\right|_{\CO(r,\beta)}<\beta^{-\frac1{32}}$ for $2\leq k \leq 16s$.
\end{lemma}


\begin{remark} The non-degeneracy condition (\ref{Dsnu_sum}) is essential through this paper.
Nevertheless, we only need the estimates of Lemma \ref{lem_deriv_a} in the rest of this section.
\end{remark}

\begin{lemma}\label{lem_a0} For every $\omega\in\CO(r,\beta) $, $|\alpha(\xi\eta)-\alpha(\omega)|_{\omega,\tilde\beta}<\frac98\tilde\beta$
and
$$|\alpha^{(j)}(\xi\eta)-\alpha^{(j)}(\omega)|_{\omega,\tilde\beta}<\frac{\beta^{\frac{17}{16}}}2,\quad 1\leq j \leq s.$$
\end{lemma}
\proof
For $(\xi,\eta)\in{\cL C}_{\omega,\tilde\beta}$, we have $|\xi\eta-\omega|<\tilde\beta$.
Developing $\alpha(\cdot)$ around $\omega$,
$$|\alpha(\xi\eta)-\alpha(\omega)|\leq \sum_{k\geq1}\frac{|\alpha^{(k)}(\omega)|}{k!}|\xi\eta-\omega|^k< \sum_{k\geq1}\frac{|\alpha^{(k)}(\omega)|\cdot \tilde\beta^k}{k!}.$$
According to \re{tildebeta}, we have $$\frac{\beta^{-\frac1{32}}\tilde\beta}{2(1-\tilde\beta)}=\frac{16\beta^{\frac54-\frac1{32}}}{2(1-\tilde\beta)} <10\beta^{\frac54-\frac{1}{32}}<\frac{1}{90}.$$
Then, in view of Lemma \ref{lem_deriv_a}, we have
\begin{equation}\label{er1}
\sum_{k=1}^{16s}\frac{|\alpha^{(k)}(\omega)|\cdot \tilde\beta^k}{k!}< \frac{11}{10}\tilde\beta+ \frac{\beta^{-\frac1{32}}}2\sum_{k=2}^{16 s}\tilde\beta^k<\frac{11}{10}\tilde\beta+ \frac{\beta^{-\frac1{32}} \tilde\beta^2}{2(1-\tilde\beta)}< \frac{10}{9}\tilde\beta.
\end{equation}
Since $4\pi+1<2^4$, Cauchy's inequality and (\ref{al_nu_sum}) lead to
\begin{equation}\label{deriv_high}
 |\alpha^{(k)}(\omega)|\leq k!\sup_{|z-w|=\frac{\beta}{2}}|\alpha(z)| \cdot\frac{2^{k}}{\beta^{k}}<\frac{k!\cdot 2^{k+4}}{\beta^{k}},\quad k\geq 16s+1.
\end{equation}
Then, we obtain
$$
\sum_{k\geq 16s+1}\frac{|\alpha^{(k)}(\omega)|\cdot\tilde\beta^k}{k!} < \sum_{k\geq 16s+1}\frac{2^{k+4}\tilde\beta^k}{\beta^{k}}
   = 2^4 \sum_{k\geq 16s+1}  \left(2^5 \beta^{\frac{1}{4}}\right)^k
   =\frac{ 2^4 \left(2^5 \beta^{\frac{1}{4}}\right)^{16s+1}}{1-2^5 \beta^{\frac{1}{4}}}.
$$
According to \re{tildebeta}, we have $\frac{1}{1-\beta^{\frac18}}<\frac43$ and
$$2^4\beta^{\frac{16s+1}8}=2^4\cdot \beta^{\frac18} \cdot \beta^{2s}< 2^{-12}\beta^2<\beta^{2}.$$ Therefore, we have
\begin{equation}\label{er2}
\sum_{k\geq 16s+1}\frac{|\alpha^{(k)}(\omega)|\cdot\tilde\beta^k}{k!}<\frac{ 2^4 \left(2^5 \beta^{\frac{1}{4}}\right)^{16s+1}}{1-2^5 \beta^{\frac{1}{4}}}<\frac{2^4\beta^{\frac{16s+1}8}}{1-\beta^{\frac18}}<\frac{4\beta^{2}}{3}<\frac{\tilde\beta}{72}.
\end{equation}
Adding estimates (\ref{er1}) and (\ref{er2}), we obtain $|\alpha(\xi\eta)-\alpha(\omega)|_{\omega,\tilde\beta}<\frac98\tilde\beta$.

\smallskip

For $(\xi,\eta)\in{\cL C}_{\omega,\tilde\beta}$, and $1\leq j \leq s$, by Lemma \ref{lem_deriv_a}, we have,
$$|\alpha^{(j)}(\xi\eta)-\alpha^{(j)}(\omega)|\leq \sum_{k\geq 1}\frac{|\alpha^{(j+k)}(\omega)|}{k!}|\xi\eta-\omega|^k
<\frac{\beta^{-\frac1{32}}}2\sum_{k=1}^{16s-j}\tilde\beta^{k}+\sum_{k\geq 16s+1-j} \frac{(j+k)!  \cdot  2^{j+k+4}\tilde\beta^{k}}{k! \cdot  \beta^{j+k}}.$$
Note that under (\ref{varepsilon_small-general}), for $k\geq 16s+1-j\geq 15s$,
$$\frac{(j+k)! \cdot 2^{j+k+4}\tilde\beta^{k} }{k! \cdot \beta^{j+k}}< 2^{j+5k+4}\beta^{\frac{k}{6}-j} \beta^{\frac{k}{15}}\leq 2^{s+4+5k} \beta^{\frac{k}{24}}\cdot\beta^{\frac{k}{8}-j}\beta^{\frac{k}{15}} <\frac12 \beta^{\frac{k}{15}}.$$
Hence, we have
$$
|\alpha^{(j)}(\xi\eta)-\alpha^{(j)}(\omega)|\leq \frac{\beta}{2}\left(\frac{16 \beta^{\frac{7}{32}}}{1-\tilde\beta}+\frac{\beta^{\frac{16(s-1)+2-j}{15}}}{1-\beta^{\frac{1}{15}}}\right)<\frac{\beta^{\frac{17}{16}}}2.\qed
$$

By Lemma \ref{lem_deriv_a}, Lemma \ref{lem_a0} and (\ref{deriv_high}), we have
\begin{cor}\label{cor_deri_a} $\sup_{\omega\in\CO(r,\beta)}|\alpha'(\xi\eta)|_{\omega,\tilde\beta}<\frac65$,
$$ \sup_{\omega\in\CO(r,\beta)}|\alpha^{(k)}(\xi\eta)|_{\omega,\tilde\beta}<
\left\{\begin{array}{cl}
 2\beta^{-\frac1{32}} &2\leq k \leq s \ if \ s\geq 2\\[2mm]
 \frac{ k! 2^{k+5}}{\beta^{k}}, & k\geq s+1,
\end{array}\right. .$$
\end{cor}

\noindent

\begin{cor}\label{lem_eiba}
For any $b\in\R$, and any $0\leq \beta'\leq\tilde\beta$, we have
$$\sup_{\omega\in\CO(r,\beta)}|e^{{\rm i}b\alpha(\xi\eta)}|_{\omega,\beta'}<e^{\frac98|b|\tilde\beta}.$$
Moreover, for $-1\leq b\leq 1$,
$$\sup_{\omega\in\CO(r,\beta)}\left||e^{{\rm i}b\alpha(\xi\eta)}|-1\right|_{\omega,\beta'}<\frac54 \tilde\beta.$$
\end{cor}

\begin{remark}\label{rem_eiba} By Corollary \ref{lem_eiba}, for $b\in\R$ with $|b|\leq\beta^{-\frac14}$, we have
$\sup_{\omega\in\CO(r,\beta)} |e^{{\rm i}b\alpha(\xi\eta)}|_{\omega,\beta'}<e^{\frac98\beta^{-\frac14}\cdot\tilde\beta}=e^{18 \beta}$, for all $0\leq \beta'\leq\tilde\beta$.
Since $\beta<2^{-16}$, we usually use the rough estimate in this paper for convenience:
$$\sup_{\omega\in\CO(r,\beta)} |e^{{\rm i}b\alpha(\xi\eta)}|_{\omega,\beta'}<\frac{101}{100}, \quad  \forall \ b\in\R \ with \  |b|\leq \beta^{-\frac14}.$$
\end{remark}

\noindent
{\it Proof of Corollary \ref{lem_eiba}.} Since $\alpha(\omega)\in\R$, $|e^{{\rm i}b\alpha(\omega)}|=1$. Then $|e^{{\rm i}b\alpha(\xi\eta)}|=|e^{{\rm i}b(\alpha(\xi\eta)-\alpha(\omega))}|$.
For $(\xi,\eta)\in {\cL C}_{\omega,\beta'}$, we have
$$e^{{\rm i}b(\alpha(\xi\eta)-\alpha(\omega))}=1+\sum_{k\geq 1}\frac{{\rm i}^k b^k}{k!}(\alpha(\xi\eta)-\alpha(\omega))^k.$$
If $|b|\leq 1$, then, under (\ref{varepsilon_small-general}), $|b||\alpha(\xi\eta)-\alpha(\omega)|< \frac98\tilde\beta$ is sufficiently small, and we have
$$ \left||e^{{\rm i}b\alpha(\xi\eta)}|-1\right| = \left||e^{{\rm i}b(\alpha(\xi\eta)-\alpha(\omega))}|-1\right|
\leq \left|\sum_{k\geq 1}\frac{{\rm i}^k b^k}{k!}(\alpha(\xi\eta)-\alpha(\omega))^k\right|\leq \frac{\frac98\tilde\beta}{1- \frac98\tilde\beta}
<\frac54\tilde\beta.$$
Moreover, for any $b\in\R$, any $0\leq \beta'\leq\tilde\beta$,
$$
\left|e^{{\rm i}b\alpha(\xi\eta)}\right|_{\omega,\beta'}\leq\left|\sum_{k\geq 0}\frac{{\rm i}^k b^k}{k!}(\alpha(\xi\eta)-\alpha(\omega))^k\right|_{\omega,\tilde\beta}
<\sum_{k\geq 0}\frac{|b|^k}{k!}\left(\frac98\tilde\beta\right)^k
=e^{\frac98|b|\tilde\beta}.\qed
$$

\begin{lemma}\label{lem_a1} Given $0<r'<r$, assume $\beta$ is sufficiently small such that
\begin{equation}\label{beta_small01}
\beta<r^2-r'^2.
\end{equation}
 Given $0\leq \beta'\leq \tilde\beta$, and $f\in{\CA}_{\beta', r'}({\CO})$ with
\begin{equation}\label{f_small}
\|f\|_{{\CO}, \beta',r'}<\beta^{24s},
\end{equation}
then $\alpha(\xi\eta+f)\in{\CA}_{\beta', r'}({\CO})$ with
 \begin{equation}\label{uniform_alpha}
 \|\alpha(\xi\eta+f)-\alpha(\xi\eta)\|_{{\CO},\beta', r'} <\frac{5}4\|f\|_{{\CO},\beta',r'},
 \end{equation}
and for $-1\leq b\leq 1$, $e^{{\rm i}b\alpha(\xi\eta+f)}\in{\CA}_{\beta', r'}({\CO})$ with
\begin{equation}\label{uniform_eialpha}
\|e^{{\rm i}b\alpha(\xi\eta+f)}-e^{{\rm i}b\alpha(\xi\eta)}\|_{{\CO},\beta', r'}<\frac43\|f\|_{{\CO},\beta',r'}.
\end{equation}
\end{lemma}
\begin{remark} Recalling $r^{(m)}$ given in (\ref{r_m-sec6}), we see that, under the assumption (\ref{varepsilon_small-general}), (\ref{beta_small01}) is satisfied for
$$
r=r^{(m+1)},\quad r'=r^{(m)},\qquad m=0,1,\cdots,7.
$$
Indeed, according to (\ref{varepsilon_small-general}), we have
$$\beta\leq \varepsilon^{\frac{1}{60s}}<\frac{(r-r_+)r_+}{(16s+1)^{16s}}
<\frac{(r-r_+)(15r_++r)}{8}=(r^{(1)})^2-(r^{(0)})^2\leq (r^{(m+1)})^2-(r^{(m)})^2.$$
\end{remark}

\proof Let $\varsigma:=\|f\|_{{\CO},\beta',r'}$.
\re{f_small} implies that, for $(\xi,\eta)\in {\cL C}_{\omega,\beta'}$, $\omega\in {\CO}(r',\beta')$,
$$|\xi\eta+f(\xi,\eta)-\omega|\leq |\xi\eta-\omega|+\|f\|_{{\CO},\beta',r'}\leq\tilde\beta+\varsigma<\beta,$$
and (\ref{beta_small01}) implies that $r'^2-\beta'<r'^2<r^2-\beta$.

Developing $\alpha(\cdot)$ around $\xi\eta$, we obtain
\begin{equation}\label{develop_alpha}
\alpha(\xi\eta+f)-\alpha(\xi\eta)
=\sum_{k\geq 1}\frac{ \alpha^{(k)}(\xi\eta)}{k!}f^k.
\end{equation}
By Lemma \ref{lem_norm00}, for every $\omega\in{\CO}(r',\beta')$,
we have $\|f^k\|_{\omega,\beta', r'}\leq \varsigma^k$ for $k\in\N^*$. Then, in view of Corollary \ref{cor_deri_a}, we have
$$\|\alpha(\xi\eta+f)-\alpha(\xi\eta)\|_{\omega,\beta', r'}
<\left\{\begin{array}{lc}
\frac65\varsigma + 2\beta^{-\frac1{32}} \sum_{k=2}^s \frac{\varsigma^k}{k!}+ 2^5\sum_{k\geq s+1} \frac{2^{k}\varsigma^k}{\beta^{k}}, & s\geq 2 \\[2mm]
\frac65\varsigma + 2^5\sum_{k\geq 2} \frac{2^{k}\varsigma^k}{\beta^{k}}, & s=1
\end{array}
\right.,$$
which can be bounded by $\frac{5}{4}\varsigma$ under assumption $\varsigma<\beta^{24s}$. Indeed,
\begin{itemize}
  \item if $s=1$, we have $\varsigma<\beta^{24s}=\beta^{24}$, then
  $$2^5 \sum_{k\geq 2}\frac{2^k \varsigma^k}{\beta^k}< 2^5 \sum_{k\geq 2} \frac{2^k \varsigma^k}{\varsigma^{\frac{k}{24}}}=2^5 \sum_{k\geq 2}\left(2\varsigma^{\frac{23}{24}}\right)^k=\frac{2^7\varsigma^{\frac{23}{12}}}{1-2\varsigma^{\frac{23}{24}}}.$$
  \item if $s\geq 2$, then, $\varsigma<\beta^{24s}$ implies that
  $$2^5\sum_{k\geq s+1}\frac{2^k \varsigma^k}{\beta^k}< 2^5 \sum_{k\geq s+1} \frac{2^k \varsigma^k}{\varsigma^{\frac{k}{24s}}}<2^5 \sum_{k\geq 2}\left(2\varsigma^{1-\frac{1}{24s}}\right)^k=\frac{2^7\varsigma^{2-\frac{1}{12s}}}{1-2\varsigma^{1-\frac{1}{24s}}},$$
$$2\beta^{-\frac1{32}}\sum_{k=2}^s\frac{\varsigma^k}{k!}< \beta^{-\frac{1}{32}}\sum_{k=2}^s\varsigma^k<\frac{\varsigma^{-\frac{1}{24s\cdot 32}} \cdot\varsigma^2}{1-\varsigma}=\varsigma\cdot\frac{\varsigma^{1-\frac{1}{24s\cdot 32}}}{1-\varsigma}.$$
\end{itemize}
According to \re{tildebeta}, we have
$$
\varsigma<\varsigma^{1-\frac{1}{32\cdot24s}}<\varsigma^{1-\frac{1}{24s}}<\varsigma^{1-\frac{1}{12s}}<\beta^{22}<2^{-16\cdot32\cdot22}.
$$
Hence, we obtain these rough estimates
$$
\frac{2^7\varsigma^{1-\frac{1}{12s}}}{1-2\varsigma^{1-\frac{1}{24s}}}<\frac{1}{200},\quad \frac{\varsigma^{1-\frac{1}{24s\cdot 32}}}{1-\varsigma}<\frac{1}{200},\quad \forall \ s\geq 1.
$$
As a consequence, we have, for every $\omega\in{\CO}(r',\beta')$, $\|\alpha(\xi\eta+f)-\alpha(\xi\eta)\|_{\omega,\beta', r'}\leq \frac65\varsigma+ \frac{\varsigma}{100}\leq \frac54\varsigma$.
By \rl{lem_norm00} and \ref{lemma_analytic_prod}, we see that $f^k\in {\CA}_{\beta', r'}({\CO})$ for every $k\in\N$, and, in view of (\ref{esi_coeff_uniform}),
$$|(f^k)_{l,j}|_{{\CO}(r',\beta')}\leq \|f\|_{{\CO},\beta',r'}^k r^{-(l+j)}=\varsigma^k r^{-(l+j)}.$$
Hence, according to (\ref{develop_alpha}), $\alpha(\xi\eta+f)-\alpha(\xi\eta)\in {\CA}_{\beta', r'}({\CO})$ with
$$\left(\alpha(\xi\eta+f)-\alpha(\xi\eta)\right)_{l,j}(\omega)=\sum_{k\geq 1}\frac{ \alpha^{(k)}(\omega)}{k!}(f^k)_{l,j}(\omega),\quad \forall \ l,j\geq 0 .$$

For $-1\leq b \leq 1$, we have
$$e^{{\rm i}b\alpha(\xi\eta+f)}-e^{{\rm i}b\alpha(\xi\eta)} = e^{{\rm i}b\alpha(\xi\eta)}\left(e^{{\rm i}b(\alpha(\xi\eta+f)-\alpha(\xi\eta))}-1 \right)
=e^{{\rm i}b\alpha(\xi\eta)}\sum_{k\geq 1} \frac{{\rm i}^kb^k}{ k!}\left(\alpha(\xi\eta+f)-\alpha(\xi\eta) \right)^k.$$
Then, by \rl{lem_norm00} and Remark \ref{rem_eiba}, we obtain, for every $\omega\in{\CO}(r',\beta')$,
$$\|e^{{\rm i}b\alpha(\xi\eta+f)}-e^{{\rm i}b\alpha(\xi\eta)}\|_{\omega,\beta', r'}\leq|e^{{\rm i}b\alpha(\xi\eta)}|_{\omega,\beta'}\sum_{k\geq 1} \frac{1}{ k!}\|\alpha(\xi\eta+f)-\alpha(\xi\eta)\|^k_{\omega,\beta', r'}
<\frac{101}{100}\sum_{k\geq 1}\frac{5^k \varsigma^k}{k!   4^k}<\frac43\varsigma,$$
which implies (\ref{uniform_eialpha}). By Lemma \ref{lemma_analytic_prod} and \ref{lemma_analytic_exp}, we see that
$$e^{{\rm i}b\alpha(\xi\eta+f)}-e^{{\rm i}b\alpha(\xi\eta)}=\left(e^{{\rm i}b(\alpha(\xi\eta+f)-\alpha(\xi\eta))}-1\right)e^{{\rm i}b\alpha(\xi\eta)}\in {\CA}_{\beta', r'}({\CO}).\qed$$

\begin{lemma}\label{lem_eia} Given $0<r''<r'\leq r<\frac{1}{4} $,
if $\beta$ is sufficiently small such that
\begin{equation}\label{small_beta02}
8\beta^{\frac12}<(r'-r'')r'', \quad e^{\frac98\beta}\frac{r''}{r'}<1-\frac{\beta^2}{16},
\end{equation}
then for $-1\leq b\leq 1$, for $0<\beta'\leq \tilde\beta$, $h\in {\CA}_{\beta',r'}({\CO})$ with $\|h\|_{{\CO},\beta',r'}< +\infty$, we have $h(e^{{\rm i}b\alpha(\xi\eta)}\xi,e^{-{\rm i}b\alpha(\xi\eta)}\eta)\in{\CA}_{\beta'',r''}({\CO})$ with
$$
\|h(e^{{\rm i}b\alpha(\xi\eta)}\xi,e^{-{\rm i}b\alpha(\xi\eta)}\eta)\|_{\CO,\beta'',r''}<  \|h\|_{\CO,\beta',r'}.
$$
\end{lemma}

\begin{remark}\label{rem_small_beta02} Under the assumption (\ref{varepsilon_small-general}), (\ref{small_beta02}) is satisfied for
\begin{equation}\label{case_small_beta02}
r'=r^{(m+1)},\quad r''=r^{(m)},\qquad m=0,1,\cdots,7.
\end{equation}
Indeed, (\ref{varepsilon_small-general}) implies that
\begin{equation}\label{1/(r1-r2)}
\frac17>-\ln\left(\frac{7}{8}+\frac{r_+}{8r}\right)=\left|\ln\left(\frac{7}{8}+\frac{r_+}{8r}\right)\right|>\frac{(16s+1)^{16s} |\ln\varepsilon| \cdot \varepsilon^{\frac{1}{2400s^2}}}{(r-r_+)r_+}>2\beta^{\frac{1}{32}}>\frac54\beta.
\end{equation}
Then we have
$$8\beta^{\frac12}<\frac{(r-r_+)r_+}{7(16s+1)^{16s}}\leq \frac{8(r^{(m+1)}-r^{(m)})r^{(m)}}{7(16s+1)^{16s}}<(r^{(m+1)}-r^{(m)})r^{(m)},$$
$$e^{\frac98\beta}\frac{r^{(7)}}{r^{(8)}}=e^{\frac98\beta}\left( \frac{7}{8}+\frac{r_+}{8r}\right)<e^{-\frac{\beta}{8}}<1-\frac{\beta^2}{16}.$$
Hence we obtain (\ref{small_beta02}) for the case (\ref{case_small_beta02}) by noting that, for $0\leq m\leq 6$,
$$\frac{r^{(m)}}{r^{(m+1)}}=\frac{8 r_++m(r-r_+)}{8 r_++(m+1)(r-r_+)}<\frac{8 r_++(m+1)(r-r_+)}{8 r_++(m+2)(r-r_+)}=\frac{r^{(m+1)}}{r^{(m+2)}}.$$
\end{remark}

\proof In view of Corollary \ref{lem_eiba}, we have, for every $\omega\in{\CO}(r'',\beta'')$,
\begin{eqnarray*}
 \left\|h(e^{{\rm i}b\alpha}\xi,e^{-{\rm i}b\alpha}\eta)\right\|_{\omega,\beta'',r''}
&\leq&
\left| h_{0,0}\right|_{\omega,\beta'}+\sum_{l\geq 1}  \left|e^{{\rm i}b l\alpha} h_{l,0}\right|_{\omega,\beta'} r''^{l}+\sum_{j\geq 1}  \left|e^{-{\rm i}b j\alpha} h_{0,j}\right|_{\omega,\beta'} r''^{j}\\
&\leq&\left| h_{0,0}\right|_{\omega,\beta'}+\sum_{l\geq 1}e^{\frac98 l|b|\tilde\beta} \left| h_{l,0}\right|_{\omega,\beta'}r''^{l}+\sum_{j\geq 1}e^{\frac98 j|b|\tilde\beta} \left| h_{0,j}\right|_{\omega,\beta'}r''^{j}\\
&\leq&\left|h_{0,0}\right|_{\omega,\beta'}+\sum_{l\geq 1}\left(e^{\frac98\tilde\beta}\frac{r''}{r'}\right)^{l} \left| h_{l,0}\right|_{\omega,\beta'}r'^{l}+\sum_{j\geq 1}\left(e^{\frac98\tilde\beta}\frac{r''}{r'}\right)^{j} \left| h_{0,j}\right|_{\omega,\beta'}r'^{j}\\
&<&\left\|h\right\|_{\omega,\beta',r'}.
\end{eqnarray*}
Noting that, for $l,j\geq 0$,
$$\left(h(e^{{\rm i}b\alpha(\xi\eta)}\xi,e^{-{\rm i}b\alpha(\xi\eta)}\eta) \right)_{l,j}(\omega)=h_{l,j}(\omega) e^{{\rm i}b(l-j)\alpha(\omega)},$$
we see that $h\left(e^{{\rm i}b\alpha(\xi\eta)}\xi,e^{-{\rm i}b\alpha(\xi\eta)}\eta\right)\in{\CA}_{\beta'',r''}({\CO})$.
\qed

\begin{lemma}\label{lem_norm}
Let $0<r''<r'\leq r<\frac{1}{4}$ and $0<2\beta''\leq\beta'\leq \tilde\beta$. If $\beta$ is small enough such that (\ref{small_beta02}) is satisfied, then for
$h \in {\CA}_{\beta',r'}\left({\CO}\right)$ with $\|h\|_{\CO,\beta',r'}< +\infty$, for $f_1,f_2,g_1,g_2$ with
\begin{equation}\label{f_g_m}
\|f_m\|_{\CO,\beta'',r''},  \|g_m\|_{\CO,\beta'',r''} < \frac{\beta'^2}{16} , \quad m=1,2,
\end{equation}
we have that, for $-1\leq b\leq 1$,
\begin{eqnarray*}
& &\|h(e^{{\rm i}b\alpha(\xi\eta)}\xi+f_1,e^{-{\rm i}b\alpha(\xi\eta)}\eta+g_1)-h(e^{{\rm i}b\alpha(\xi\eta)}\xi+f_2,e^{-{\rm i}b\alpha(\xi\eta)}\eta+g_2)\|_{{\CO},\beta'',r''}\nonumber\\
&<&\frac{3r'\|h\|_{\CO,\beta',r'}}{(r'-r'')\beta'}\max\left\{ \|f_1-f_2\|_{\CO,\beta'',r''}, \|g_1-g_2\|_{\CO,\beta'',r''}\right\}.
\end{eqnarray*}
Moreover, if $f_1,f_2,g_1,g_2\in {\CA}_{\beta'',r''}\left({\CO}\right)$, then
$$h(e^{{\rm i}b\alpha(\xi\eta)}\xi+f_1,e^{-{\rm i}b\alpha(\xi\eta)}\eta+g_1)-h(e^{{\rm i}b\alpha(\xi\eta)}\xi+f_2,e^{-{\rm i}b\alpha(\xi\eta)}\eta+g_2)
 \in   {\CA}_{\beta'',r''}\left({\CO}\right).$$
\end{lemma}

\begin{remark} We deduce Lemma \ref{lem_norm-pre} from Lemma \ref{lem_norm} by taking $b=0$, with (\ref{smallness_beta}) verified by (\ref{small_beta02}).
\end{remark}
We postpone the detailed proof of Lemma \ref{lem_norm} in Appendix \ref{proof_lemma_norm}.

\smallskip

\subsection{Properties of the perturbation}\label{sub_pert}

Recall that $\tau_1$ and $\tau_2$ are given by \re{gen-invol} with  $p,q\in{\CA}_{\beta,r}(\CO)$ satisfying \re{esti_pq}.

\begin{lemma}\label{lemma_real_res}
For $0\leq k \leq 16s$,
\begin{equation}\label{e_p_01}
|(e^{-\frac{\rm i}2\alpha}p_{0,1})^{(k)}|_{{\cL O}(r^{(7)},\tilde\beta)} \, ,\  |(e^{\frac{\rm i}2\alpha}\bar p_{0,1})^{(k)}|_{{\cL O}(r^{(7)},\tilde\beta)}<\frac{\varepsilon^\frac13}{10}.
\end{equation}
\end{lemma}
\proof By (\ref{esi_coeff}), Corollary \ref{lem_eiba} and the assumption $ \|p\|_{\CO,\beta,r}<\frac{\varepsilon}{10}$, we see that
$$\sup_{|z-\omega|<\beta} |(e^{-\frac{\rm i}2\alpha}p_{0,1})(z)|< \frac{101\varepsilon}{1000 r},\quad \omega\in {\cL O}(r^{(7)},\tilde\beta).$$
Hence \re{e_p_01} is true for $k=0$.
According to (\ref{varepsilon_small-general}) and \re{exptbeta}, we have
$$  \frac{(16s)! e^{\frac9{16}\tilde\beta}\varepsilon^{\frac12}}{10r}<\varepsilon^{\frac12-\frac1{1920}}\cdot\frac{101(16s+1)^{16s}\varepsilon^{\frac1{1920s}}}{1000r_+}<\varepsilon^{\frac12-\frac1{1920}}< \frac{\varepsilon^\frac13}{10}.
$$
Indeed, according to (\ref{varepsilon_small-general}), we have
$$\frac{\varepsilon^{\frac{1}{1920 s}}}{r_+}< \left|\ln\left(\frac{7}{8}+\frac{r_+}{8r}\right)\right|(16s+1)^{-16s}(r-r_+)<\frac{r-r_+}{7(16s+1)^{16s}}<1.$$
Hence, $\varepsilon^{\frac16-\frac{1}{1920}}<\varepsilon^{\frac{1}{1920 s}}<r_+.$
Then, applying Cauchy's inequality and recalling that $\tilde\beta\in [16\varepsilon^{\frac{1}{32s}},16\varepsilon^{\frac{1}{48s}}]$, we have, for $1\leq k \leq 16s$, $\omega\in {\cL O}(r^{(7)},\tilde\beta)$,
$$|(e^{-\frac{\rm i}2\alpha}p_{0,1})^{(k)}(\omega)| \leq  k!\cdot\frac{\sup_{|z-\omega|=\tilde\beta}|(e^{-\frac{\rm i}2\alpha}p_{0,1})(z)|}{16^k\varepsilon^{\frac{k}{32s}}}
 < \frac{(16s)! e^{\frac9{16}\tilde\beta}\varepsilon^{\frac12}}{10r}
 < \frac{\varepsilon^\frac13}{10}.$$
The proof for $e^{\frac{\rm i}2\alpha}\bar p_{0,1}$ is similar.\qed

\smallskip

For $\sigma=\tau_1\circ\tau_2$, we have
$$
\sigma(\xi,\eta)=\left(\begin{array}{l}
\exp\left\{{\frac{\rm i}2\left(\alpha(\xi\eta+e^{-\frac{\rm i}2\alpha}\eta \bar q+e^{\frac{\rm i}2\alpha}\xi \bar p+\bar p\bar q)+\alpha(\xi\eta)\right)}\right\}\xi\\[2mm]
+ \, \exp\left\{{\frac{\rm i}2\alpha(\xi\eta+e^{-\frac{\rm i}2\alpha}\eta \bar q+e^{\frac{\rm i}2\alpha}\xi \bar p+\bar p\bar q)}\right\}\bar q(\xi,\eta)\\[2mm]
+ \,  p(e^{-\frac{\rm i}2\alpha}\eta+\bar p,e^{\frac{\rm i}2\alpha}\xi+\bar q)\\[3mm]
\exp\left\{-{\frac{\rm i}2\left(\alpha(\xi\eta+e^{-\frac{\rm i}2\alpha}\eta \bar q+e^{\frac{\rm i}2\alpha}\xi \bar p+\bar p\bar q)+\alpha(\xi\eta)\right)}\right\}\eta\\[2mm]
+ \,  \exp\left\{-{\frac{\rm i}2\alpha(\xi\eta+e^{-\frac{\rm i}2\alpha}\eta \bar q+e^{\frac{\rm i}2\alpha}\xi \bar p+\bar p\bar q)}\right\}\bar p(\xi,\eta)\\[2mm]
+ \,  q(e^{-\frac{\rm i}2\alpha}\eta+\bar p,e^{\frac{\rm i}2\alpha}\xi+\bar q)
\end{array}
\right)=:\left(\begin{array}{l}
	e^{\frac{\rm i}2\alpha(\xi\eta)}\xi+ f(\xi,\eta)\\
	e^{-\frac{\rm i}2\alpha(\xi\eta)}\eta+ g(\xi,\eta)
\end{array}\right).
$$
Hence, we have
\begin{align}
  f =&\left(\exp\left\{\frac{\rm i}2\left(\alpha(\xi\eta+e^{-\frac{\rm i}2\alpha}\eta \bar q+e^{\frac{\rm i}2\alpha}\xi \bar p+\bar p\bar q)-\alpha(\xi\eta)\right)\right\}-1 \right) e^{{\rm i}\alpha}\xi \label{f}\\
   & + \exp\left\{{\frac{\rm i}2\alpha(\xi\eta+e^{-\frac{\rm i}2\alpha}\eta \bar q+e^{\frac{\rm i}2\alpha}\xi \bar p+\bar p\bar q)}\right\}\bar q +   p(e^{-\frac{\rm i}2\alpha}\eta+\bar p,e^{\frac{\rm i}2\alpha}\xi+\bar q),\nonumber\\[2mm]
  g =& \left( \exp\left\{-\frac{\rm i}2\left(\alpha(\xi\eta+e^{-\frac{\rm i}2\alpha}\eta \bar q+e^{\frac{\rm i}2\alpha}\xi \bar p+\bar p\bar q)-\alpha(\xi\eta)\right)\right\}-1\right) e^{-{\rm i}\alpha}\eta \label{g}\\
   & + \exp\left\{-{\frac{\rm i}2\alpha(\xi\eta+e^{-\frac{\rm i}2\alpha}\eta \bar q+e^{\frac{\rm i}2\alpha}\xi \bar p+\bar p\bar q)}\right\}\bar p +  q(e^{-\frac{\rm i}2\alpha}\eta+\bar p,e^{\frac{\rm i}2\alpha}\xi+\bar q).\nonumber
\end{align}

\begin{lemma}\label{lem_appro}
$f,g\in{\CA}_{\tilde\beta, r^{(7)}}({\cL O})$ with
\begin{align}
\left\|f-\frac{{\rm i}\alpha'(\xi\eta)}2 (e^{-\frac{\rm i}2\alpha}\eta \bar q+e^{\frac{\rm i}2\alpha}\xi \bar p)e^{{\rm i}\alpha}\xi-e^{\frac{\rm i}2\alpha}\bar q-p(e^{-\frac{\rm i}2\alpha}\eta,e^{\frac{\rm i}2\alpha}\xi)\right\|_{\CO,\tilde\beta,r^{(7)}}&<\frac{\varepsilon^{\frac{31}{16}}}{80},\label{appro_f}\\
\left\|g+\frac{{\rm i}\alpha'(\xi\eta)}2 (e^{-\frac{\rm i}2\alpha}\eta \bar q+e^{\frac{\rm i}2\alpha}\xi \bar p)e^{-{\rm i}\alpha}\eta-e^{-\frac{\rm i}2\alpha}\bar p-q(e^{-\frac{\rm i}2\alpha}\eta,e^{\frac{\rm i}2\alpha}\xi)\right\|_{\CO,\tilde\beta,r^{(7)}}&<\frac{\varepsilon^{\frac{31}{16}}}{80}.\label{appro_g}
\end{align}
\end{lemma}
\proof According to Lemma \ref{lem_a1} and \ref{lem_norm}, we see that $p(e^{-\frac{\rm i}2\alpha}\eta+\bar p,e^{\frac{\rm i}2\alpha}\xi+\bar q)$, $q(e^{-\frac{\rm i}2\alpha}\eta+\bar p,e^{\frac{\rm i}2\alpha}\xi+\bar q)$ and $\alpha(\xi\eta+e^{-\frac{\rm i}2\alpha}\eta \bar q+e^{\frac{\rm i}2\alpha}\xi \bar p+\bar p\bar q)$ are elements of ${\CA}_{\tilde\beta, r^{(7)}}({\cL O})$.
Then, by Lemma \ref{lemma_analytic_exp},
\begin{align*}
\exp\left\{\pm{\frac{\rm i}2\alpha(\xi\eta+e^{-\frac{\rm i}2\alpha}\eta \bar q+e^{\frac{\rm i}2\alpha}\xi \bar p+\bar p\bar q)}\right\}\in&{\CA}_{\tilde\beta, r^{(7)}}({\cL O}), \\
\exp\left\{\pm\frac{\rm i}2\left(\alpha(\xi\eta+e^{-\frac{\rm i}2\alpha}\eta \bar q+e^{\frac{\rm i}2\alpha}\xi \bar p+\bar p\bar q)-\alpha(\xi\eta)\right)\right\}\in&{\CA}_{\tilde\beta, r^{(7)}}({\cL O}).
\end{align*}
In view of (\ref{f}) and (\ref{g}), combining with Lemma \ref{lemma_analytic_prod}, we obtain that $f,g\in{\CA}_{\tilde\beta, r^{(7)}}({\cL O})$.

For $(\xi,\eta)\in {\cL C}^{r^{(7)}}_{\omega,\tilde\beta}$, $\omega\in {\cL O}(r^{(7)},\tilde\beta)$, we have
\begin{eqnarray}
&&\left(\exp\left\{\frac{\rm i}2(\alpha(\xi\eta+e^{-\frac{\rm i}2\alpha}\eta \bar q+e^{\frac{\rm i}2\alpha}\xi \bar p+\bar p\bar q)-\alpha(\xi\eta))\right\}-1 \right) e^{{\rm i}\alpha}\xi\nonumber\\
&=&e^{{\rm i}\alpha}\xi\sum_{k\geq 1}\frac{{\rm i}^k}{2^k k!}\left(\alpha(\xi\eta+e^{-\frac{\rm i}2\alpha}\eta \bar q+e^{\frac{\rm i}2\alpha}\xi \bar p+\bar p\bar q)-\alpha(\xi\eta)\right)^k\nonumber\\
&=&\frac{\rm i}2\alpha'(\xi\eta) \left( e^{-\frac{\rm i}2\alpha}\eta \bar q+e^{\frac{\rm i}2\alpha}\xi \bar p\right) e^{{\rm i}\alpha}\xi\label{1st-term}\\
& &+\, \frac{\rm i}2  \alpha'(\xi\eta) e^{{\rm i}\alpha}\xi\cdot \bar p\bar q + \frac{{\rm i}e^{{\rm i}\alpha}\xi}{2} \sum_{j\geq 2} \frac{ \alpha^{(j)}(\xi\eta)}{j!}(e^{\frac{\rm i}2\alpha}\eta \bar q+e^{-\frac{\rm i}2\alpha}\xi \bar p+\bar p\bar q)^j\label{H1}\\
& &+ \, e^{{\rm i}\alpha}\xi\sum_{k\geq 2}\frac{{\rm i}^k}{2^k k!}\left(\alpha(\xi\eta+e^{\frac{\rm i}2\alpha}\eta \bar q+e^{-\frac{\rm i}2\alpha}\xi \bar p+\bar p\bar q)-\alpha(\xi\eta)\right)^k.\label{H2}
\end{eqnarray}
Since $\|p\|_{\CO,\beta,r}$, $ \|q\|_{\CO,\beta,r}<\frac{\varepsilon}{10}$, by Corollary \ref{lem_eiba} and \re{exptbeta}, we obtain
$$
\|e^{-\frac{\rm i}2\alpha}\eta \bar q+e^{\frac{\rm i}2\alpha}\xi \bar p+\bar p\bar q\|_{\CO,\tilde\beta,r^{(7)}}<\frac{e^{\frac98\tilde\beta}r^{(7)}\varepsilon}{10}+\frac{\varepsilon^2}{100}<\frac{\varepsilon}{10}.
$$
Hence, applying Corollary \ref{cor_deri_a} and Lemma \ref{lemma_analytic_prod}, we have, if $s\geq 2$, then
\begin{multline*}
	 \left\|\frac{{\rm i}e^{{\rm i}\alpha}\xi}{2} \sum_{j\geq 2} \frac{ \alpha^{(j)}(\xi\eta)}{j!}(e^{-\frac{\rm i}2\alpha}\eta \bar q+e^{\frac{\rm i}2\alpha}\xi \bar p+\bar p\bar q)^j\right\|_{\CO,\tilde\beta,r^{(7)}} \\
	\leq r^{(7)} e^{\frac98\tilde\beta}\cdot\beta^{-\frac1{32}}\sum_{j=2}^{s} \frac{1}{j! }\left(\frac{\varepsilon}{10}\right)^j+  \frac{r^{(7)} e^{\frac98\tilde\beta}}{2}  \sum_{j\geq s+1}  \frac{2^{j+5}\varepsilon^j}{\beta^{j}10^j}
	<\frac{\varepsilon^{\frac{63}{32}}}{100}.
\end{multline*}
Otherwise, for $s=1$, we have
$$ \left\|\frac{{\rm i}e^{{\rm i}\alpha}\xi}{2} \sum_{j\geq 2} \frac{ \alpha^{(j)}(\xi\eta)}{j!}(e^{-\frac{\rm i}2\alpha}\eta \bar q+e^{\frac{\rm i}2\alpha}\xi \bar p+\bar p\bar q)^j\right\|_{\CO,\tilde\beta,r^{(7)}}
 \leq\frac{r^{(7)} e^{\frac98\tilde\beta}}{2} \sum_{j\geq 2}  \frac{2^{j+5}\varepsilon^j}{\beta^{j}10^j}
 <\frac{\varepsilon^{\frac{31}{16}}}{200}.$$
Then, the $\|\cdot\|_{\CO,\tilde\beta,r^{(7)}}-$norm of terms in (\ref{H1}) is bounded by
\begin{equation}\label{est_f1}
\frac{\varepsilon^{\frac{31}{16}}}{200}+ \frac{3 r^{(7)} e^{\frac98\tilde\beta}}{5} \left(\frac{\varepsilon}{10}\right)^2<\frac{ \varepsilon^{\frac{31}{16}}}{160}.
\end{equation}
Applying Lemma \ref{lem_a1}, we obtain $$ \|\alpha(\xi\eta+e^{-\frac{\rm i}2\alpha}\eta \bar q+e^{\frac{\rm i}2\alpha}\xi \bar p+\bar p\bar q)-\alpha(\xi\eta)\|_{\CO,\tilde\beta,r^{(7)}}<\frac{\varepsilon}{8}.$$ Hence the $\|\cdot\|_{\CO,\tilde\beta,r^{(7)}}-$norm of terms (\ref{H2}) is bounded by
\begin{equation}\label{est_f2}
r^{(7)} e^{{\frac98}\tilde\beta} \sum_{k\geq 2}\frac{\varepsilon^k}{k! 2^k 8^k}<\frac{\varepsilon^{2}}{100}.
\end{equation}
Since $\|p\|_{\CO,\tilde\beta,r^{(7)}}\leq\|p\|_{\CO,\beta,r}<\frac{\varepsilon}{10}$, by Lemma \ref{lem_norm}, we have
\begin{multline}\label{est_f3}
\|p(e^{-\frac{\rm i}2\alpha}\eta+\bar p,e^{\frac{\rm i}2\alpha}\xi+\bar q)-p(e^{-\frac{\rm i}2\alpha}\eta,e^{\frac{\rm i}2\alpha}\xi)\|_{\CO,\tilde\beta,r^{(7)}} \\
<\frac{3r}{(r-r^{(7)})\tilde\beta} \max\{\|p\|_{\CO,\tilde\beta,r^{(7)}}, \|q\|_{\CO,\tilde\beta,r^{(7)}}\}\cdot\|p\|_{\CO,\beta,r}\leq \frac{3r \varepsilon^2}{100\cdot 16(r-r^{(7)})\varepsilon^{\frac{1}{32s}}}
<\frac{\varepsilon^{\frac{31}{16}}}{200}.
\end{multline}
Moreover, by (\ref{uniform_eialpha}) in Lemma \ref{lem_a1}, we have
\begin{multline}\label{est_f4}
	\left\|\exp\left\{{\frac{\rm i}2\alpha(\xi\eta+e^{-\frac{\rm i}2\alpha}\eta \bar q+e^{\frac{\rm i}2\alpha}\xi \bar p+\bar p\bar q)}\right\}\bar q- e^{\frac{\rm i}2\alpha(\xi\eta)}\bar q\right\|_{\CO,\tilde\beta,r^{(7)}}\\
	\leq\left\|\exp\left\{{\frac{\rm i}2\alpha(\xi\eta+e^{-\frac{\rm i}2\alpha}\eta \bar q+e^{\frac{\rm i}2\alpha}\xi \bar p+\bar p\bar q)}\right\} - e^{\frac{\rm i}2\alpha(\xi\eta)}\right\|_{\CO,\tilde\beta,r^{(7)}}\|q\|_{\CO,\tilde\beta,r^{(7)}}
	<\frac43\cdot\frac{\varepsilon}{10}\cdot \frac{\varepsilon}{10}
	=\frac{\varepsilon^2}{75}.
\end{multline}
Hence,
(\ref{appro_f}) is shown by combining $(\ref{est_f1})-(\ref{est_f4})$.
The proof for (\ref{appro_g}) is similar.\qed

\begin{cor}\label{cor_fg}
$\|f\|_{\CO,\tilde\beta,r^{(7)}}, \ \|g\|_{\CO,\tilde\beta,r^{(7)}}< \frac\varepsilon4$.
\end{cor}
\proof
Lemma \ref{lem_eia} implies that
$$\|p(e^{-\frac{\rm i}2\alpha}\eta,e^{\frac{\rm i}2\alpha}\xi)\|_{\CO,\tilde\beta,r^{(7)}}, \   \|q(e^{-\frac{\rm i}2\alpha}\eta,e^{\frac{\rm i}2\alpha}\xi)\|_{\CO,\tilde\beta,r^{(7)}}< \frac{\varepsilon}{10}.$$
Moreover, we have
\begin{gather*}
	\left\|\frac{{\rm i}\alpha'(\xi\eta)}2 (e^{-\frac{\rm i}2\alpha}\eta \bar q+e^{\frac{\rm i}2\alpha}\xi \bar p)e^{{\rm i}\alpha}\xi\right\|_{\CO,\tilde\beta,r^{(7)}},
	\left\|\frac{{\rm i}\alpha'(\xi\eta)}2 (e^{-\frac{\rm i}2\alpha}\eta \bar q+e^{\frac{\rm i}2\alpha}\xi \bar p)e^{-{\rm i}\alpha}\eta\right\|_{\CO,\tilde\beta,r^{(7)}} < \frac{3e^{2\tilde\beta} r^2\varepsilon}{50},  \\
	\|e^{\frac{\rm i}2\alpha}  \bar q\|_{\CO,\tilde\beta,r^{(7)}} , \  \|e^{-\frac{\rm i}2\alpha}  \bar p\|_{\omega,\tilde\beta,r^{(7)}} <\frac{e^{\frac{9}{16}\tilde\beta}\varepsilon}{10}.
\end{gather*}
By Lemma \ref{lem_appro}, the corollary is shown.\qed

\smallskip

Let $C(\xi,\eta):=\frac{\rm i}{2}\alpha'(\xi\eta)\left(e^{-\frac{\rm i}{2}\alpha(\xi\eta)}\eta \bar q (\xi,\eta)+e^{\frac{\rm i}{2}\alpha(\xi\eta)}\xi \bar p (\xi,\eta)\right)$. Applying (\ref{esi_coeff}), we have
\begin{cor}\label{cor_fg_coeff}
$C \in{\CA}_{\tilde\beta,r^{(7)}}({\cL O})$ and we have
\begin{align}
\left\|f_{l,0}-e^{\frac{\rm i}{2}\alpha} \bar q_{l,0}-e^{\frac{\rm i}{2}l\alpha} p_{0,l}-e^{{\rm i}\alpha}C_{l-1,0}\right\|_{\CO,\tilde\beta,r^{(7)}} &< \frac{\varepsilon^{\frac{31}{16}}}{80(r^{(7)})^{l}},\quad l\geq 2, \label{appro_f_coeff_l}\\
\left\|f_{0,j}-e^{\frac{\rm i}{2}\alpha} \bar q_{0,j}-e^{-\frac{\rm i}{2}j\alpha} p_{j,0}-(\xi\eta)e^{{\rm i}\alpha}C_{0,j+1}\right\|_{\CO,\tilde\beta,r^{(7)}} &< \frac{\varepsilon^{\frac{31}{16}}}{80(r^{(7)})^{j}},\quad j\geq 0, \label{appro_f_coeff_j}\\
\left\|g_{l,0}-e^{-\frac{\rm i}{2}\alpha} \bar p_{l,0}-e^{\frac{\rm i}{2}l\alpha} q_{0,l}+(\xi\eta)e^{-{\rm i}\alpha} C_{l+1,0}\right\|_{\CO,\tilde\beta,r^{(7)}} &< \frac{\varepsilon^{\frac{31}{16}}}{80(r^{(7)})^{l}},\quad l\geq 0, \label{appro_g_coeff_l}\\
\left\|g_{0,j}-e^{-\frac{\rm i}{2}\alpha} \bar p_{0,j} - e^{-\frac{\rm i}{2}j\alpha} q_{j,0}+e^{-{\rm i}\alpha} C_{0,j-1} \right\|_{\CO,\tilde\beta,r^{(7)}} &< \frac{\varepsilon^{\frac{31}{16}}}{80(r^{(7)})^{j}},\quad j\geq 2. \label{appro_g_coeff_j}
\end{align}
\end{cor}

\smallskip

Recall that $\tau_1\circ\tau_1={\rm Id}$, which means that
$$\left(\begin{array}{l}
                          \xi \\
                           \eta
                          \end{array}
\right)=
\left(\begin{array}{l}
\exp\left\{{\frac{\rm i}2\left(\alpha(\xi\eta+e^{\frac{\rm i}2\alpha}\eta q+e^{-\frac{\rm i}2\alpha}\xi p+pq)-\alpha(\xi\eta)\right)}\right\}\xi\\[2mm]
+ \, \exp\left\{{\frac{\rm i}2\alpha(\xi\eta+e^{\frac{\rm i}2\alpha}\eta q+e^{-\frac{\rm i}2\alpha}\xi p+pq)}\right\}q\\[2mm]
+ \, p(e^{\frac{\rm i}2\alpha}\eta+p,e^{-\frac{\rm i}2\alpha}\xi+q), \\[3mm]
\exp\left\{-{\frac{\rm i}2\left(\alpha(\xi\eta+e^{\frac{\rm i}2\alpha}\eta q+e^{-\frac{\rm i}2\alpha}\xi p+pq)-\alpha(\xi\eta)\right)}\right\}\eta\\[2mm]
+ \,  \exp\left\{-{\frac{\rm i}2\alpha(\xi\eta+e^{\frac{\rm i}2\alpha}\eta q+e^{-\frac{\rm i}2\alpha}\xi p+pq)}\right\}p\\[2mm]
+ \,  q(e^{\frac{\rm i}2\alpha}\eta+p,e^{-\frac{\rm i}2\alpha}\xi+q)
\end{array}
\right).$$
Then, similarly to Lemma \ref{lem_appro}, we have
\begin{lemma}We have
\begin{align}
   \left\|\frac{{\rm i}\alpha'(\xi\eta)}{2} (e^{\frac{\rm i}2\alpha}\eta q+e^{-\frac{\rm i}2\alpha}\xi p)\xi+e^{\frac{\rm i}2\alpha}q+p(e^{\frac{\rm i}2\alpha}\eta,e^{-\frac{\rm i}2\alpha}\xi)\right\|_{\CO,\tilde\beta,r^{(7)}}&< \frac{\varepsilon^{\frac{31}{16}}}{80},\label{appro_p_}\\
\left\|-\frac{{\rm i}\alpha'(\xi\eta)}{2} (e^{\frac{\rm i}2\alpha}\eta q+e^{-\frac{\rm i}2\alpha}\xi p) \eta+e^{-\frac{\rm i}2\alpha}p+q(e^{\frac{\rm i}2\alpha}\eta,e^{-\frac{\rm i}2\alpha}\xi)\right\|_{\CO,\tilde\beta,r^{(7)}} &<\frac{\varepsilon^{\frac{31}{16}}}{80}.\label{appro_q_}
\end{align}
\end{lemma}

\begin{cor}\label{cor_pq_lj} We have, for $l,j\geq 1$
	\begin{align}
\label{pq_0110_}\|(\xi\eta)(e^{\frac{\rm i}2\alpha} q_{1,0}+e^{-\frac{\rm i}2\alpha}p_{0,1})\|_{\CO,\tilde\beta,r^{(7)}}&<\frac{r^{(7)}\varepsilon^{\frac{31}{16}}}{80},\\
\label{pq_l+-1}\left\|(\xi\eta)(e^{\frac{\rm i}2\alpha} q_{l+1,0}+e^{-\frac{\rm i}2(l+1)\alpha} p_{0,l+1})+e^{-\frac{\rm i}2\alpha}p_{l-1,0}+e^{-\frac{\rm i}2(l-1)\alpha} q_{0,l-1}\right\|_{\CO,\tilde\beta,r^{(7)}}&<\frac{\varepsilon^{\frac{31}{16}}}{40(r^{(7)})^{l-1}},\\
\label{pq_j+-1}\left\|(\xi\eta)(e^{-\frac{\rm i}2\alpha}p_{0,j+1}+e^{\frac{\rm i}2(j+1)\alpha}
		q_{j+1,0})+e^{\frac{\rm i}2\alpha} q_{0,j-1}+e^{\frac{\rm i}2(j-1)\alpha} p_{j-1,0}\right\|_{\CO,\tilde\beta,r^{(7)}}&<\frac{\varepsilon^{\frac{31}{16}}}{40(r^{(7)})^{j-1}}.
	\end{align}
\end{cor}
\proof Due to cancellation of terms, we have that
\begin{eqnarray*}
& &  \eta \left(\frac{{\rm i}\alpha'(\xi\eta)}{2} (e^{\frac{\rm i}2\alpha}\eta q+e^{-\frac{\rm i}2\alpha}\xi p)\xi+e^{\frac{\rm i}2\alpha}q+p(e^{\frac{\rm i}2\alpha}\eta,e^{-\frac{\rm i}2\alpha}\xi)\right)\\
& &+ \,  \xi \left(-\frac{{\rm i}\alpha'(\xi\eta)}{2} (e^{\frac{\rm i}2\alpha}\eta q+e^{-\frac{\rm i}2\alpha}\xi p) \eta+e^{-\frac{\rm i}2\alpha}p+q(e^{\frac{\rm i}2\alpha}\eta,e^{-\frac{\rm i}2\alpha}\xi)\right)\\
&=&e^{\frac{\rm i}2\alpha}\eta q+e^{-\frac{\rm i}2\alpha}\xi p+ \eta p(e^{\frac{\rm i}2\alpha}\eta,e^{-\frac{\rm i}2\alpha}\xi)+\xi q(e^{\frac{\rm i}2\alpha}\eta,e^{-\frac{\rm i}2\alpha}\xi).
\end{eqnarray*}
Hence, by (\ref{appro_p_})  and (\ref{appro_q_}),
$$\left\|e^{\frac{\rm i}2\alpha}\eta q+e^{-\frac{\rm i}2\alpha}\xi p+ \eta p(e^{\frac{\rm i}2\alpha}\eta,e^{-\frac{\rm i}2\alpha}\xi)+\xi q(e^{\frac{\rm i}2\alpha}\eta,e^{-\frac{\rm i}2\alpha}\xi)\right\|_{\CO,\tilde\beta,r^{(7)}}
\leq 2 r^{(7)}\cdot \frac{\varepsilon^{\frac{31}{16}}}{80}=\frac{r^{(7)}\varepsilon^{\frac{31}{16}}}{40}.$$
The corresponding coefficients under the decomposition (\ref{decomp}) satisfy \re{pq_0110_} -- \re{pq_j+-1}.\qed

\begin{cor}\label{cor_coeff_res_pq_fg} We have
\begin{eqnarray}
\|e^{\frac{\rm i}2\alpha} q_{1,0}+e^{-\frac{\rm i}2\alpha}p_{0,1}\|_{\CO,\tilde\beta,r^{(7)}}&<&\frac{\varepsilon^{\frac{61}{32}}}{60r^{(7)}},\label{coeff_res_pq}\\
\|e^{\frac{\rm i}{2}\alpha}\bar q_{1,0}+e^{\frac{\rm i}{2}\alpha}p_{0,1}- f_{1,0}\|_{\CO,\tilde\beta,r^{(7)}}&<&\frac{\varepsilon^{\frac{61}{32}}}{60r^{(7)}}, \label{coeff_res_f}\\
\|e^{-\frac{\rm i}{2}\alpha}\bar p_{0,1}+e^{-\frac{\rm i}{2}\alpha}q_{1,0}- g_{1,0}\|_{\CO,\tilde\beta,r^{(7)}}&<&\frac{\varepsilon^{\frac{61}{32}}}{60r^{(7)}}.\label{coeff_res_g}
\end{eqnarray}
\end{cor}

\proof Note that (\ref{pq_0110_}) actually means that
$$
\|(e^{\frac{\rm i}2\alpha}\eta q+e^{-\frac{\rm i}2\alpha}\xi p)_{0,0}\|_{\CO,\tilde\beta,r^{(7)}}<\frac{r^{(7)}\varepsilon^{\frac{31}{16}}}{80}.$$
Hence we obtain (\ref{coeff_res_pq}), since in (\ref{appro_p_}), the coefficients of the term $\xi^1\eta^0$ satisfies
$$
\left\|\frac{{\rm i}\alpha'(\xi\eta)}{2}  (e^{\frac{\rm i}2\alpha}\eta q+e^{-\frac{\rm i}2\alpha}\xi p)_{0,0}+e^{\frac{\rm i}2\alpha} q_{1,0}+e^{-\frac{\rm i}2\alpha}p_{0,1}\right\|_{\CO,\tilde\beta,r^{(7)}}<\frac{\varepsilon^{\frac{31}{16}}}{80r^{(7)}}.
$$
In (\ref{appro_f}), the coefficients of the term $\xi$ satisfies
$$\left\|f_{1,0}-\frac{{\rm i}\alpha'(\xi\eta)}2 (e^{-\frac{\rm i}2\alpha}\eta \bar q+e^{\frac{\rm i}2\alpha}\xi \bar p)_{0,0}e^{{\rm i}\alpha}-e^{\frac{\rm i}2\alpha}\bar q_{1,0}-e^{\frac{\rm i}2\alpha}p_{0,1}\right\|_{\CO,\tilde\beta,r^{(7)}}<\frac{\varepsilon^{\frac{31}{16}}}{80r^{(7)}},$$
then we have (\ref{coeff_res_f}). The proof for (\ref{coeff_res_g}) is similar by applying (\ref{appro_g}).
\qed

\subsection{The skew terms}

In the KAM-like (or Newton) scheme stated in Section \ref{sec_KAM}, we also consider the skew term $e^{\frac{\rm i}2\alpha}\eta q+e^{-\frac{\rm i}2\alpha}\xi p$ of $\tau_1$, and the skew term $e^{-{\rm i}\alpha}\eta f+e^{{\rm i}\alpha}\xi g$ of $\sigma$.
By Lemma \ref{lemma_analytic_prod} and \ref{lem_appro}, we see that
both skew terms belong to ${\CA}_{\tilde\beta,r^{(7)}}({\cL O})$.

\begin{lemma}\label{lem_cond_F} We have
 \begin{equation}\label{esti_fg}
 \|e^{-{\rm i}\alpha}\eta f+e^{{\rm i}\alpha}\xi g\|_{\CO,\tilde\beta,r^{(7)}}<2\|e^{\frac{\rm i}2\alpha}\eta q+e^{-\frac{\rm i}2\alpha}\xi p\|_{\CO,\beta,r}+\frac{\varepsilon^{\frac{31}{16}}}{40}.
 \end{equation}
\end{lemma}
\proof According to (\ref{appro_f}) and (\ref{appro_g}), $e^{-{\rm i}\alpha}\eta f+e^{{\rm i}\alpha}\xi g$ can be approximated by
\begin{eqnarray*}
 & &  e^{-{\rm i}\alpha}\eta\left( \frac{{\rm i}\alpha'}2 \left(e^{-\frac{\rm i}2\alpha}\eta \bar q+e^{\frac{\rm i}2\alpha}\xi \bar p\right)e^{{\rm i}\alpha}\xi+e^{\frac{\rm i}2\alpha}\bar q+p(e^{-\frac{\rm i}2\alpha}\eta,e^{\frac{\rm i}2\alpha}\xi)\right)\\
 & &    +  \,   e^{{\rm i}\alpha}\xi\left(- \frac{{\rm i}\alpha'}2\left(e^{-\frac{\rm i}2\alpha}\eta \bar q+e^{\frac{\rm i}2\alpha}\xi \bar p\right)e^{-{\rm i}\alpha}\eta+ e^{-\frac{\rm i}2\alpha}\bar p+q(e^{-\frac{\rm i}2\alpha}\eta,e^{\frac{\rm i}2\alpha}\xi)\right) \\
 &=&  e^{-\frac{\rm i}2\alpha}\eta\bar q+e^{\frac{\rm i}2\alpha}\xi\bar p+ e^{-\frac{\rm i}2\alpha}\left(e^{-\frac{\rm i}2\alpha}\eta\right) p(e^{-\frac{\rm i}2\alpha}\eta,e^{\frac{\rm i}2\alpha}\xi)+e^{\frac{\rm i}2\alpha}\left(e^{\frac{\rm i}2\alpha}\xi\right) q(e^{-\frac{\rm i}2\alpha}\eta,e^{\frac{\rm i}2\alpha}\xi),
\end{eqnarray*}
up to an error smaller than $\displaystyle \frac{\varepsilon^{\frac{31}{16}}}{40}$. In view of (\ref{esi_coeff}) and (\ref{esti_mono}), we have
 $$\|e^{-\frac{\rm i}2\alpha}\eta\bar q+e^{\frac{\rm i}2\alpha}\xi\bar p\|_{\CO,\tilde\beta,r^{(7)}}=\|e^{\frac{\rm i}2\alpha}\eta q+e^{-\frac{\rm i}2\alpha}\xi p\|_{\CO,\tilde\beta,r^{(7)}}\leq \|e^{\frac{\rm i}2\alpha}\eta q+e^{-\frac{\rm i}2\alpha}\xi p\|_{\CO,\beta,r},$$
By (\ref{esi_coeff}) and Lemma \ref{lem_eia}, we have
$$ \| e^{-\frac{\rm i}2\alpha}(e^{-\frac{\rm i}2\alpha}\eta) p(e^{-\frac{\rm i}2\alpha}\eta,e^{\frac{\rm i}2\alpha}\xi)+e^{\frac{\rm i}2\alpha}(e^{\frac{\rm i}2\alpha}\xi) q(e^{-\frac{\rm i}2\alpha}\eta,e^{\frac{\rm i}2\alpha}\xi) \|_{\CO,\tilde\beta,r^{(7)}}
 \leq \, \|e^{\frac{\rm i}2\alpha}\eta q+e^{-\frac{\rm i}2\alpha}\xi p\|_{\CO,\beta,r}.$$
Then (\ref{esti_fg}) is shown by combining the error smaller than $\displaystyle \frac{\varepsilon^{\frac{31}{16}}}{40}$.\qed

\section{Transformations on crowns}\label{sec_trans}

Fix $s\in\N^*$, $0<r_+<r<\frac{1}{4}$, $0<\varepsilon,\beta<r^2$ as in Section \ref{sec_rev}.
In this section, we introduce two types of transformations on the ``crown"
${\cL C}^r_{\omega,\beta}$, which will be used in the KAM-like scheme.

\subsection{Product-preserving scaling transformation}\label{sub_product_preserving}

Consider the map
\begin{equation}\label{tau_PD}
\tau(\xi,\eta)=\left(\begin{array}{l}
(e^{\frac{\rm i}2\theta(\xi\eta)}+A(\xi\eta))\eta+ p(\xi,\eta) \\
(e^{\frac{\rm i}2\theta(\xi\eta)}+A(\xi\eta))^{-1}\xi+ q(\xi,\eta)
\end{array}
\right),\quad (\xi,\eta)\in {\cL C}^{r'}_{\omega,\beta'},
\end{equation}
with $0\leq\beta'<\tilde\beta$ and $0<r'<(r^2-\beta)^{\frac12}$,
where $ A=A(\xi\eta)\in{\CA}_{\beta,r}(\CO)$ with $\| A\|_{\CO,\beta,r}<\frac1{16}$, $p,q\in{\CA}_{\beta,r}(\CO)$,
$\theta=\theta(\xi\eta)\in{\CA}_{\beta,r}^\R(\CO)$, with $\theta(\omega)\in \, ]-\frac12,4\pi+\frac12[$ for $\omega\in{\CO}(r,\beta)$ and  $\|\theta\|_{\CO,\beta,r}<4\pi+1$. Furthermore, we assume $|\theta^{(s)}-s!|_{{\CO}(r,\beta)}<\frac{s!}{15}$ and
$$
|\theta^{(k)}|_{{\CO}(r,\beta)}<\begin{cases}\frac{1}{15},\quad 1\leq k\leq s-1 \ {\rm and} \ s\geq 2,\\[2mm]
\frac{r^{-\frac{1}{2}}}2,\quad s+1\leq k\leq 16s.\end{cases}
$$
$\theta$ satisfies all the hypothesis of $\alpha$ in Section \ref{sec_rev} (see \re{al_nu_sum} -- \re{Ds+nu_sum} in Section \ref{sec_KAM}). Hence all the lemmas and corollaries in Subsection \ref{sub_alpha} are applicable on $\theta$. In particular, by Remark \ref{rem_eiba}, we have
  \begin{equation}\label{eibtheta1}
\|e^{{\rm i}b\theta}\|_{\CO,\beta',r'} < \frac{101}{100},\quad \forall \  -\beta^{-\frac14}\leq b\leq \beta^{-\frac14}.
\end{equation}

For $(\xi,\eta)\in {\cL C}^{r'}_{\omega,\beta'}$, $\omega\in{\CO}(r',\beta') \, $, let
\begin{equation}\label{Lambda}
\Theta(\xi\eta):=\left((e^{\frac{\rm i}2\theta(\xi\eta)}+A(\xi\eta))(e^{-\frac{\rm i}2\theta(\xi\eta)}+\bar A(\xi\eta)) \right)^{\frac14}
\end{equation}
be the fourth root,
and let us set
\begin{equation}\label{varphi}
\varphi:(\xi,\eta)\mapsto \left(\Theta(\xi\eta)\xi, \Theta^{-1}(\xi\eta)\eta\right).
\end{equation}
It is easy to show that $\rho\circ\varphi=\varphi\circ\rho$.

 \begin{lemma}\label{lemma_Lambda_k}
For $k=\pm 1$, $\pm 2$, $\Theta^k\in{\CA}^\R_{\beta',r'}(\CO)$ satisfies that
\begin{equation}\label{Lambda_k}
\|\Theta^k-1\|_{\CO,\beta',r'}<\frac{3|k|}4 \| A\|_{\CO,\beta',r'}.
\end{equation}
\end{lemma}
\proof
Since $ A\in{\CA}_{\beta,r}(\CO)$ with $\| A\|_{\CO,\beta,r}<\frac1{16}$, by \re{eibtheta1}, we have that $e^{\frac{\rm i}2\theta}\bar A+e^{-\frac{\rm i}2\theta} A+ A\bar A\in{\CA}^\R_{\beta',r'}(\CO)$ with
 $$\|e^{\frac{\rm i}2\theta}\bar A+e^{-\frac{\rm i}2\theta} A+ A\bar A\|_{\CO,\beta',r'}<\frac{101}{50}\| A\|_{\CO,\beta,r}+\| A\|_{\CO,\beta,r}^2< \frac{21}{10} \| A\|_{\CO,\beta,r}.$$
For $k=\pm 1, \ \pm 2$,
we have
\begin{eqnarray*}
 \|(1+e^{\frac{\rm i}2\theta}\bar A+e^{-\frac{\rm i}2\theta} A+ A\bar A)^{\frac{k}4}-1\|_{\CO,\beta',r'}&\leq& \frac{|k|}4  \frac{\|e^{\frac{\rm i}2\theta}\bar A+e^{-\frac{\rm i}2\theta} A+ A\bar A\|_{\CO,\beta',r'}}{\left(1-\|e^{\frac{\rm i}2\theta}\bar A+e^{-\frac{\rm i}2\theta} A+ A\bar A\|_{\CO,\beta',r'}\right)^{1-\frac{k}{4}}} \\
 &\leq& \frac{|k|}4 \cdot \frac{21}{10}  \| A\|_{\CO,\beta,r} \cdot \left(1-\frac{21}{10} \| A\|_{\CO,\beta,r}\right)^{-\frac{3}{2}}  \\
&\leq& \frac{3|k|}4 \| A\|_{\CO,\beta,r}.
\end{eqnarray*}
By the expression of $\Theta$ in (\ref{Lambda}), we obtain that $\Theta^k\in{\CA}^\R_{\beta',r'}(\CO)$ and (\ref{Lambda_k}) is satisfied.\qed

\smallskip

\begin{prop}\label{prop_transform_product}
If $\| A\|_{\CO,\beta,r}$ satisfies that
\begin{equation}\label{r_Lambda}
\| A\|_{{\CO},\beta,r}<\frac{4(r-r')}{3r'},
\end{equation}
then $\varphi$, as defined by \re{varphi}, satisfies $\|\varphi-{\rm Id}\|_{\CO,\beta', r'}<\frac34\| A\|_{\CO,\beta,r}$ and, for every $\omega\in{\CO}(r',\beta')$, we have $\varphi({\cL C}^{r'}_{\omega,\beta'})\subset {\cL C}^{r}_{\omega,\beta}$. Furthermore, there are $\theta_+\in{\CA}^\R_{\beta', r'}({\CO})$,  $p_+, \, q_+\in{\CA}_{\beta', r'}({\CO})$ such that
$$(\varphi^{-1}\circ\tau\circ\varphi)(\xi,\eta)=\left(\begin{array}{l}
e^{\frac{\rm i}2\theta_+(\xi\eta)}\eta+p_+(\xi,\eta) \\
e^{-\frac{\rm i}2\theta_+(\xi\eta)}\xi+q_+(\xi,\eta)
\end{array}
\right),$$
\begin{eqnarray}
\theta_+(\xi\eta)-\theta(\xi\eta)&=& -{\rm i}(e^{-\frac{\rm i}2\theta(\xi\eta)}A(\xi\eta)-e^{\frac{\rm i}2\theta(\xi\eta)}\bar A(\xi\eta)),\label{theta_+}\\
\|p_+\|_{{\CO},\beta', r'}&<&\left(1+\frac34\| A\|_{\CO,\beta,r}\right)\|p\|_{\CO,\beta,r}+\| A\|_{\CO,\beta,r}^2,\label{p_+}\\
\|q_+\|_{{\CO},\beta', r'}&<&\left(1+\frac34\| A\|_{\CO,\beta,r}\right)\|q\|_{\CO,\beta,r}+\| A\|_{\CO,\beta,r}^2.\label{q_+}
\end{eqnarray}
Moreover, we have
\begin{eqnarray}
& &\|e^{-\frac{\rm i}2\theta_+}\xi p_++e^{-\frac{\rm i}2\theta_+}\eta q_+\|_{{\CO},\beta', r'}\nonumber\\
&<&\|e^{-\frac{\rm i}2\theta} \xi p+ e^{\frac{\rm i}2\theta} \eta  q\|_{\CO,\beta,r}+\| A\|_{\CO,\beta,r}(\|p\|_{\CO,\beta,r}+\|q\|_{\CO,\beta,r})+\| A\|_{\CO,\beta,r}^2.\label{crossing_theta_+}
\end{eqnarray}
\end{prop}

\begin{remark}
At each KAM step, we always work with the involutions of the form \re{gen-invol} with  $p,q\in{\CA}_{\beta,r}(\CO)$ satisfying  \re{esti_pq}.
After conjugacy by the KAM transformation, the new involution, denoted by $\tilde\tau_1$ (resp.  $\tilde\tau_2=\rho\circ\tilde\tau_1\circ\rho$), has the form
$$\tilde\tau_1(\xi,\eta)=\left(\begin{array}{l}
 \lambda(\xi\eta)\eta+\tilde p(\xi,\eta) \\
 \lambda^{-1}(\xi\eta)\xi+\tilde q(\xi,\eta)
\end{array}\right),$$
with new perturbations $\tilde p$ and $\tilde q$ of much smaller size than $p$, $q$. Nevertheless, the new principal part usually does not satisfy $|\lambda(\omega)|= 1$ (but close to $1$). Hence, an additional transformation $\varphi$ (as in (\ref{varphi}), see also \re{scaling_MW}) is needed in order to preserve the same form as $\tau_1$ in \re{gen-invol}.
This is similar to the role of Theorem 3.4 in \cite{moser-webster} for the formal hyperbolic non-exceptional manifold case.
\end{remark}

Before the proof of Proposition \ref{prop_transform_product}, we show the following lemma similar to Lemma \ref{lem_eia}.

\begin{lemma}\label{lemma_h_Lambda}
For $0<r'<r$ and $h\in {\CA}_{\beta,r}(\CO)$ with $\|h\|_{\CO,\beta,r}< +\infty$, if $\| A\|_{\CO,\beta,r}$ is sufficiently small such that (\ref{r_Lambda}) is satisfied,
then we have $h(\Theta\xi,\Theta^{-1}\eta)\in {\CA}_{\beta',r'}(\CO)$ with
\begin{equation}\label{h_Lambda}
\|h(\Theta\xi,\Theta^{-1}\eta)\|_{\CO,\beta',r'}<  \|h\|_{\CO,\beta,r}.
\end{equation}
\end{lemma}
\proof According to (\ref{r_Lambda}), we have $\left(1+\frac{3}{4}| A|_{\omega,\beta'}\right)\frac{r'}{r}=1+\frac{\frac{3}{4}| A|_{\omega,\beta'}r'-(r-r')}{r}<1$. Therefore, with the help of (\ref{Lambda_k}), we obtain for $\omega \in {\CO}(r',\beta')$~:
\begin{eqnarray*}
\left\|h(\Theta\xi,\Theta^{-1}\eta)\right\|_{\omega,\beta',r'}
&\leq&  |h_{0,0}|_{\omega,\beta'}+\sum_{l\geq 1} |h_{l,0}|_{\omega,\beta'} |\Theta|_{\omega,\beta'}^l r'^l+\sum_{j\geq 1} |h_{0,j}|_{\omega,\beta'}|\Theta^{-1}|_{\omega,\beta'}^j r'^j\nonumber\\
&<&  |h_{0,0}|_{\omega,\beta'}+\sum_{l\geq 1} \left(1+\frac{3}{4}| A|_{\omega,\beta'}\right)^l\left(\frac{r'}{r}\right)^l(|h_{l,0}|_{\omega,\beta'}+|h_{0,l}|_{\omega,\beta'})r^l \\
&<&   \|h\|_{\omega,\beta,r}.\end{eqnarray*}
Since for $l,j\geq 0$ with $lj=0$ and  $\omega \in {\CO}(r',\beta')$, $\left(h(\Theta\xi,\Theta^{-1}\eta)\right)_{l,j}(\omega)=h_{l,j}(\omega) \cdot \Theta^{l-j}(\omega)$, we see that $h(\Theta\xi,\Theta^{-1}\eta) \in{\CA}_{\beta',  r'}({\CO})$.\qed

\smallskip

\noindent {\it Proof of Proposition \ref{prop_transform_product}.}
With $\varphi$ defined in (\ref{varphi}), we have
$$\varphi^{-1}\circ\tau\circ\varphi=\left(\begin{array}{c}
     \Theta^{-2}(\xi\eta)(e^{\frac{\rm i}2\theta(\xi\eta)}+A(\xi\eta))\eta\\[1mm]
     \Theta^2(\xi\eta)(e^{\frac{\rm i}2\theta(\xi\eta)}+A(\xi\eta))^{-1} \xi
   \end{array}\right)+\left(\begin{array}{c}
    \Theta^{-1}(\xi\eta)p(\Theta\xi,\Theta^{-1}\eta)\\[1mm]
     \Theta(\xi\eta)q(\Theta\xi,\Theta^{-1}\eta)
   \end{array}\right).$$
In view of (\ref{h_Lambda}), we obtain
$$\left\|p(\Theta\xi,\Theta^{-1}\eta)\right\|_{\CO, \beta',r'}<\|p\|_{\CO,\beta,r}, \quad \left\|q(\Theta\xi,\Theta^{-1}\eta)\right\|_{\CO,\beta',r'}<\|q\|_{\CO,\beta,r}.$$
Combining with (\ref{Lambda_k}) together with $B:=1+\frac34\| A\|_{\CO,\beta,r}$, we have
$$\left\|\Theta^{-1} p(\Theta\xi,\Theta^{-1}\eta)\right\|_{\CO,\beta',r'}<B\|p\|_{\CO,\beta,r},\quad\left\|\Theta q(\Theta\xi,\Theta^{-1}\eta)\right\|_{\CO,\beta',r'}<B\|q\|_{\CO,\beta,r}.
$$

By a direct computation, we have, for $(\xi,\eta)\in{\cL C}^{r'}_{\omega,\beta'}$, $\omega \in {\CO}(r',\beta')$,
\begin{eqnarray*}
\Theta^{-2}(\xi\eta)(e^{\frac{\rm i}2\theta(\xi\eta)}+A(\xi\eta))
&=&\frac{e^{\frac{\rm i}2\theta(\xi\eta)}+A(\xi\eta)}{(e^{\frac{\rm i}2\theta(\xi\eta)}+A(\xi\eta))^{\frac12}(e^{-\frac{\rm i}2\theta(\xi\eta)}+\bar A(\xi\eta))^{\frac12}}\\
&=&(e^{\frac{\rm i}2\theta(\xi\eta)}+A(\xi\eta))^{\frac12}(e^{-\frac{\rm i}2\theta(\xi\eta)}+\bar A(\xi\eta))^{-\frac12}\\
&=&e^{\frac{\rm i}2\theta(\xi\eta)}(1+e^{-\frac{\rm i}2\theta(\xi\eta)}A(\xi\eta))^{\frac12}(1+e^{\frac{\rm i}2\theta(\xi\eta)}\bar A(\xi\eta))^{-\frac12},\\
\Theta^{2}(\xi\eta)(e^{\frac{\rm i}2\theta(\xi\eta)}+A(\xi\eta))^{-1}&=&e^{-\frac{\rm i}2\theta(\xi\eta)}(1+e^{-\frac{\rm i}2\theta(\xi\eta)}A(\xi\eta))^{-\frac12}(1+e^{\frac{\rm i}2\theta(\xi\eta)}\bar A(\xi\eta))^{\frac12}.
\end{eqnarray*}
According to (\ref{eibtheta1}) and (\ref{r_Lambda}), we see that $e^{-\frac{\rm i}2\theta(\xi\eta)}A(\xi\eta)$ and $e^{\frac{\rm i}2\theta(\xi\eta)}\bar A(\xi\eta)$ are small enough. Then, by the expansion of power series, we have
\begin{eqnarray*}
& &\ln\left((1+e^{-\frac{\rm i}2\theta(\xi\eta)}A(\xi\eta))^{\frac12}(1+e^{\frac{\rm i}2\theta(\xi\eta)}\bar A(\xi\eta))^{-\frac12}\right) \\
&=& \ln\left(1+\frac12\left(e^{-\frac{\rm i}2\theta(\xi\eta)}A(\xi\eta)-e^{\frac{\rm i}2\theta(\xi\eta)}\bar A(\xi\eta)\right)+\tilde{\CP}(\xi\eta))\right)\\
&=&\frac12\left(e^{-\frac{\rm i}2\theta(\xi\eta)}A(\xi\eta)-e^{\frac{\rm i}2\theta(\xi\eta)}\bar A(\xi\eta)\right)+{\CP}(\xi\eta),
\end{eqnarray*}
where $\tilde{\CP}$, ${\CP}$ are sums of powers of $(e^{-\frac{\rm i}2\theta}A, e^{\frac{\rm i}2\theta}\bar A)$ of order greater than $2$, satisfying that
$$\|\tilde{\CP}\|_{\CO,\beta',r'}, \  \|{\CP}\|_{\CO,\beta',r'}<\frac12\|A\|_{\CO,\beta',r'}^2.$$
Now we define
$\theta_+(\xi\eta):=\theta(\xi\eta)-{\rm i}(e^{-\frac{\rm i}2\theta(\xi\eta)}A(\xi\eta)-e^{\frac{\rm i}2\theta(\xi\eta)}\bar A(\xi\eta))$
as in \re{theta_+},
which belongs to ${\CA}^\R_{\beta',r'}(\CO)$.
Then we rewrite $(1+e^{-\frac{\rm i}2\theta}A)^{\frac12}(1+e^{\frac{\rm i}2\theta}\bar A)^{-\frac12}$ as
$$(1+e^{-\frac{\rm i}2\theta}A)^{\frac12}(1+e^{\frac{\rm i}2\theta}\bar A)^{-\frac12} =e^{\frac{\rm i}2(\theta_+-\theta)}e^{\CP}.$$
Since (\ref{eibtheta1}) implies that
$$\|e^{-\frac{\rm i}2\theta} A-e^{\frac{\rm i}2\theta}\bar A\|_{\CO,\beta',r'}\leq 2\|e^{-\frac{\rm i}2\theta} A\|_{\CO,\beta',r'} <\frac{101}{50}\| A\|_{\CO,\beta,r},$$
we have
\begin{eqnarray}
\left\|\Theta^{-2}\cdot(e^{\frac{\rm i}2\theta}+ A)-e^{\frac{\rm i}2\theta_+}\right\|_{\CO,\beta',r'}&=&\left\|e^{\frac{\rm i}2\theta}(1+e^{-\frac{\rm i}2\theta}A)^{\frac12}(1+e^{\frac{\rm i}2\theta}\bar A)^{-\frac12}-e^{\frac{\rm i}2\theta_+}\right\|_{\CO,\beta',r'}\label{err_fre_1}\\
&=&\left\|e^{\frac{\rm i}2\theta} e^{\frac{\rm i}2(\theta_+-\theta)}e^{\CP}   -e^{\frac{\rm i}2\theta_+}\right\|_{\CO,\beta',r'}\nonumber\\
&\leq &  \|e^{\frac{\rm i}2\theta}\|_{\CO,\beta',r'}\|e^{\frac{\rm i}2(\theta_+-\theta)}\|_{\CO,\beta',r'} \|e^{\CP}-1\|_{\CO,\beta',r'} \nonumber\\
 &<& \frac{101}{100}\left(1+\frac{101}{40}\| A\|_{\CO,\beta,r}\right)\cdot \frac23\| A\|^2_{\CO,\beta',r'}  \nonumber\\
  &<& \| A\|^2_{\CO,\beta',r'}.  \nonumber
\end{eqnarray}
Similarly, we have
\begin{equation}\label{err_fre_2}
\left\|\Theta^{2}\cdot(e^{\frac{\rm i}2\theta}+ A)^{-1}-e^{-\frac{\rm i}2\theta_+}\right\|_{\CO,\beta',r'}<\| A\|^2_{\CO,\beta,r}.
\end{equation}
Then we obtain (\ref{p_+}) and (\ref{q_+}) by letting
$$\left(\begin{array}{c}
p_+(\xi,\eta)\\[1mm]
q_+(\xi,\eta)
   \end{array}\right):=\left(\begin{array}{c}
   (\Theta^{-2}\cdot(e^{\frac{\rm i}2\theta}+ A)-e^{\frac{\rm i}2\theta_+})\eta\\[1mm]
    (\Theta^{2}\cdot(e^{\frac{\rm i}2\theta}+ A)^{-1}-e^{-\frac{\rm i}2\theta_+})\xi
   \end{array}\right)+\left(\begin{array}{c}
    \Theta^{-1}\cdot p(\Theta\xi,\Theta^{-1}\eta)\\[1mm]
    \Theta\cdot q(\Theta\xi,\Theta^{-1}\eta)
   \end{array}\right).$$
Since $\| A\|_{\CO,\beta,r}\leq \frac1{16}$, we see that $(e^{\frac{\rm i}2\theta}+ A)^{\pm1}\in{\CA}_{\beta',r'}\left({\CO}\right)$.
According to Lemma \ref{lemma_analytic_prod}, \rl{lemma_Lambda_k} and \rl{lemma_h_Lambda}, we have that $p_+, \, q_+\in{\CA}_{\beta',  r'}({\CO})$.

It remains to prove \re{crossing_theta_+}.
With the above $p_+$ and $q_+$, we have that
\begin{eqnarray}
& & e^{-\frac{\rm i}2\theta_+}\xi p_++e^{\frac{\rm i}2\theta_+}\eta q_+\nonumber\\
&=&e^{-\frac{\rm i}2\theta_+}\cdot (\xi\eta) \left(\Theta^{-2}\cdot(e^{\frac{\rm i}2\theta}+ A)-e^{\frac{\rm i}2\theta_+}\right)\label{crossing_part01}\\
& & + \, e^{-\frac{\rm i}2\theta_+}\cdot (\xi\eta)  \left(\Theta^{2}\cdot(e^{\frac{\rm i}2\theta}+ A)^{-1}-e^{-\frac{\rm i}2\theta_+}\right)\label{crossing_part02}\\
& & + \, e^{-\frac{\rm i}2\theta_+}\xi  \cdot \Theta^{-1}\cdot p(\Theta\xi,\Theta^{-1}\eta)+ e^{\frac{\rm i}2\theta_+}\eta \cdot  \Theta\cdot q(\Theta\xi,\Theta^{-1}\eta).\label{crossing_part1}
\end{eqnarray}
Since \re{theta_+} implies $\|\theta_+-\theta\|_{\CO,\beta',r'} <\frac{101}{50}\| A\|_{\CO,\beta,r}$, so that
$$\|e^{{\rm i}b\theta_+}-e^{{\rm i}b\theta}\|_{\CO,\beta',r'}< \frac72 |b| \| A\|_{\CO,\beta,r},\quad -1\leq b \leq 1,$$
which implies that $\|e^{{\rm i}b\theta_+}\|_{\CO,\beta',r'}<2$,
we see that, in view of \re{err_fre_1} and \re{err_fre_2}, the sum of terms in (\ref{crossing_part01}) and (\ref{crossing_part02}) is smaller than $\| A\|_{\CO,\beta,r}^2$.
Let us focus on (\ref{crossing_part1}), which equals to
\begin{eqnarray*}
& &e^{-\frac{\rm i}2\theta_+}  \cdot \Theta^{-1}\xi\cdot p(\Theta\xi,\Theta^{-1}\eta)+ e^{\frac{\rm i}2\theta_+} \cdot  \Theta \eta \cdot q(\Theta\xi,\Theta^{-1}\eta)\\
&=&e^{-\frac{\rm i}2\theta}  \cdot (\Theta\xi) p(\Theta\xi,\Theta^{-1}\eta)+ e^{\frac{\rm i}2\theta} \cdot  (\Theta^{-1} \eta ) q(\Theta\xi,\Theta^{-1}\eta)\\
& & + \,  (e^{-\frac{\rm i}2\theta_+} \Theta^{-1}-e^{-\frac{\rm i}2\theta} \Theta)\xi p(\Theta\xi,\Theta^{-1}\eta)+(e^{\frac{\rm i}2\theta_+}\Theta-e^{\frac{\rm i}2\theta}\Theta^{-1} )\eta q(\Theta\xi,\Theta^{-1}\eta)
\end{eqnarray*}
Since Lemma \ref{lemma_Lambda_k} implies that $\|\Theta^{-1}\|_{\CO,\beta',r'}<1+\frac{3}{4}\| A\|_{\CO,\beta,r}<\frac{19}{16}$ and
$$\|\Theta-\Theta^{-1}\|_{\CO,\beta',r'}\leq\|\Theta-1\|_{\CO,\beta',r'}+\|\Theta^{-1}-1\|_{\CO,\beta',r'}<\frac{3}{2}\| A\|_{\CO,\beta,r},$$
we obtain that
\begin{eqnarray*}
\|e^{-\frac{\rm i}2\theta_+} \Theta^{-1}-e^{-\frac{\rm i}2\theta} \Theta\|_{\CO,\beta',r'}
&\leq&\|(e^{-\frac{\rm i}2\theta_+} - e^{-\frac{\rm i}2\theta})\Theta^{-1}\|_{\CO,\beta',r'}+\|e^{-\frac{\rm i}2\theta}
 (\Theta^{-1}-\Theta)\|_{\CO,\beta',r'}\\
&<&\frac{19}{16}\cdot\frac{7}{4}\| A\|_{\CO,\beta,r}+\frac{101}{100}\cdot\frac{3}{2}\| A\|_{\CO,\beta,r}<4\| A\|_{\CO,\beta,r}.
\end{eqnarray*}
Similarly, we have
$$\|e^{\frac{\rm i}2\theta_+} \Theta-e^{\frac{\rm i}2\theta} \Theta^{-1}\|_{\CO,\beta',r'}<4\| A\|_{\CO,\beta,r}.$$
On the other hand, by Lemma \ref{lemma_h_Lambda}, we see that
$$\|e^{-\frac{\rm i}2\theta}  \cdot (\Theta\xi) p(\Theta\xi,\Theta^{-1}\eta)+ e^{\frac{\rm i}2\theta} \cdot  (\Theta^{-1} \eta ) q(\Theta\xi,\Theta^{-1}\eta)\|_{\CO,\beta',r'}<\|e^{-\frac{\rm i}2\theta} \xi p+ e^{\frac{\rm i}2\theta} \eta  q\|_{\CO,\beta,r}.$$
Hence (\ref{crossing_part1}) is bounded by $\|e^{-\frac{\rm i}2\theta} \xi p+ e^{\frac{\rm i}2\theta} \eta  q\|_{\CO,\beta',r'}+\| A\|_{\CO,\beta,r}(\|p\|_{\CO,\beta,r}+\|q\|_{\CO,\beta,r})$. Combining with the errors in (\ref{crossing_part01}) and (\ref{crossing_part02}), \re{crossing_theta_+} is shown.\qed

\subsection{Approximated cohomological equations -- Proof of Theorem \ref{thm_trans_cohomo}}\label{subsec_cohomo}

For $0<r_+<r<\frac14$, let $\varepsilon>0$ be sufficiently small such that (\ref{varepsilon_small-general}) is satisfied,
let $\beta$, $\beta_+$, $\tilde\beta$ be defined as in (\ref{new_def_beta}), and let $r^{(m)}$, $m=0,1,\cdots,8$, be defined as in (\ref{r_m}).
Consider the holomorphic involution
$$\tau_1(\xi,\eta)=\left(\begin{array}{l}
e^{\frac{\rm i}2\alpha(\xi\eta)}\eta+ p(\xi,\eta) \\
e^{-\frac{\rm i}2\alpha(\xi\eta)}\xi+ q(\xi,\eta)
\end{array}
\right),\quad (\xi,\eta)\in {\cL C}^r_{\omega,\beta}, \quad \omega\in\CO(r,\beta)$$
given as in Section \ref{sec_rev} (same with that in Theorem \ref{thm_trans_cohomo}), with $\|p\|_{\CO,\beta,r}$, $\|q\|_{\CO,\beta,r}<\frac{\varepsilon}{10}$.

The rest of this subsection is devoted to the proof of Theorem \ref{thm_trans_cohomo}. The core of proof is the resolution of the approximated cohomological equations (see Lemma \ref{lem_cohomo}).

%

\smallskip

At first, we see that the definition of $\CO_{\delta}$ in (\ref{SD_general}) implies that
\begin{equation}\label{SD_general2}
|e^{{\rm i}n\alpha(\xi\eta)}-1|\geq \frac{\delta}2,\quad \forall \ 0<|n|\leq K+1,\quad \forall \ (\xi,\eta)\in {\cL C}_{\omega,\tilde\beta}, \ \omega\in\CO_{\delta}.
\end{equation}
Indeed, recalling that $\tilde\beta=16\beta_+$ and $\delta>80 \varepsilon^{\frac1{60s}}$, by 
(\ref{varepsilon_small-general}), we have
\begin{equation}\label{K_beta}
K+1=\frac{|\ln\varepsilon|}{\left|\ln(r^{(7)}/r)\right|}+1= \frac{|\ln\varepsilon|}{\left|\ln(\frac78+\frac{r_+}{8r})\right|}+1< \frac{\varepsilon^{-\frac{1}{2400s^2}}}2.
\end{equation}
Hence, by Lemma \ref{lem_a0}, for $0<|n|\leq K+1$, we have
$$|n| |\alpha(\xi\eta)-\alpha(\omega)|_{\omega,\tilde\beta}< (K+1) \cdot 18\beta_+
\leq 18(K+1)\varepsilon^{\frac1{48s}}<20 \varepsilon^{\frac1{60s}}<\frac{\delta}{4}.$$
In view of Remark \ref{rem_eiba}, we have
$$|e^{{\rm i}n\alpha(\xi\eta)}-e^{{\rm i}n\alpha(\omega)}|_{\omega,\tilde\beta} \leq e^{\frac98|n|\tilde\beta}|e^{{\rm i}n(\alpha(\xi\eta)-\alpha(\omega))}-1|_{\omega,\tilde\beta}
\leq  \frac{101}{100}\sum_{k\geq 1} \frac{|n|^k|\alpha(\xi\eta)-\alpha(\omega)|^k_{\omega,\tilde\beta}}{k!}<\frac{\delta}{2},$$
and hence, combining with \re{SD_general}, (\ref{SD_general2}) is obtained.

\

Define $p_K$ and $q_K\in{\CA}_{\beta,r}(\CO)$ as
\beq\label{pK}\left(\begin{array}{l}
p_K\\
q_K
\end{array}\right):=\left(\begin{array}{l}
p_{0,0}(\xi\eta)+\sum_{1\leq l\leq K} p_{l,0}(\xi\eta)\xi^l+\sum_{1\leq j\leq K} p_{0,j}(\xi\eta)\eta^j\\
q_{0,0}(\xi\eta)+\sum_{1\leq l\leq K} q_{l,0}(\xi\eta)\xi^l+\sum_{1\leq j\leq K} q_{0,j}(\xi\eta)\eta^j
\end{array}\right).
\eeq
Since $\|p\|_{\CO,\beta, r}$, $\|q\|_{\CO,\beta, r}< \frac\varepsilon{10}$, we have
\begin{equation}\label{rest_pq}
\|p-p_{K}\|_{\CO,\tilde\beta,r^{(7)}},\  \|q-q_K\|_{\CO,\tilde\beta,r^{(7)}}
< \frac\varepsilon{10}  \left(\frac{r^{(7)}}{r}\right)^{\frac{|\ln\varepsilon|}{\left|\ln(r^{(7)}/r)\right|}}=\frac{\varepsilon^2}{10}.
\end{equation}

\begin{lemma}\label{lem_cohomo}
There is $\hat{\cL U}=\left(\begin{array}{c}
                                                              \hat u \\
                                                              \hat v
                                                            \end{array}
\right)\in({\CA}^\R_{\tilde\beta,r^{(7)}}(\CO_\delta))^2$, with $\hat u_{1,0}=\hat v_{0,1}=0$ satisfying
\begin{eqnarray}
\left\|\hat u\right\|_{\CO_{\delta} , \tilde\beta,r^{(7)}}, \left\|\hat v\right\|_{\CO_{\delta},\tilde\beta,r^{(7)}} &<& \frac{\varepsilon^{\frac{49}{50}}}{20},\label{esti_uv}\\
\left\|\eta \hat u+\xi \hat v\right\|_{\CO_{\delta},\tilde\beta,r^{(7)}}&<&\frac{\varepsilon^{\frac{61}{32}}}{16}+5(K+1)\delta^{-1}\|e^{\frac{\rm i}2\alpha}\eta q+e^{-\frac{\rm i}2\alpha}\xi p\|_{\CO,\beta,r}\label{esti_Luv}
\end{eqnarray}
 such that
 \begin{equation}\label{cohomo1}
 \|e^{\frac{\rm i}2\alpha} \hat v-\hat u(e^{\frac{\rm i}{2}\alpha}\eta,e^{-\frac{\rm i}{2}\alpha}\xi)+p_{K}-p_{0,1}\eta \|_{\CO_{\delta},\tilde\beta,r^{(7)}}
 < \frac{\varepsilon^{\frac{61}{32}}}{80}+\frac{6(K+1)}{\delta}\|e^{\frac{\rm i}2\alpha}\eta q+e^{-\frac{\rm i}2\alpha}\xi p\|_{\CO,\beta, r},
 \end{equation}
  \begin{equation}\label{cohomo2}
  \|e^{-\frac{\rm i}2\alpha} \hat u-\hat v(e^{\frac{\rm i}{2}\alpha}\eta,e^{-\frac{\rm i}{2}\alpha}\xi)+q_{K}-q_{1,0}\xi\|_{\CO_{\delta},\tilde\beta,r^{(7)}}
 < \frac{\varepsilon^{\frac{61}{32}}}{80}+\frac{6(K+1)}{\delta}\|e^{\frac{\rm i}2\alpha}\eta q+e^{-\frac{\rm i}2\alpha}\xi p\|_{\CO,\beta, r},
 \end{equation}
   \begin{eqnarray}
  & & \left\|e^{-\frac{\rm i}{2}\alpha}\xi\left(e^{\frac{\rm i}2\alpha} \hat v-\hat u(e^{\frac{\rm i}{2}\alpha}\eta,e^{-\frac{\rm i}{2}\alpha}\xi)+p_{K}-p_{0,1}\eta\right)\right. \label{crossing_cohomo}\\
    & & \  \  \  \   \left.+ \, e^{\frac{\rm i}{2}\alpha}\eta \left( e^{-\frac{\rm i}2\alpha} \hat u-\hat v(e^{\frac{\rm i}{2}\alpha}\eta,e^{-\frac{\rm i}{2}\alpha}\xi)+q_{K}-q_{1,0}\xi\right)\right\|_{\CO_{\delta},\tilde\beta,r^{(7)}}<\frac{\varepsilon^{\frac{61}{32}}}{20}.\nonumber
   \end{eqnarray}
\end{lemma}

\begin{remark} As mentioned in Section 2.1, the aim of the KAM-like process is to eliminate the main part of the perturbation and get a much smaller new perturbation. 

In view of \re{rest_pq}, we would like to construct $\hat{\cL U}=\left(\begin{array}{c}
                                                              \hat u \\
                                                              \hat v
                                                            \end{array}
\right)$ so that the change of variables ${\rm Id}+\hat{\cL U}$ eliminates $p_{K}$ and $q_{K}$. This amounts to solve the cohomological equations (arising as the ``linearized equation" of the conjugacy equation),
\begin{eqnarray*}
e^{\frac{\rm i}2\alpha} \hat v-\hat u(e^{\frac{\rm i}{2}\alpha}\eta,e^{-\frac{\rm i}{2}\alpha}\xi)+p_{K}-p_{0,1}\eta&=&0, \\
e^{-\frac{\rm i}2\alpha} \hat u-\hat v(e^{\frac{\rm i}{2}\alpha}\eta,e^{-\frac{\rm i}{2}\alpha}\xi)+q_{K}-q_{1,0}\xi&=&0,
\end{eqnarray*}
with $p_{0,1}\eta$ and $q_{1,0}\xi$ added to the new principal part.
Based on the fact that $\tau_1$ and $\tau_2$ are involutions and $\tau_1\circ\tau_2$ is reversible w.r.t. $\rho$,
we solve the above equations approximately, with the errors estimated as in (\ref{cohomo1}) and (\ref{cohomo2}).
Then, with the transformation ${\rm Id}+\hat{\CU}$, the main parts of perturbations of $\tau_1$, $\tau_2$ and $\tau_1\circ\tau_2$ are all eliminated approximately.
\end{remark}

\proof Let us define $\hat u, \hat v$ by giving the coefficients: $\hat u_{1,0}=\hat v_{0,1}=0$,
\begin{eqnarray}
\hat u_{l,0}(\xi\eta) &:=&\frac12\cdot\frac{f_{l,0}(\xi\eta)-e^{{\rm i}(l+1)\alpha(\xi\eta)}\bar f_{l,0}(\xi\eta)}{e^{{\rm i}l \alpha(\xi\eta)}-e^{{\rm i}\alpha(\xi\eta)}}  , \quad 2\leq l\leq K,\label{ul} \\
\hat u_{0,j}(\xi\eta) &:=&\frac12\cdot\frac{f_{0,j}(\xi\eta)-e^{-{\rm i}(j-1)\alpha(\xi\eta)}\bar f_{0,j}(\xi\eta)}{e^{-{\rm i}j\alpha(\xi\eta)}-e^{{\rm i}\alpha(\xi\eta)}}, \quad 0\leq j\leq K,\label{uj}\\
\hat v_{l,0}(\xi\eta) &:=&\frac12\cdot\frac{g_{l,0}(\xi\eta)-e^{{\rm i}(l-1)\alpha(\xi\eta)}\bar g_{l,0}(\xi\eta)}{e^{{\rm i}l \alpha(\xi\eta)}-e^{-{\rm i}\alpha(\xi\eta)}}  , \quad 0\leq l\leq K, \label{vl}\\
\hat v_{0,j}(\xi\eta) &:=&\frac12\cdot\frac{g_{0,j}(\xi\eta)-e^{-{\rm i}(j+1)\alpha(\xi\eta)}\bar g_{0,j}(\xi\eta)}{e^{-{\rm i}j\alpha(\xi\eta)}-e^{-{\rm i}\alpha(\xi\eta)}}, \quad 2\leq j\leq K,\label{vj}
\end{eqnarray}
with other coefficients being $0$. Here, $f,g$ are defined by \re{f} and \re{g}.
In view of (\ref{ul}) -- (\ref{vj}), we see that if $\hat u_{l,j}\neq 0$, then
$$\overline{\hat u}_{l,j}(\xi\eta)=\frac12\cdot\frac{\bar f_{l,j}(\xi\eta)-e^{-{\rm i}(l-j+1)\alpha(\xi\eta)} f_{l,j}(\xi\eta)}{e^{-{\rm i}(l-j) \alpha(\xi\eta)}-e^{-{\rm i}\alpha(\xi\eta)}}=\hat u_{l,j}(\xi\eta),$$
and if $\hat v_{l,j}\neq 0$, then
$$\overline{\hat v}_{l,j}(\xi\eta)=\frac12\cdot \frac{\bar g_{l,j}(\xi\eta)-e^{-{\rm i}(l-j-1)\alpha(\xi\eta)} g_{l,j}(\xi\eta)}{e^{-{\rm i}(l-j) \alpha(\xi\eta)}-e^{{\rm i}\alpha(\xi\eta)}}= \hat v_{l,j}(\xi\eta).$$
Hence, according to the definition of $\CO_{\delta}$ in (\ref{SD_general}) and Lemma \ref{lem_appro}, $\hat u, \hat v\in{\CA}_{\tilde\beta,r^{(7)}}^\R\left(\CO_{\delta}\right)$.

With the coefficients in (\ref{ul}) -- (\ref{vj}),  we have, for $lj= 0$, $l,j\leq K$ and $(l,j)\neq (0,1)$,
\begin{equation}\label{evu}
e^{\frac{\rm i}{2}\alpha}\hat v_{l,j}-e^{\frac{\rm i}{2}(j-l)\alpha}\hat u_{j,l}=\frac{e^{\frac{\rm i}{2}\alpha}}{2}\frac{g_{l,j}-e^{{\rm i}(l-j-1)\alpha} \bar g_{l,j}}{e^{{\rm i}(l-j)\alpha}-e^{-{\rm i}\alpha}}-\frac{e^{\frac{\rm i}{2}(j-l)\alpha}}{2}\frac{f_{j,l}-e^{{\rm i}(j-l+1)\alpha}\bar f_{j,l}}{e^{{\rm i}(j-l)\alpha}-e^{{\rm i}\alpha}}, \end{equation}
and for $lj= 0$, $l,j\leq K$, $(l,j)\neq (1,0)$,
\begin{equation}\label{euv}
e^{-\frac{\rm i}{2}\alpha}\hat u_{l,j}-e^{\frac{\rm i}{2}(j-l)\alpha}\hat v_{j,l}= \frac{e^{-\frac{\rm i}{2}\alpha}}{2}\frac{f_{l,j}-e^{{\rm i}(l-j+1)\alpha} \bar f_{l,j}}{e^{{\rm i}(l-j)\alpha}-e^{{\rm i}\alpha}}-\frac{e^{\frac{\rm i}{2}(j-l)\alpha}}{2}\frac{g_{j,l}-e^{{\rm i}(j-l-1)\alpha}\bar g_{j,l}}{e^{{\rm i}(j-l)\alpha}-e^{-{\rm i}\alpha}}.
\end{equation}
According to Corollary \ref{cor_fg_coeff}, we have
$$C=\frac{\rm i}{2}\alpha'(\xi\eta)\left(e^{-\frac{\rm i}{2}\alpha(\xi\eta)}\eta \bar q +e^{\frac{\rm i}{2}\alpha(\xi\eta)}\xi \bar p \right) \in {\CA}_{\tilde\beta,r^{(7)}}(\CO_{\delta}).$$
Replacing the coefficients of $f$ and $g$ in (\ref{evu}) and (\ref{euv}) and according to (\ref{appro_f_coeff_l}) -- (\ref{appro_g_coeff_j}), let us show that, for $lj= 0$ with $l,j\leq K$,
\begin{align}
\|e^{\frac{\rm i}{2}\alpha}\hat v_{l,j}-e^{\frac{\rm i}{2}(j-l)\alpha}\hat u_{j,l}+p_{l,j}-\hat p_{l,j}\|_{\CO_{\delta},\tilde\beta,r^{(7)}} &< \frac{\delta^{-1} \varepsilon^{\frac{31}{16}}}{16 (r^{(7)})^{l+j}},\quad (l,j)\neq (0,1),\label{appro_hat_p} \\
\|e^{-\frac{\rm i}{2}\alpha}\hat u_{l,j}-e^{\frac{\rm i}{2}(j-l)\alpha}\hat v_{j,l}+q_{l,j}-\hat q_{l,j}\|_{\CO_{\delta},\tilde\beta,r^{(7)}} &<\frac{\delta^{-1}  \varepsilon^{\frac{31}{16}}}{16 (r^{(7)})^{l+j}},\quad  (l,j)\neq (1,0).\label{appro_hat_q}
\end{align}
Here, $\hat p, \hat q\in \CA_{\tilde\beta,r^{(7)}}(\CO_{\delta})$ with $\hat p_{l,j}$ and $\hat q_{l,j}$  defined by
\begin{eqnarray}
\hat p_{l,0} &:=&(\xi\eta) \cdot \frac{-e^{-\frac{\rm i}2\alpha}C_{l+1,0}+e^{{\rm i}l\alpha}e^{\frac{\rm i}2\alpha}\bar C_{l+1,0}+e^{\frac{\rm i}2l\alpha}C_{0,l+1}-e^{-\frac{\rm i}2l\alpha}e^{-{\rm i}\alpha}\bar C_{0,l+1}}{2(e^{{\rm i}l\alpha}-e^{-{\rm i}\alpha})}\label{hat_p_l}\\
& & + \, \frac12\left(p_{l,0}+e^{-\frac{\rm i}{2}(l-1)\alpha} q_{0,l} \right),\quad 0\leq l \leq K,\nonumber
\end{eqnarray}
\begin{eqnarray}
\hat p_{0,j} &:=&\frac{-e^{-\frac{\rm i}2\alpha}C_{0,j-1}+e^{-{\rm i}j\alpha}e^{\frac{\rm i}2\alpha}\bar C_{0,j-1}+e^{-\frac{\rm i}2j\alpha}C_{j-1,0}-e^{\frac{\rm i}2j\alpha}e^{-{\rm i}\alpha}\bar C_{j-1,0}}{2 (e^{-{\rm i}j\alpha}-e^{-{\rm i}\alpha})} \label{hat_p_j}\\
& & + \, \frac12\left(p_{0,j}+e^{\frac{\rm i}{2}(j+1)\alpha} q_{j,0} \right),\quad 2\leq j\leq K, \nonumber
\end{eqnarray}
\begin{eqnarray}
\hat q_{l,0}&:=&\frac{e^{\frac{\rm i}2\alpha}C_{l-1,0}-e^{{\rm i}l\alpha}e^{-\frac{\rm i}2\alpha}\bar C_{l-1,0}-e^{\frac{\rm i}2l\alpha}C_{0,l-1}+e^{-\frac{\rm i}2l\alpha}e^{{\rm i}\alpha}\bar C_{0,l-1}}{2(e^{{\rm i}l\alpha}-e^{{\rm i}\alpha})}\label{hat_q_l}\\
& & + \, \frac12\left(q_{l,0}+e^{-\frac{\rm i}{2}(l+1)\alpha} p_{0,l} \right),\quad 2\leq l\leq K, \nonumber
\end{eqnarray}
\begin{eqnarray}
\hat q_{0,j}&:=& (\xi\eta) \cdot\frac{e^{\frac{\rm i}2\alpha}C_{0,j+1}-e^{-{\rm i}j\alpha}e^{-\frac{\rm i}2\alpha}\bar C_{0,j+1}-e^{-\frac{\rm i}2j\alpha}C_{j+1,0}+e^{\frac{\rm i}2j\alpha}e^{{\rm i}\alpha}\bar C_{j+1,0}}{2(e^{-{\rm i}j\alpha}-e^{{\rm i}\alpha})}\label{hat_q_j}\\
& & + \, \frac12\left(q_{0,j}+e^{\frac{\rm i}{2}(j-1)\alpha} p_{j,0} \right),\quad 0\leq j \leq K,\nonumber
\end{eqnarray}
and all other coefficients of $\hat p$ and $\hat q$ are defined as $0$.
Indeed, by (\ref{evu}), we have
$$e^{\frac{\rm i}{2}\alpha}\hat v_{l,0}-e^{-\frac{\rm i}{2}l\alpha}\hat u_{0,l}
=  \frac{e^{\frac{\rm i}{2}\alpha}}{2}\frac{g_{l,0}-e^{{\rm i}(l-1)\alpha} \bar g_{l,0}}{e^{{\rm i}l\alpha}-e^{-{\rm i}\alpha}}-\frac{e^{-\frac{\rm i}{2}l\alpha}}{2}\frac{f_{0,l}-e^{-{\rm i}(l-1)\alpha}\bar f_{0,l}}{e^{-{\rm i}l\alpha}-e^{{\rm i}\alpha}}
=\sum_{k=1}^4 ({\CM}_k+{\CN}_k),$$
where, 
${\CM}_2:=\bar {\CM}_1e^{{\rm i}\alpha}$,  ${\CM}_4:=\bar {\CM}_3e^{-{\rm i}l\alpha}$, ${\CN}_2:=\bar {\CN}_1e^{{\rm i}\alpha}$,  ${\CN}_4:=\bar {\CN}_3e^{-{\rm i}l\alpha}$ together with\begin{eqnarray*}
{\CM}_1 &:=&   \frac{e^{\frac{3{\rm i}}{2}\alpha}}{2(e^{{\rm i}(l+1)\alpha}-1)}\left( e^{-\frac{\rm i}{2}\alpha} \bar p_{l,0}+e^{\frac{\rm i}{2}l\alpha} q_{0,l}-(\xi\eta)e^{-{\rm i}\alpha} C_{l+1,0}\right),\\
{\CN}_1 &:=&  \frac{e^{\frac{3{\rm i}}{2}\alpha}}{2(e^{{\rm i}(l+1)\alpha}-1)} \left(g_{l,0}-e^{-\frac{\rm i}{2}\alpha} \bar p_{l,0}-e^{\frac{\rm i}{2}l\alpha} q_{0,l}+(\xi\eta)e^{-{\rm i}\alpha} C_{l+1,0}\right),\\
{\CM}_3 &:=& - \, \frac{e^{-\frac{\rm i}{2}(l+2)\alpha}}{2(e^{-{\rm i}(l+1)\alpha}-1)} \left( e^{\frac{\rm i}{2}\alpha} \bar q_{0,l}+e^{-\frac{\rm i}{2}l\alpha} p_{l,0}+(\xi\eta)e^{{\rm i}\alpha}C_{0,l+1}\right),\\
{\CN}_3  &:=& - \, \frac{e^{-\frac{\rm i}{2}(l+2)\alpha}}{2(e^{-{\rm i}(l+1)\alpha}-1)}  \left(f_{0,l}-e^{\frac{\rm i}{2}\alpha} \bar q_{0,l}-e^{-\frac{\rm i}{2}l\alpha} p_{l,0}-(\xi\eta)e^{{\rm i}\alpha}C_{0,l+1} \right).
\end{eqnarray*}
By Corollary \ref{lem_eiba} and Remark \ref{rem_eiba}, we obtain, for $0\leq l \leq K$, $\omega\in \CO(r^{(7)},\tilde\beta)$,
$$|e^{\frac{3{\rm i}}{2}\alpha}|_{\omega,\tilde\beta},\ |e^{{\rm i}(l+\frac12)\alpha}|_{\omega,\tilde\beta}, \ |e^{-\frac{\rm i}{2}(l+2)\alpha}|_{\omega,\tilde\beta}, \ |e^{-\frac{3{\rm i}}2l\alpha}|_{\omega,\tilde\beta}<\frac{101}{100},$$
and, by (\ref{SD_general2}), for $0\leq l \leq K$,$\left\|(e^{\pm{\rm i}(l+1)\alpha}-1)^{-1}\right\|_{\CO_{\delta},\tilde\beta,r^{(7)}}<2\delta^{-1}$.
Then, by (\ref{appro_f_coeff_j}), (\ref{appro_g_coeff_l}), we have
\begin{equation}
 \sum_{k=1}^4\|{\CN}_k\|_{\CO_{\delta}, \tilde\beta, r^{(7)}}
\leq 4 \cdot\frac{101}{100}\delta^{-1}\cdot\frac{\varepsilon^{\frac{31}{16}}}{80(r^{(7)})^{l}}
<\frac{\delta^{-1} \varepsilon^{\frac{31}{16}}}{16 (r^{(7)})^l}.\label{com_esti_1}
\end{equation}
On the other hand, we have
\begin{eqnarray*}
{\CM}_1&=&\frac{e^{\frac{{\rm i}}{2}\alpha}}{2(e^{{\rm i}l\alpha}-e^{-{\rm i}\alpha})}\left( e^{-\frac{\rm i}{2}\alpha} \bar p_{l,0}+e^{\frac{\rm i}{2}l\alpha} q_{0,l}-(\xi\eta)e^{-{\rm i}\alpha} C_{l+1,0}\right)\\
&=&\frac{1}{2(e^{{\rm i}l\alpha}-e^{-{\rm i}\alpha})}\left( \bar p_{l,0}+e^{\frac{\rm i}{2}(l+1)\alpha} q_{0,l}-(\xi\eta)e^{-\frac{\rm i}2\alpha} C_{l+1,0}\right),\\
{\CM}_2
&=& - \, \frac{1}{2(e^{{\rm i}l\alpha}-e^{-{\rm i}\alpha})}\left( e^{{\rm i}l\alpha} p_{l,0}+e^{\frac{\rm i}{2}(l-1)\alpha} \bar q_{0,l}-(\xi\eta)e^{{\rm i}(l+\frac12)\alpha}\bar C_{l+1,0} \right),\\
{\CM}_3
&=& \frac{1}{2(e^{{\rm i}l\alpha}-e^{-{\rm i}\alpha})}\left( e^{\frac{\rm i}{2}(l-1)\alpha} \bar q_{0,l}+e^{-{\rm i}\alpha} p_{l,0}+(\xi\eta)e^{\frac{\rm i}2l\alpha}C_{0,l+1}\right),\\
{\CM}_4
&=&- \, \frac{1}{2(e^{{\rm i}l\alpha}-e^{-{\rm i}\alpha})}\left(e^{-\frac{\rm i}2(l+1)\alpha} q_{0,l}+ \bar p_{l,0}+(\xi\eta)e^{-\frac{\rm i}2(l+2)\alpha}\bar C_{0,l+1} \right).
\end{eqnarray*}
Hence, adding the above terms, we have
\begin{eqnarray}
 p_{l,0}+ \sum_{k=1}^4 {\CM}_k
&=& p_{l,0}+\frac{1}{2(e^{{\rm i}l\alpha}-e^{-{\rm i}\alpha})} \left[\left(\bar p_{l,0} -\bar p_{l,0}\right)+\left(e^{\frac{\rm i}{2}(l-1)\alpha} \bar q_{0,l}-e^{\frac{\rm i}{2}(l-1)\alpha} \bar q_{0,l}\right) \right] \nonumber\\
& & + \, \frac{(e^{-{\rm i}\alpha} p_{l,0}-e^{{\rm i}l\alpha} p_{l,0})+(e^{\frac{\rm i}{2}(l+1)\alpha} q_{0,l}-e^{-\frac{\rm i}2(l+1)\alpha} q_{0,l})}{2(e^{{\rm i}l\alpha}-e^{-{\rm i}\alpha})}\nonumber\\
& &  + \, (\xi\eta)\cdot\frac{e^{\frac{\rm i}2l\alpha}C_{0,l+1}-e^{-\frac{\rm i}2l\alpha}e^{-{\rm i}\alpha}\bar C_{0,l+1}-e^{-\frac{\rm i}2\alpha} C_{l+1,0} +e^{{\rm i}l\alpha}e^{\frac{\rm i}2\alpha}\bar C_{l+1,0}}{2(e^{{\rm i}l\alpha}-e^{-{\rm i}\alpha})}\nonumber\\
&=& \frac{p_{l,0}+e^{-\frac{\rm i}{2}(l-1)\alpha}q_{0,l}}{2}\nonumber \\
& & + \, (\xi\eta) \cdot \frac{-e^{-\frac{\rm i}2\alpha}C_{l+1,0}+e^{{\rm i}l\alpha}e^{\frac{\rm i}2\alpha}\bar C_{l+1,0}+e^{\frac{\rm i}2l\alpha}C_{0,l+1}-e^{-\frac{\rm i}2l\alpha}e^{-{\rm i}\alpha}\bar C_{0,l+1}}{2(e^{{\rm i}l\alpha}-e^{-{\rm i}\alpha})} \nonumber\\
&=& \hat p_{l,0}. \label{com_esti_2}
\end{eqnarray}
By \re{com_esti_1} and \re{com_esti_2}, we obtain \re{appro_hat_p} for the case $(l,j)=(l,0)$, $l\geq 0$. The proof is similar for the case $(l,j)=(0,j)$, $j\geq 2$. It is also similar for proving \re{appro_hat_q}.

Now we are going to show \re{cohomo1} -- \re{crossing_cohomo}. Since \re{appro_hat_p} and \re{appro_hat_q} imply that
\begin{eqnarray*}
\|e^{\frac{\rm i}2\alpha} \hat v-\hat u(e^{\frac{\rm i}{2}\alpha}\eta,e^{-\frac{\rm i}{2}\alpha}\xi)+p_{K}-p_{0,1}\eta -\hat p \|_{\CO_{\delta},\tilde\beta,r^{(7)}}&<&\frac{K\delta^{-1} \varepsilon^{\frac{31}{16}}}{8},\\
\|e^{-\frac{\rm i}2\alpha} \hat u-\hat v(e^{\frac{\rm i}{2}\alpha}\eta,e^{-\frac{\rm i}{2}\alpha}\xi)+q_{K}-q_{1,0}\xi -\hat q\|_{\CO_{\delta},\tilde\beta,r^{(7)}}&<&\frac{K\delta^{-1} \varepsilon^{\frac{31}{16}}}{8},
\end{eqnarray*}
it is sufficient to prove
\begin{align}
\|\hat p\|_{\CO_{\delta},\tilde\beta,r^{(7)}},\ \|\hat q\|_{\CO_{\delta},\tilde\beta,r^{(7)}}&<\frac{\varepsilon^{\frac{61}{32}}}{100}+6(K+1)\delta^{-1}\|e^{\frac{\rm i}2\alpha}\eta q+e^{-\frac{\rm i}2\alpha}\xi p\|_{\CO,\beta,r} ,\label{esti_hat_pq}\\
 \|e^{-\frac{\rm i}{2}\alpha}\xi \hat p +e^{\frac{\rm i}{2}\alpha}\eta \hat q\|_{\CO_{\delta},\tilde\beta,r^{(7)}}&<\frac{\varepsilon^{\frac{61}{32}}}{40} \label{crossing_hat_pq}.
\end{align}
By (\ref{esi_coeff}), Lemma \ref{lem_norm00} and Corollary \ref{cor_deri_a}, we have
$$
\|C\|_{\CO_{\delta},\tilde\beta,r^{(7)}}\leq \frac{\|\alpha'\|_{\CO,\beta,r}}{2} \|e^{\frac{\rm i}2\alpha}\eta q+e^{-\frac{\rm i}2\alpha}\xi p\|_{\CO,\tilde\beta,r^{(7)}}<\frac35 \|e^{\frac{\rm i}2\alpha}\eta q+e^{-\frac{\rm i}2\alpha}\xi p\|_{\CO,\beta,r},
$$
which gives the estimates for the coefficients: for $lj=0$,
$$\|C_{l,j}\|_{\CO_{\delta},\tilde\beta,r^{(7)}}< \frac{3}{5(r^{(7)})^{l+j}} \|e^{\frac{\rm i}2\alpha}\eta q+e^{-\frac{\rm i}2\alpha}\xi p\|_{\CO,\beta,r}.$$
Then, in view of (\ref{appro_p_}), we have that
\begin{eqnarray*}
\|e^{\frac{\rm i}{2}\alpha}q+p(e^{\frac{\rm i}{2}\alpha}\eta, e^{-\frac{\rm i}{2}\alpha}\xi)\|_{\CO_{\delta},\tilde\beta,r^{(7)}}
&\leq& \|e^{\frac{\rm i}{2}\alpha}q+p(e^{\frac{\rm i}{2}\alpha}\eta, e^{-\frac{\rm i}{2}\alpha}\xi)-\bar C \xi\|_{\CO_{\delta},\tilde\beta,r^{(7)}}+\|C\xi\|_{\CO_{\delta},\tilde\beta,r^{(7)}} \\
&\leq&  \frac{\varepsilon^{\frac{31}{16}}}{80}+\frac{3r^{(7)}}{5}\|e^{\frac{\rm i}2\alpha}\eta q+e^{-\frac{\rm i}2\alpha}\xi p\|_{\CO,\beta,r},
\end{eqnarray*}
which gives the estimates for the coefficients: for $lj=0$,
\begin{eqnarray*}
& & \|e^{\frac{\rm i}{2}\alpha}q_{l,j}+e^{\frac{\rm i}{2}(j-l)\alpha}p_{j,l}\|_{\CO_{\delta},  \tilde\beta,r^{(7)}}, \ \|e^{-\frac{\rm i}{2}\alpha}p_{l,j}+e^{\frac{\rm i}{2}(j-l)\alpha}q_{j,l}\|_{\CO_{\delta},  \tilde\beta,r^{(7)}}  \\
&<&
\left(\frac{\varepsilon^{\frac{31}{16}}}{80}+\frac{3r^{(7)}}{5}\|e^{\frac{\rm i}2\alpha}\eta q+e^{-\frac{\rm i}2\alpha}\xi p\|_{\CO,\beta,r}\right)(r^{(7)})^{-(l+j)}.
\end{eqnarray*}
Recalling the coefficients in (\ref{hat_p_l}) -- (\ref{hat_q_j}),
and combining with (\ref{SD_general2}), we obtain
\begin{eqnarray*}
 \|\hat p_{l,j}\|_{\CO_{\delta},  \tilde\beta,r^{(7)}}, \ \|\hat q_{l,j}\|_{\CO_{\delta},  \tilde\beta,r^{(7)}}
&<& \frac{101}{200}\left(\frac{\varepsilon^{\frac{31}{16}}}{80}+\frac{3r^{(7)}}{5} \|e^{-\frac{\rm i}2\alpha}\eta \bar q+e^{\frac{\rm i}2\alpha}\xi \bar p\|_{\CO_{\delta}, \tilde\beta,r^{(7)}}\right)(r^{(7)})^{-(l+j)}\\
& & + \, 4\cdot\frac{101}{100}\delta^{-1} \cdot \frac35 \|e^{-\frac{\rm i}2\alpha}\eta \bar q+e^{\frac{\rm i}2\alpha}\xi \bar p\|_{\CO_{\delta},  \tilde\beta,r^{(7)}} (r^{(7)})^{-(l+j)}\\
&<&\left(\frac{\varepsilon^{\frac{31}{16}}}{100}+3\delta^{-1}\|e^{\frac{\rm i}2\alpha}\eta q+e^{-\frac{\rm i}2\alpha}\xi p\|_{\CO,\beta,r}\right) (r^{(7)})^{-(l+j)},
\end{eqnarray*}
which implies \re{esti_hat_pq}. By Corollary \ref{cor_pq_lj}, we have, for $l\geq 1$,
\begin{equation}\label{cr_0}
\|e^{-\frac{\rm i}{2}\alpha}(p_{l-1,0}+e^{-\frac{\rm i}{2}(l-2)\alpha} q_{0,l-1})+(\xi\eta) e^{\frac{\rm i}{2}\alpha} (q_{l+1,0}+e^{-\frac{\rm i}{2}(l+2)\alpha} p_{0,l+1}  )\|_{\CO_{\delta},  \tilde\beta,r^{(7)}}<\frac{\varepsilon^{\frac{31}{16}}}{40(r^{(7)})^{l-1}} .
\end{equation}
Moreover, we have
\begin{eqnarray}
e^{-\frac{\rm i}{2}\alpha}\cdot\frac{ -(\xi\eta)e^{-\frac{\rm i}2\alpha}C_{l,0} }{2(e^{{\rm i}(l-1)\alpha}-e^{-{\rm i}\alpha})}+(\xi\eta)e^{\frac{\rm i}{2}\alpha}\cdot \frac{e^{\frac{\rm i}2\alpha}C_{l,0}}{2(e^{{\rm i}(l+1)\alpha}-e^{{\rm i}\alpha})} &=& 0,\label{cr_1}\\
e^{-\frac{\rm i}{2}\alpha}\cdot\frac{ (\xi\eta)e^{{\rm i}(l-1)\alpha}e^{\frac{\rm i}2\alpha}\bar C_{l,0} }{2(e^{{\rm i}(l-1)\alpha}-e^{-{\rm i}\alpha})}+(\xi\eta)e^{\frac{\rm i}{2}\alpha}\cdot \frac{-e^{{\rm i}(l+1)\alpha}e^{-\frac{\rm i}2\alpha}\bar C_{l,0}}{2(e^{{\rm i}(l+1)\alpha}-e^{{\rm i}\alpha})} &=& 0,\label{cr_2}\\
e^{-\frac{\rm i}{2}\alpha}\cdot\frac{(\xi\eta)e^{\frac{\rm i}2(l-1)\alpha}C_{0,l} }{2(e^{{\rm i}(l-1)\alpha}-e^{-{\rm i}\alpha})}+(\xi\eta)e^{\frac{\rm i}{2}\alpha}\cdot \frac{-e^{\frac{\rm i}2(l+1)\alpha}C_{0,l}}{2(e^{{\rm i}(l+1)\alpha}-e^{{\rm i}\alpha})} &=& 0,\label{cr_3}\\
e^{-\frac{\rm i}{2}\alpha}\cdot\frac{ -(\xi\eta)e^{-\frac{\rm i}2(l-1)\alpha}e^{-{\rm i}\alpha}\bar C_{0,l} }{2(e^{{\rm i}(l-1)\alpha}-e^{-{\rm i}\alpha})}+(\xi\eta)e^{\frac{\rm i}{2}\alpha}\cdot \frac{e^{-\frac{\rm i}2(l+1)\alpha}e^{{\rm i}\alpha}\bar C_{0,l}}{2(e^{{\rm i}(l+1)\alpha}-e^{{\rm i}\alpha})} &=& 0.\label{cr_4}
\end{eqnarray}
In view of the definition of coefficients of $\hat p$ and $\hat q$ in (\ref{hat_p_l}) -- (\ref{hat_q_j}), we have, for $1\leq l \leq K$,
\begin{eqnarray}
 & &  e^{-\frac{\rm i}{2}\alpha}\hat p_{l-1,0}+(\xi\eta)e^{\frac{\rm i}{2}\alpha}\hat q_{l+1,0} \nonumber\\
 &=& \frac{e^{-\frac{\rm i}{2}\alpha}}{2}(p_{l-1,0}+e^{-\frac{\rm i}{2}(l-2)\alpha} q_{0,l-1} )+ \frac{(\xi\eta)e^{\frac{\rm i}{2}\alpha}}2(q_{l+1,0}+e^{-\frac{\rm i}{2}(l+2)\alpha} p_{0,l+1} )\nonumber\\
 & &  + \, (\xi\eta) e^{-\frac{\rm i}{2}\alpha} \cdot\frac{-e^{-\frac{\rm i}2\alpha}C_{l,0}+e^{{\rm i}(l-1)\alpha}e^{\frac{\rm i}2\alpha}\bar C_{l,0}+e^{\frac{\rm i}2(l-1)\alpha}C_{0,l}-e^{-\frac{\rm i}2(l-1)\alpha}e^{-{\rm i}\alpha}\bar C_{0,l}}{2(e^{{\rm i}(l-1)\alpha}-e^{-{\rm i}\alpha})}\label{C1}\\
 & &  + \,(\xi\eta)e^{\frac{\rm i}{2}\alpha} \frac{e^{\frac{\rm i}2\alpha}C_{l,0}-e^{{\rm i}(l+1)\alpha}e^{-\frac{\rm i}2\alpha}\bar C_{l,0}-e^{\frac{\rm i}2l\alpha}C_{0,l}+e^{-\frac{\rm i}2(l+1)\alpha}e^{{\rm i}\alpha}\bar C_{0,l}}{2(e^{{\rm i}(l+1)\alpha}-e^{{\rm i}\alpha})}.\label{C2}
\end{eqnarray}
Combining (\ref{cr_0}) -- (\ref{cr_4}), we have $\re{C1} + \re{C2}=0$, so that
\begin{equation*}
\|e^{-\frac{\rm i}{2}\alpha}\hat p_{l-1,0}+(\xi\eta)e^{\frac{\rm i}{2}\alpha}\hat q_{l+1,0}\|_{\CO_{\delta},  \tilde\beta,r^{(7)}} < \frac{\varepsilon^{\frac{31}{16}}}{2(r^{(7)})^l},\quad 1\leq l \leq K+1.
\end{equation*}
Similarly, we have
\begin{equation*}
\|(\xi\eta)e^{-\frac{\rm i}{2}\alpha}\hat p_{0,j+1}+ e^{\frac{\rm i}{2}\alpha}\hat q_{0,j-1}\|_{\CO_{\delta},  \tilde\beta,r^{(7)}} < \frac{\varepsilon^{\frac{31}{16}}}{2 (r^{(7)})^j}\quad 1\leq j \leq K+1.
\end{equation*}
Then, according to \re{K_beta}, we obtain \re{crossing_hat_pq} by
\begin{eqnarray*}
\|e^{-\frac{\rm i}{2}\alpha}\xi \hat p +e^{\frac{\rm i}{2}\alpha}\eta \hat q\|_{\CO_{\delta},  \tilde\beta,r^{(7)}}&\leq&\sum_{l=1}^{K+1} |e^{-\frac{\rm i}{2}\alpha}\hat p_{l-1,0}+(\xi\eta)e^{\frac{\rm i}{2}\alpha}\hat q_{l+1,0}|_{\CO_{\delta},  \tilde\beta,r^{(7)}} (r^{(7)})^l\\
& &  + \, \sum_{j=1}^{K+1} |(\xi\eta)e^{-\frac{\rm i}{2}\alpha}\hat p_{0,j+1}+ e^{\frac{\rm i}{2}\alpha}\hat q_{0,j-1}|_{\CO^{(6)}_{\delta},8\beta_+} (r^{(7)})^j\\
&<&2(K+1)\varepsilon^{\frac{31}{16}}<\frac{\varepsilon^{\frac{61}{32}}}{40}.
\end{eqnarray*}

In view of (\ref{ul}) and (\ref{uj}), and recalling that $\|f\|_{\CO_{\delta},  \tilde\beta,r^{(7)}}<\frac{\varepsilon}{2}$ in Corollary \ref{cor_fg}, we have
$$\|\hat u_{l,j}\|_{\CO_{\delta},  \tilde\beta, r^{(7)}}<\frac{101}{100}\delta^{-1}\cdot \left(1+\frac{101}{100}\right)\frac{\varepsilon}{2(r^{(7)})^{l+j}} <\frac{6\delta^{-1}\varepsilon}{5(r^{(7)})^{l+j}}.$$
Hence, recalling that $\delta>80\varepsilon^{\frac{1}{60s}}$ and \re{K_beta}, we have
$$\|\hat u\|_{\CO_{\delta},  \tilde\beta,r^{(7)}}\leq\sum_{l=2}^K\|\hat u_{l,0}\|_{\CO_{\delta},  \tilde\beta,r^{(7)}}(r^{(7)})^l+ \sum_{j=0}^K\|\hat u_{0,j}\|_{\CO_{\delta},  \tilde\beta,r^{(7)}}(r^{(7)})^j< 2K\cdot \frac65\delta^{-1}\varepsilon
< \frac{\varepsilon^{\frac{49}{50}}}{20}.$$
We have similarly, $\displaystyle\|\hat v\|_{\CO_{\delta},  \tilde\beta,r^{(7)}}<\frac{\varepsilon^{\frac{49}{50}}}{20}$.

In $\eta \hat u+\xi \hat v$, $(\eta \hat u+\xi \hat v)_{0,0}=0$ since $\hat u_{1,0}=\hat v_{0,1}=0$. For other terms, we have
\begin{eqnarray*}
& &\sum_{l=1}^{K+1}\left\|(\eta \hat u+\xi \hat v)_{l,0}\xi^l\right\|_{\CO_{\delta}, \tilde\beta, r^{(7)}}\\
&=&\sum_{l=1}^{K+1}\left\| \left(\hat u_{l+1,0}\cdot \xi\eta + \hat v_{l-1,0} \right) \xi^l \right\|_{\CO_{\delta}, \tilde\beta, r^{(7)}}\\
&=&  \frac{1}2 \sum_{l=1}^{K+1}\left\|\frac{f_{l+1,0}-e^{{\rm i}(l+2)\alpha}\bar f_{l+1,0}}{e^{{\rm i}(l+1)\alpha}-e^{{\rm i}\alpha}}\xi\eta +\frac{g_{l-1,0}-e^{{\rm i}(l-2)\alpha}\bar g_{l-1,0}}{e^{{\rm i}(l-1)\alpha}-e^{-{\rm i}\alpha}}\right\|_{\CO_{\delta}, \tilde\beta, r^{(7)}}(r^{(7)})^l  \\
&=&\frac{1}2 \sum_{l=1}^{K+1} \left\| \frac{((\xi\eta)e^{-{\rm i}\alpha} f_{l+1,0}+e^{{\rm i}\alpha} g_{l-1,0})-e^{{\rm i}l\alpha}((\xi\eta)e^{{\rm i}\alpha}\bar f_{l+1,0}+e^{-{\rm i}\alpha}\bar g_{l-1,0})}{e^{{\rm i}l\alpha}-1} \right\|_{\CO_{\delta}, \tilde\beta, r^{(7)}} (r^{(7)})^l ,
\end{eqnarray*}
and similarly,
\begin{eqnarray*}
& &\sum_{j=1}^{K+1}\left\|(\eta \hat u+\xi \hat v)_{0,j}\eta^j\right\|_{\CO_{\delta}, \tilde\beta, r^{(7)}}\\
&=& \frac{1}2 \sum_{j=1}^{K+1} \left\| \frac{(e^{-{\rm i}\alpha} f_{0,j-1}+(\xi\eta)e^{{\rm i}\alpha} g_{0,j+1})-e^{-{\rm i}j\alpha}(e^{{\rm i}\alpha}\bar f_{0,j-1}+ (\xi\eta)e^{-{\rm i}\alpha}\bar g_{0,j+1})}{e^{-{\rm i}j\alpha}-1} \right\|_{\CO_{\delta}, \tilde\beta, r^{(7)}} (r^{(7)})^j.
\end{eqnarray*}
Note that $(\xi\eta)e^{-{\rm i}\alpha} f_{l+1,0}+e^{{\rm i}\alpha} g_{l-1,0}$ and $e^{-{\rm i}\alpha} f_{0,j-1}+(\xi\eta)e^{{\rm i}\alpha} g_{0,j+1}$ are respectively the coefficients of $\xi^l$ and $\eta^j$ in $e^{-{\rm i}\alpha}\eta f+e^{{\rm i}\alpha}\xi g$ .
Hence, in view of Lemma \ref{lem_cond_F}, we have
\begin{eqnarray*}
& &\sum_{l=1}^{K+1}\left\|(\eta \hat u+\xi \hat v)_{l,0}\xi^l\right\|_{\CO_{\delta}, \tilde\beta, r^{(7)}}+ \sum_{j=1}^{K+1}\left\|(\eta \hat u+\xi \hat v)_{0,j}\eta^j\right\|_{\CO_{\delta}, \tilde\beta, r^{(7)}}\\
&<&2(K+1)\delta^{-1}\left(1+\frac{101}{100}\right)\| e^{-{\rm i}\alpha}\eta f+e^{{\rm i}\alpha}\xi g \|_{\CO,\beta, r}\\
&<&2(K+1)\delta^{-1}\left(1+\frac{101}{100}\right)\left(\frac{\varepsilon^{\frac{31}{16}}}{40}+2\|e^{\frac{\rm i}2\alpha}\eta q+e^{-\frac{\rm i}2\alpha}\xi p\|_{\CO,\beta,r}\right)\\
&<&\frac{\varepsilon^{\frac{61}{32}}}{16}+5(K+1)\delta^{-1}\|e^{\frac{\rm i}2\alpha}\eta q+e^{-\frac{\rm i}2\alpha}\xi p\|_{\CO,\beta,r}.
\end{eqnarray*}
This finishes the proof of Lemma \ref{lem_cohomo}.\qed

\smallskip

By Lemma \ref{lem_cohomo}, we see that $\phi={\rm Id}+\hat{\cL U}$ is invertible on ${\cL C}_{\omega,\tilde\beta}^{r^{(7)}}$ and $\phi({\cL C}_{\omega,\tilde\beta}^{r^{(7)}})\subset {\cL C}_{\omega,\beta}^{r}$ for $\omega\in{\CO}_\delta(r^{(7)},\tilde\beta)= {\CO}_\delta \, \cap \, ]-(r^{(7)})^2+\tilde\beta, (r^{(7)})^2-\tilde\beta [$. Indeed, according to \re{esti_uv}, \re{esti_Luv} and \re{K_beta}, we have
\begin{eqnarray*}
\left\|\eta \hat u+\xi \hat v+\hat u\hat v\right\|_{\CO_{\delta}, \tilde\beta, r^{(7)}}&<&\frac{\varepsilon^{\frac{61}{32}}}{16}+5(K+1)\delta^{-1}\|e^{\frac{\rm i}2\alpha}\eta q+e^{-\frac{\rm i}2\alpha}\xi p\|_{\CO,\beta,r}+\frac{\varepsilon^{\frac{49}{25}}}{400}\\
&<& \frac{\varepsilon^{\frac{61}{32}}}{16}+5(K+1)\delta^{-1}\cdot\frac{101}{200}\cdot\frac{\varepsilon}{10}+\frac{\varepsilon^{\frac{49}{25}}}{400} \, < \, \varepsilon^{\frac{49}{50}}.
\end{eqnarray*}
Then, in view of the definition of the set given in (\ref{cwb}), we have for any $(\xi,\eta)\in{\cL C}_{\omega,\tilde\beta}^{r^{(7)}}$ with $\omega\in {\CO}_\delta(r^{(7)},\tilde\beta)$,
\begin{eqnarray*}
|(\xi+\hat u(\xi,\eta))(\eta+\hat v(\xi,\eta)) -\omega|&\leq&|\xi\eta-\omega|+|\eta \hat u(\xi,\eta)+\xi \hat v(\xi,\eta)+\hat u(\xi,\eta)\hat v(\xi,\eta)|\\
&<&\tilde\beta+\|\eta \hat u+\xi \hat v+\hat u\hat v\|_{\CO_{\delta}, \tilde\beta, r^{(7)}}<\tilde\beta+\varepsilon^{\frac{49}{50}}<\beta,\\
|\xi+\hat u(\xi,\eta)|, \ |\eta+\hat v(\xi,\eta)|
&<&r^{(7)}+\varepsilon^{\frac{49}{50}}<r,
\end{eqnarray*}
noting that, under the condition (\ref{varepsilon_small-general}),
 $$\varepsilon^{\frac{49}{50}}<\varepsilon^{\frac{1}{40s}}\leq \frac\beta2<\beta-\tilde\beta,\quad \varepsilon^{\frac{49}{50}}<\varepsilon^{\frac{1}{2400s^2}}<\frac{r-r_+}8=r-r^{(7)}.$$
Moreover, since $\hat u,\hat v\in{\CA}^\R_{\tilde\beta, r^{(7)}}(\CO_\delta)$, we have, by Lemma \ref{lem_reve_psi}, $\rho\circ\phi=\phi\circ\rho$.

\smallskip

With $\phi={\rm Id}+\hat{\cL U}$ constructed above, we define $\tilde\tau_1:=\phi^{-1}\circ\tau_1\circ\phi$.
For any $h=h(\xi\eta)$, we define the linear operator $L_h$ by:
\begin{equation}\label{Lh}
\left[L_{h}\left(\begin{array}{c}
                  p_1 \\
                  p_2
                \end{array}
\right)\right](\xi,\eta):=e^{-{\rm i}h(\xi\eta)}\xi p_1(\xi,\eta)+ e^{{\rm i}h(\xi\eta)} \eta p_2(\xi,\eta).
\end{equation}

\begin{lemma}\label{prop_pq+}
For $\tilde\tau_{1}=\phi^{-1}\circ\tau_1\circ\phi$, we have
$$\tilde\tau_{1}(\xi,\eta)-\left(\begin{array}{l}
(e^{\frac{\rm i}2\alpha(\xi\eta)}+p_{0,1})\eta\\
(e^{\frac{\rm i}2\alpha(\xi\eta)}+p_{0,1})^{-1}\xi
\end{array}\right)\in ({\CA}_{2\beta_+,\tilde r}(\CO_{\delta}))^2,
$$
satisfying that
\beq
\left\|\tilde\tau_{1}(\xi,\eta)-\left(\begin{array}{l}
(e^{\frac{\rm i}2\alpha(\xi\eta)}+p_{0,1})\eta\\
(e^{\frac{\rm i}2\alpha(\xi\eta)}+p_{0,1})^{-1}\xi
\end{array}\right)
\right\|_{\CO_\delta,2\beta_+,\tilde r}<\frac{\varepsilon^{\frac{61}{32}}}{3}+\frac{22(K+1)}{\delta}\|e^{\frac{\rm i}2\alpha}\eta q+e^{-\frac{\rm i}2\alpha}\xi p\|_{\CO,\beta, r},\label{esti_pq+}
\eeq
\begin{equation}
\left\|L_{\frac{\alpha}2}\left(\tilde\tau_{1}(\xi,\eta)-\left(\begin{array}{l}
(e^{\frac{\rm i}2\alpha(\xi\eta)}+p_{0,1})\eta\\
(e^{\frac{\rm i}2\alpha(\xi\eta)}+p_{0,1})^{-1}\xi
\end{array}\right) \right)
\right\|_{\CO_\delta,2\beta_+,\tilde r}  < \frac{\varepsilon^{\frac{61}{32}}}{2}.\label{esti_Lpq+}
\end{equation}
\end{lemma}
\proof A direct computation yields that
\begin{eqnarray*}
\tilde\tau_1&=&\phi^{-1}\circ\tau_1\circ\phi\nonumber\\
&=&\left(\begin{array}{l}
e^{\frac{\rm i}2\alpha(\xi\eta+\eta{\hat u} +\xi \hat v+\hat u\hat v)}\eta+e^{\frac{\rm i}2\alpha(\xi\eta+\eta{\hat u} +\xi \hat v+\hat u\hat v)} \hat v+p\circ\phi
-\hat u\circ\tau_1\circ\phi\\[2mm]
e^{-\frac{\rm i}2\alpha(\xi\eta+\eta{\hat u} +\xi \hat v+\hat u\hat v)}\xi+e^{-\frac{\rm i}2\alpha(\xi\eta+\eta{\hat u} +\xi \hat v+\hat u\hat v)} \hat u+q\circ\phi
-\hat v\circ\tau_1\circ\phi
\end{array}\right) +  \tilde{\cL U}\circ\tau_1\circ\phi,
\end{eqnarray*}
with $\tilde{\cL U}=\phi^{-1}-{\rm Id}+\hat{\cL U}$ defined as in the proof of Lemma \ref{lem_tildeU}.
Recalling \re{pK}, we see that
\begin{eqnarray}
& & \tilde\tau_{1}(\xi,\eta)-\left(\begin{array}{c}
(e^{\frac{\rm i}2\alpha(\xi\eta)}+p_{0,1})\eta\\[2mm]
(e^{\frac{\rm i}2\alpha(\xi\eta)}+p_{0,1})^{-1}\xi
\end{array}\right)\nonumber\\
&=&\left(\begin{array}{c}
0\\[2mm]
\left((e^{-\frac{\rm i}2\alpha(\xi\eta)}+q_{1,0})-(e^{\frac{\rm i}2\alpha(\xi\eta)}+p_{0,1})^{-1}\right)\xi\end{array}\right)\label{pq+0}\\
& & + \, \left(\begin{array}{l}
e^{\frac{\rm i}2\alpha(\xi\eta)} \hat v-\hat u(e^{\frac{\rm i}{2}\alpha}\eta,e^{-\frac{\rm i}{2}\alpha}\xi)+p_{K}-p_{0,1}\eta\\[2mm]
e^{-\frac{\rm i}2\alpha(\xi\eta)} \hat u-\hat v(e^{\frac{\rm i}{2}\alpha}\eta,e^{-\frac{\rm i}{2}\alpha}\xi)+q_{K}-q_{1,0}\xi
\end{array}\right) \label{pq+1} \\
& & + \, \left(\begin{array}{l}
p-p_{K}\\[1mm]
q-q_{K}
\end{array}\right)+ \left(\begin{array}{l}
p\circ\phi-p\\[1mm]
q\circ\phi-q
\end{array}\right)\label{pq+2}\\
& & + \, \left(\begin{array}{l}
(e^{\frac{\rm i}2\alpha(\xi\eta+\eta \hat u +\xi \hat v+\hat u\hat v)}-e^{\frac{\rm i}2\alpha(\xi\eta)})(\eta+\hat v)\\[2mm]
(e^{-\frac{\rm i}2\alpha(\xi\eta+\eta \hat u +\xi \hat v+\hat u\hat v)}-e^{-\frac{\rm i}2\alpha(\xi\eta)})(\xi+\hat u)
\end{array}\right) \label{pq+3} \\
& & - \,  \left(\hat{\cL U}\circ\tau_1\circ\phi-\hat{\cL U}(e^{\frac{\rm i}{2}\alpha}\eta,e^{-\frac{\rm i}{2}\alpha}\xi)\right)+ \  \tilde{\cL U}\circ\tau_1\circ\phi.  \label{pq+4}
\end{eqnarray}

\smallskip
In what follows, we shall estimate the norms of terms \re{pq+0} -- \re{pq+4} as well as their image under $L_{\frac{\alpha}2}$. We emphasize that if a given term $T\in ({\CA}_{2\beta_+, \tilde r}(\CO_\delta))^2$ satisfies $\|T\|_{\CO_\delta, 2\beta_+, \tilde r}< D \varepsilon^{\varsigma}$ with $\varsigma>1$ and some constant $D>0$, then $\|L_{\frac{\alpha}2}(T)\|_{\CO_\delta, 2\beta_+, \tilde r}< \frac{101}{200} D \varepsilon^{\varsigma}< D \varepsilon^{\varsigma}$.

\smallskip

\begin{itemize}
\item {\bf Terms in (\ref{pq+0})}
\end{itemize}

In view of (\ref{coeff_res_pq}) in Corollary \ref{cor_coeff_res_pq_fg}, we obtain, for $\omega\in \CO_\delta(\tilde r,2\beta_+)$,
\begin{equation}
|(e^{\frac{\rm i}2\alpha}+p_{0,1})(e^{-\frac{\rm i}2\alpha}+q_{1,0})-1|_{\omega,2\beta_+}=|e^{\frac{\rm i}2\alpha}q_{1,0}+e^{-\frac{\rm i}2\alpha}p_{0,1}+p_{0,1}q_{1,0}|_{\omega,2\beta_+}<\frac{\varepsilon^{\frac{61}{32}}}{50r^{(7)}}.\label{ineq_LambdaLambda}
\end{equation}
Since Corollary \ref{lem_eiba} implies that, for any $(\xi,\eta)\in {\cL C}_{\omega,2\beta_+}$, $\omega\in \CO_\delta(\tilde r,2\beta_+)$,
$$|e^{\frac{\rm i}2\alpha(\xi\eta)}+p_{0,1}(\xi\eta)|\geq  |e^{\frac{\rm i}2\alpha(\xi\eta)}|-|p_{0,1}|_{\omega,2\beta_+}>\frac45,$$
we see that
$((e^{-\frac{\rm i}2\alpha}+q_{1,0})-(e^{\frac{\rm i}2\alpha}+p_{0,1})^{-1})(\omega)$ is analytic on $ \CO_\delta(\tilde r,2\beta_+)$,
and (\ref{ineq_LambdaLambda}) implies that
$$
|(e^{-\frac{\rm i}2\alpha}+q_{1,0})-(e^{\frac{\rm i}2\alpha}+p_{0,1})^{-1}|_{\omega,2\beta_+}<\frac{\varepsilon^{\frac{61}{32}}}{40r^{(7)}}.
$$
Hence, $(\ref{pq+0})\in(\CA_{2\beta_+,\tilde r}(\CO_{\delta}))^2$,
\begin{equation}\label{pq+0_}
\left\|\left(\begin{array}{c}
0\\[2mm]
\left((e^{-\frac{\rm i}2\alpha(\xi\eta)}+q_{1,0})-(e^{\frac{\rm i}2\alpha(\xi\eta)}+p_{0,1})^{-1}\right)\xi\end{array}\right)\right\|_{\CO_{\delta},2\beta_+,\tilde r}<\frac{\varepsilon^{\frac{61}{32}}}{40}.
\end{equation}

\smallskip

\begin{itemize}
  \item {\bf Terms in (\ref{pq+1})}
\end{itemize}

According to Lemma \ref{lem_cohomo}, $\hat u, \, \hat v\in {\CA}^\R_{\tilde\beta,r^{(7)}}(\CO_\delta)$, which implies that
$$\left(\begin{array}{l}
e^{\frac{\rm i}2\alpha(\xi\eta)} \hat v-\hat u(e^{\frac{\rm i}{2}\alpha}\eta,e^{-\frac{\rm i}{2}\alpha}\xi)+p_{K}-p_{0,1}(\xi\eta)\eta\\[2mm]
e^{-\frac{\rm i}2\alpha(\xi\eta)} \hat u-\hat v(e^{\frac{\rm i}{2}\alpha}\eta,e^{-\frac{\rm i}{2}\alpha}\xi)+q_{K}-q_{1,0}(\xi\eta) \xi
\end{array}\right)\in ({\CA}_{2\beta_+,\tilde r}(\CO_\delta))^2.$$
By (\ref{cohomo1}) -- (\ref{crossing_cohomo}), we obtain
\begin{multline}
\left\|\left(\begin{array}{l}
e^{\frac{\rm i}2\alpha(\xi\eta)} \hat v-\hat u(e^{\frac{\rm i}{2}\alpha}\eta,e^{-\frac{\rm i}{2}\alpha}\xi)+p_{K}-p_{0,1}(\xi\eta)\eta\\[2mm]
e^{-\frac{\rm i}2\alpha(\xi\eta)} \hat u-\hat v(e^{\frac{\rm i}{2}\alpha}\eta,e^{-\frac{\rm i}{2}\alpha}\xi)+q_{K}-q_{1,0}(\xi\eta) \xi
\end{array}\right)
 \right\|_{\CO_\delta,2\beta_+ , \tilde r}\\
 <\frac{\varepsilon^{\frac{61}{32}}}{40}+6(K+1)\delta^{-1}\|e^{\frac{\rm i}2\alpha}\eta q+e^{-\frac{\rm i}2\alpha}\xi p\|_{\CO,\beta, r},\label{pq+1_}
\end{multline}
\begin{equation}\label{Lpq+1_}
\left\|
 L_{\frac\alpha2}\left(\begin{array}{l}
e^{\frac{\rm i}2\alpha(\xi\eta)} \hat v-\hat u(e^{\frac{\rm i}{2}\alpha}\eta,e^{-\frac{\rm i}{2}\alpha}\xi)+p_{K}-p_{0,1}(\xi\eta)\eta\\[2mm]
e^{-\frac{\rm i}2\alpha(\xi\eta)} \hat u-\hat v(e^{\frac{\rm i}{2}\alpha}\eta,e^{-\frac{\rm i}{2}\alpha}\xi)+q_{K}-q_{1,0}(\xi\eta) \xi
\end{array}\right)
 \right\|_{\CO_\delta,2\beta_+ ,\tilde r} <\frac{\varepsilon^{\frac{61}{32}}}{20}.
\end{equation}

\smallskip

\begin{itemize}
\item {\bf Terms in (\ref{pq+2})}
\end{itemize}

In view of (\ref{rest_pq}), we have $p-p_{K}$, $q-q_{K}\in {\CA}_{2\beta_+,\tilde r}(\CO_\delta)$, with
\begin{equation}\label{pq+2_}
\left\|  \left(\begin{array}{l}
p-p_{K}\\[1mm]
q-q_{K}
\end{array}\right)\right\|_{\CO_{\delta},2\beta_+,\tilde r}<\frac{\varepsilon^2}{5}.
\end{equation}
Since $\tilde r=\frac{r+r_+}{2}$ and $\beta\geq \varepsilon^{\frac{1}{40s}}$, we have, by Lemma \ref{lem_norm}, $$\left(\begin{array}{l}
p\circ\phi-p\\[1mm]
q\circ\phi-q
\end{array}\right)\in({\CA}_{2\beta_+,\tilde r}(\CO_\delta))^2,$$ satisfying that
\begin{equation}\label{pq+2__}
\left\|  \left(\begin{array}{l}
p\circ\phi-p\\[1mm]
q\circ\phi-q
\end{array}\right)\right\|_{\CO_\delta,2\beta_+ ,\tilde r} <\frac{6r}{(r-\tilde r)\beta} \cdot \frac{\varepsilon^{\frac{49}{50}}}{20}\cdot\frac{\varepsilon}{10}
\leq\frac{3r\varepsilon^{\frac{99}{50}-\frac{1}{40s}}}{50(r-r_+)}<\frac{\varepsilon^{\frac{61}{32}}}{80}.
\end{equation}
The last inequality follows from \re{1/(r1-r2)}.

\smallskip

\begin{itemize}
  \item {\bf Terms in (\ref{pq+3})}
\end{itemize}

In view of Lemma \ref{lem_cohomo},
we see that $\eta \hat u +\xi \hat v+\hat u\hat v\in \CA_{2\beta_+,\tilde r}(\CO_\delta)$ and
\begin{eqnarray}
\|\eta \hat u +\xi \hat v+\hat u\hat v\|_{\CO_\delta, 2\beta_+,\tilde r}
&<& \frac{\varepsilon^{\frac{61}{32}}}{16}+5(K+1)\delta^{-1}\|e^{\frac{\rm i}2\alpha}\eta q+e^{-\frac{\rm i}2\alpha}\xi p\|_{\CO,\beta, r}+\frac{\varepsilon^{\frac{49}{25}}}{400}\nonumber\\
&<&\frac{\varepsilon^{\frac{61}{32}}}{12}+5(K+1)\delta^{-1}\|e^{\frac{\rm i}2\alpha}\eta q+e^{-\frac{\rm i}2\alpha}\xi p\|_{\CO,\beta, r}.\label{eta_hat_u_uv}
\end{eqnarray}
By Lemma \ref{lem_a1}, we obtain $\alpha(\xi\eta+\eta \hat u +\xi \hat v+\hat u\hat v)-\alpha(\xi\eta)\in  \CA_{2\beta_+,\tilde r}(\CO_\delta)$, and for $-1\leq  b \leq 1$,
$e^{{\rm i}b\alpha(\xi\eta+\eta \hat u +\xi \hat v+\hat u\hat v)}-e^{{\rm i}b\alpha(\xi\eta)}\in \CA_{2\beta_+,\tilde r}(\CO_\delta)$,
with
\begin{equation*}
\left\|\alpha(\xi\eta+\eta \hat u +\xi \hat v+\hat u\hat v)-\alpha(\xi\eta)\right\|_{\CO_\delta,2\beta_+,\tilde r}<\frac{\varepsilon^{\frac{61}{32}}}{8}+7(K+1)\delta^{-1}\|e^{\frac{\rm i}2\alpha}\eta q+e^{-\frac{\rm i}2\alpha}\xi p\|_{\CO,\beta, r},
\end{equation*}
\begin{equation*}
\left\|e^{{\rm i}b\alpha(\xi\eta+\eta \hat u +\xi \hat v+\hat u\hat v)}-e^{{\rm i}b\alpha(\xi\eta)}\right\|_{\CO_\delta,2\beta_+,\tilde r}
<   \frac{\varepsilon^{\frac{61}{32}}}{9}+7(K+1)\delta^{-1}\|e^{\frac{\rm i}2\alpha}\eta q+e^{-\frac{\rm i}2\alpha}\xi p\|_{\CO,\beta, r},\;\  -1\leq  b \leq 1.
\end{equation*}
Hence, we have
$$\left(\begin{array}{l}(e^{\frac{\rm i}2\alpha(\xi\eta+\eta \hat u +\xi \hat v+\hat u\hat v)}-e^{\frac{\rm i}2\alpha(\xi\eta)})(\eta+\hat v)\\[2mm]
(e^{-\frac{\rm i}2\alpha(\xi\eta+\eta \hat u +\xi \hat v+\hat u\hat v)}-e^{-\frac{\rm i}2\alpha(\xi\eta)})(\xi+\hat u)
\end{array}\right) \in (\CA_{2\beta_+,\tilde r}(\CO_\delta))^2,$$
satisfying that
\begin{eqnarray}
& & \left\|\left(\begin{array}{l}(e^{\frac{\rm i}2\alpha(\xi\eta+\eta \hat u +\xi \hat v+\hat u\hat v)}-e^{\frac{\rm i}2\alpha(\xi\eta)})(\eta+\hat v)\\[2mm]
(e^{-\frac{\rm i}2\alpha(\xi\eta+\eta \hat u +\xi \hat v+\hat u\hat v)}-e^{-\frac{\rm i}2\alpha(\xi\eta)})(\xi+\hat u)
\end{array}\right)\right\|_{\CO_\delta,2\beta_+,\tilde r}\nonumber  \\
&<&\left(\frac{\varepsilon^{\frac{61}{32}}}{9}+7(K+1)\delta^{-1}\|e^{\frac{\rm i}2\alpha}\eta q+e^{-\frac{\rm i}2\alpha}\xi p\|_{\CO,\beta, r}\right)\cdot\left(\frac14+\frac{\varepsilon^{\frac{49}{50}}}{20}\right) \nonumber \\
&<&\frac{\varepsilon^{\frac{61}{32}}}{30}+2(K+1)\delta^{-1}\|e^{\frac{\rm i}2\alpha}\eta q+e^{-\frac{\rm i}2\alpha}\xi p\|_{\CO,\beta, r}\label{pq+3_}
\end{eqnarray}

In order to obtain the estimate of the image under $L_{\frac{\alpha}2}$ of this term, we shall follow the scheme of the proof of \rl{lem_appro}. Developing $e^{\frac{\rm i}2\alpha(\cdot)}$ around $\xi\eta$,
\begin{align}
 e^{\frac{\rm i}2\alpha(\xi\eta+\eta \hat u +\xi \hat v+\hat u\hat v)}-e^{\frac{\rm i}2\alpha(\xi\eta)}
&= e^{\frac{\rm i}2\alpha(\xi\eta)}\sum_{k\geq 1}\frac{{\rm i}^k}{2^k\cdot k!}(\alpha(\xi\eta+\eta \hat u +\xi \hat v+\hat u\hat v)-\alpha(\xi\eta))^k\label{ei2alpha}\\
&= \frac{{\rm i}}2 e^{\frac{\rm i}2\alpha(\xi\eta)}\alpha'(\xi\eta)(\eta \hat u +\xi \hat v) \nonumber\\
& + \, \frac{\rm i}2 e^{\frac{\rm i}2\alpha(\xi\eta)}\alpha'(\xi\eta) \hat u\hat v+e^{\frac{\rm i}2\alpha(\xi\eta)}\sum_{j\geq 2}  \frac{\alpha^{(j)}(\xi\eta)}{j!}(\eta \hat u +\xi \hat v+\hat u\hat v)^j\nonumber\\
& + \, e^{\frac{\rm i}2\alpha(\xi\eta)}\sum_{k\geq 2}\frac{{\rm i}^k}{2^k\cdot k!}(\alpha(\xi\eta+\eta \hat u +\xi \hat v+\hat u\hat v)-\alpha(\xi\eta))^k.\nonumber
\end{align}
Noting that $\displaystyle \|\hat u\hat v\|_{\CO_\delta,\tilde\beta,r^{(7)}}<\frac{\varepsilon^{\frac{49}{25}}}{400}$, and (\ref{eta_hat_u_uv}) gives the rough estimates via (\ref{esti_uv}):
\begin{eqnarray}
\|\eta \hat u +\xi \hat v+\hat u\hat v\|_{\CO_\delta,2\beta_+,\tilde r} & <   &\frac{\varepsilon^{\frac{49}{50}}}{20}, \label{eta_hat_u_uv_rough}\\
\left\|\alpha(\xi\eta+\eta \hat u +\xi \hat v+\hat u\hat v)-\alpha(\xi\eta)\right\|_{\CO_\delta,2\beta_+,\tilde r} & <   & \frac{\varepsilon^{\frac{49}{50}}}{16},\label{a_eta_hat_u_uv_rough}  \\
\left\|e^{{\rm i}b\alpha(\xi\eta+\eta \hat u +\xi \hat v+\hat u\hat v)}-e^{{\rm i}b\alpha(\xi\eta)}\right\|_{\CO_\delta,2\beta_+,\tilde r} & <   &\frac{\varepsilon^{\frac{49}{50}}}{15},\quad -1\leq  b \leq 1,\label{eia_eta_hat_u_uv_rough}
\end{eqnarray}
we have, in view of \re{ei2alpha},
\begin{equation}\label{er_1}
\left\|e^{\frac{\rm i}2\alpha(\xi\eta+\eta \hat u +\xi \hat v+\hat u\hat v)}-e^{\frac{\rm i}2\alpha(\xi\eta)}-\frac{{\rm i}}2 e^{\frac{\rm i}2\alpha(\xi\eta)}\alpha'(\xi\eta)(\eta \hat u +\xi \hat v)\right\|_{\CO_\delta,2\beta_+,\tilde r}<\frac{\varepsilon^{\frac{61}{32}}}3.
\end{equation}
Indeed,
according to Corollary \ref{cor_deri_a}, we have, for $\omega\in \CO_\delta(\tilde r,2\beta_+)$,
$$|\alpha'(\xi\eta)|_{\omega,2\beta_+}<\frac{6}{5},\qquad |\alpha^{(k)}(\xi\eta)|_{\omega,2\beta_+}<\left\{\begin{array}{ll}
2\beta^{-\frac1{32}}, & 2\leq k \leq s, \ {\rm if}  \  s\geq 2\\[2mm]
\frac{ k! 2^{k+5}}{\beta^{k}}, &k\geq s+1
\end{array} \right. .$$
As a consequence, \re{ei2alpha} implies that, if $s\geq2$, then,
\begin{eqnarray*}
& &\left\| e^{\frac{\rm i}2\alpha(\xi\eta+\eta \hat u +\xi \hat v+\hat u\hat v)}-e^{\frac{\rm i}2\alpha(\xi\eta)}-\frac{{\rm i}}2 e^{\frac{\rm i}2\alpha(\xi\eta)}\alpha'(\xi\eta)(\eta \hat u +\xi \hat v)\right\|_{\CO_\delta,2\beta_+,\tilde r} \\
&\leq&\frac{101}{200}\cdot\frac{6}{5}\cdot \frac{\varepsilon^{\frac{49}{25}}}{400} +\frac{101}{100}\cdot2\beta^{-\frac1{32}}\sum_{k=2}^s  \frac{1}{k!}\left(\frac{\varepsilon^{\frac{49}{50}}}{20}\right)^k\\
& & + \, \frac{101}{100}\sum_{k\geq s+1}\frac{2^{k+5}}{\beta^k}\left(\frac{\varepsilon^{\frac{49}{50}}}{20}\right)^k
+ \frac{101}{100}\sum_{k\geq 2}\frac{1}{2^k\cdot k!}\left(\frac{\varepsilon^{\frac{49}{50}}}{16}\right)^k<\frac{\varepsilon^{\frac{61}{32}}}{3}.
\end{eqnarray*}
Otherwise, for $s=1$, it is bounded by
$$\frac{101}{200}\cdot\frac{6}{5}\cdot \frac{\varepsilon^{\frac{49}{25}}}{400}+\frac{101}{100}\sum_{k\geq 2}\frac{2^{k+5}}{\beta^k}\left(\frac{\varepsilon^{\frac{49}{50}}}{20}\right)^k
+ \frac{101}{100}\sum_{k\geq 2}\frac{1}{2^k\cdot k!}\left(\frac{\varepsilon^{\frac{49}{50}}}{16}\right)^k<\frac{\varepsilon^{\frac{61}{32}}}{3}. $$
Similarly,
\begin{equation}\label{er_2}
\left\|e^{-\frac{\rm i}2\alpha(\xi\eta+\eta \hat u +\xi \hat v+\hat u\hat v)}-e^{-\frac{\rm i}2\alpha(\xi\eta)}+\frac{{\rm i}}2 e^{\frac{\rm i}2\alpha(\xi\eta)}\alpha'(\xi\eta)(\eta \hat u +\xi \hat v)\right\|_{\CO_\delta, 2\beta_+, \tilde r}<\frac{\varepsilon^{\frac{61}{32}}}{3}.
\end{equation}
Noting that $ L_{\frac\alpha2}\left(\begin{array}{l}
-\frac{\rm i}2 e^{\frac{\rm i}2\alpha(\xi\eta)}\alpha'(\xi\eta) (\eta \hat u +\xi \hat v)\eta\\[2mm]
\frac{\rm i}2 e^{-\frac{\rm i}2\alpha(\xi\eta)}\alpha'(\xi\eta)  (\eta \hat u +\xi \hat v)\xi
\end{array}\right)=0$, we have
\begin{eqnarray*}
& & L_{\frac\alpha2}\left(\begin{array}{l}
\left(e^{\frac{\rm i}2\alpha(\xi\eta+\eta \hat u +\xi \hat v+\hat u\hat v)}-e^{\frac{\rm i}2\alpha(\xi\eta)}\right)(\eta+\hat v)\\[2mm]
\left(e^{-\frac{\rm i}2\alpha(\xi\eta+\eta \hat u +\xi \hat v+\hat u\hat v)}-e^{-\frac{\rm i}2\alpha(\xi\eta)}\right)(\xi+\hat u)
\end{array}\right)\nonumber \\
&=& L_{\frac\alpha2}\left(\begin{array}{l}
\left(e^{\frac{\rm i}2\alpha(\xi\eta+\eta \hat u +\xi \hat v+\hat u\hat v)}-e^{\frac{\rm i}2\alpha(\xi\eta)}\right)(\eta+\hat v)-\frac{\rm i}2 e^{\frac{\rm i}2\alpha(\xi\eta)}\alpha'(\xi\eta) (\eta \hat u +\xi \hat v)\eta\\[2mm]
\left(e^{-\frac{\rm i}2\alpha(\xi\eta+\eta \hat u +\xi \hat v+\hat u\hat v)}-e^{-\frac{\rm i}2\alpha(\xi\eta)}\right)(\xi+\hat u)+\frac{\rm i}2 e^{-\frac{\rm i}2\alpha(\xi\eta)}\alpha'(\xi\eta)  (\eta \hat u +\xi \hat v)\xi
\end{array}\right).
\end{eqnarray*}
Hence, by (\ref{er_1}), (\ref{er_2}) and Lemma \ref{lem_a1}, we have
\begin{equation}
 \left \|L_{\frac\alpha2}\left(\begin{array}{l}\left(e^{\frac{\rm i}2\alpha(\xi\eta+\eta \hat u +\xi \hat v+\hat u\hat v)}-e^{\frac{\rm i}2\alpha(\xi\eta)}\right)(\eta+\hat v)\\[2mm]
		\left(e^{-\frac{\rm i}2\alpha(\xi\eta+\eta \hat u +\xi \hat v+\hat u\hat v)}-e^{-\frac{\rm i}2\alpha(\xi\eta)}\right)(\xi+\hat u)\end{array}\right) \right\|_{\CO_\delta,2\beta_+,\tilde r}\\
	<\frac{2\tilde r}3\varepsilon^{\frac{61}{32}}+\frac{\varepsilon^{\frac{49}{25}}}{150}<\frac{\varepsilon^{\frac{61}{32}}}{5}.\label{Lpq+3_}
\end{equation}
\smallskip

\begin{itemize}
  \item {\bf Terms in (\ref{pq+4})}
\end{itemize}

Since $\displaystyle \|\hat{\cL U}\|_{\CO_\delta,8\beta_+,r^{(7)}}<\frac{\varepsilon^{\frac{49}{50}}}{10}$, and according to Lemma \ref{lem_tildeU}, $\tilde{\cL U}=\Phi^{-1}-{\rm Id}+{\cL U}$ satisfies
$$\|\tilde{\cL U}\|_{\CO_\delta,4\beta_+, r^{(6)}}\leq \frac{r^{(7)}\|\hat{\cL U}\|^2_{\CO_\delta,8\beta_+,r^{(7)}}}{(r^{(7)}-r^{(6)})\beta_+}\leq \frac{\varepsilon^{\frac{61}{32}}}{8}.$$
Then, by Lemma \ref{lem_eia}, we have that
$$\|\hat{\cL U}(e^{\frac{\rm i}2\alpha(\xi\eta)}\eta,e^{-\frac{\rm i}2\alpha(\xi\eta)}\xi)\|_{\CO_\delta,4\beta_+,r^{(6)}}<\frac{\varepsilon^{\frac{49}{50}}}{10}, \quad
\|\tilde{\cL U}(e^{\frac{\rm i}2\alpha(\xi\eta)}\eta,e^{-\frac{\rm i}2\alpha(\xi\eta)}\xi)\|_{\CO_\delta,2\beta_+,r^{(5)}}<\frac{\varepsilon^{\frac{61}{32}}}{8}.$$
By Lemma \ref{lem_norm}, we have
$$
 \|p(\xi+\hat u,\eta+\hat v)\|_{\CO_\delta,8\beta_+,r^{(7)}}, \ \|q(\xi+\hat u,\eta+\hat v)\|_{\CO_\delta,8\beta_+,r^{(7)}}
< \frac{\varepsilon}{10}+\frac{3r}{(r-r^{(7)})\beta}\cdot\frac{\varepsilon^{1+\frac{49}{50}}}{200}<\frac{\varepsilon}{8},
$$
which, together with (\ref{eia_eta_hat_u_uv_rough}), implies that
$$ \left\|\left(\begin{array}{c}
(e^{\frac{\rm i}2\alpha(\xi\eta+\eta \hat u +\xi \hat v + \hat u\hat v)}- e^{\frac{\rm i}2\alpha(\xi\eta)})\eta+p(\xi+\hat u,\eta+\hat v)\\
(e^{-\frac{\rm i}2\alpha(\xi\eta+\eta \hat u +\xi \hat v + \hat u\hat v)}- e^{-\frac{\rm i}2\alpha(\xi\eta)})\xi+q(\xi+\hat u,\eta+\hat v)
\end{array}\right)\right\|_{\CO_\delta,2\beta_+, \tilde r}
 <\frac{2r^{(7)}\varepsilon^{\frac{49}{50}}}{15} + \frac{\varepsilon}{4}
<\frac{\varepsilon^{\frac{49}{50}}}{20}.$$
Therefore, by Lemma \ref{lem_norm}, we obtain
\begin{equation}
	\| \hat{\cL U}\circ\tau_{1}\circ\phi- \hat{\cL U}(e^{\frac{\rm i}2\alpha(\xi\eta)}\eta,e^{-\frac{\rm i}2\alpha(\xi\eta)}\xi)\|_{\CO_\delta,2\beta_+, \tilde r}
	< \frac{6r^{(7)}}{8(r^{(7)}-\tilde r)\beta_+}\cdot \frac{\varepsilon^{\frac{49}{50}}}{10}\cdot \frac{\varepsilon^{\frac{49}{50}}}{20}<\frac{\varepsilon^{\frac{61}{32}}}{200}.\label{pq+4_}
\end{equation}
In a similar way, we have,
\begin{eqnarray}
 \|\tilde{\cL U}\circ\tau^{(1)}\circ\phi\|_{\CO_\delta,2\beta_+,\tilde r}
&\leq&\|\tilde{\cL U}(e^{\frac{\rm i}2\alpha(\xi\eta)}\eta,e^{-\frac{\rm i}2\alpha(\xi\eta)}\xi)\|_{\omega,4\beta_+,r^{(5)}}\nonumber\\
 &&+\| \tilde{\cL U}\circ\tau^{(1)}\circ\phi- \tilde{\cL U}(e^{\frac{\rm i}2\alpha(\xi\eta)}\eta,e^{-\frac{\rm i}2\alpha(\xi\eta)}\xi)\|_{\omega,2\beta_+,\tilde r}\nonumber\\
&<&\frac{\varepsilon^{\frac{61}{32}}}{8} +\frac{6 r^{(6)}}{4(r^{(6)}-\tilde r)\beta_+}\cdot \frac{\varepsilon^{\frac{49}{50}}}{20}\cdot \frac{\varepsilon^{\frac{61}{32}}}{8}<\frac{\varepsilon^{\frac{61}{32}}}{6}.\label{pq+4__}
\end{eqnarray}


\smallskip

With (\ref{pq+0}) -- (\ref{pq+4}) estimated as above, we have
\re{esti_pq+} by combining \re{pq+0_}, \re{pq+1_}, \re{pq+2_}, \re{pq+2__}, \re{pq+3_}, \re{pq+4_}, \re{pq+4__}, and get \re{esti_Lpq+} by combining \re{pq+0_}, \re{Lpq+1_}, \re{pq+2_}, \re{pq+2__}, \re{Lpq+3_}, \re{pq+4_}, \re{pq+4__}.
Hence Lemma \ref{prop_pq+} is shown.
\qed

\smallskip

\noindent
{\it Proof of Theorem \ref{thm_trans_cohomo}.}
By Lemma \ref{prop_pq+}, we see that $\tilde\tau_{1}=\phi^{-1}\circ\tau_1\circ\phi$ can be written as
$$\tilde\tau_{1}=\left(\begin{array}{l}
(e^{\frac{\rm i}2\alpha(\xi\eta)}+p_{0,1})\eta\\
(e^{\frac{\rm i}2\alpha(\xi\eta)}+p_{0,1})^{-1}\xi
\end{array}\right)+\left(\begin{array}{l}
\tilde p(\xi, \eta)\\
\tilde q(\xi, \eta)
\end{array}\right)$$
with $\tilde p$, $\tilde q\in {\CA}_{2\beta_+,\tilde r}(\CO_\delta)$ satisfying that
\begin{eqnarray}
\left\|\left(\begin{array}{l}
\tilde p\\
\tilde q
\end{array}\right)
\right\|_{\CO_\delta,2\beta_+,\tilde r}  & <&\frac{\varepsilon^{\frac{61}{32}}}{3}+22(K+1)\delta^{-1}\|e^{\frac{\rm i}2\alpha}\eta q+e^{-\frac{\rm i}2\alpha}\xi p\|_{\CO,\beta, r},\label{tilde_p_tilde_q}\\
\left\|L_{\frac\alpha2}\left(\begin{array}{l}
\tilde p\\
\tilde q
\end{array}\right)
\right\|_{\CO_\delta,2\beta_+,\tilde r} &<&\frac{\varepsilon^{\frac{61}{32}}}{2}.\label{Ltilde_p_tilde_q}
\end{eqnarray}
Since $\|p\|_{\CO,\beta, r}$, $\|q\|_{\CO,\beta, r}<\frac{\varepsilon}{10}$, \re{tilde_p_tilde_q} implies a rough estimate for $\tilde p$ and $\tilde q$:
\begin{equation}\label{rough_tilde_pq}
\|\tilde p\|_{\CO_\delta,2\beta_+,\tilde r}, \|\tilde q\|_{\CO_\delta,2\beta_+,\tilde r}<\frac{\varepsilon^{\frac{61}{32}}}{3}+\frac{\varepsilon^{\frac{31}{32}}}{20}<\frac{\varepsilon^{\frac{31}{32}}}{16}.
\end{equation}

\smallskip

As in Subsection \ref{sub_product_preserving}, with the well-defined fourth root
$$
\Theta(\xi\eta):=\left((e^{\frac{\rm i}2\alpha(\xi\eta)}+p_{0,1}(\xi\eta))(e^{-\frac{\rm i}2\alpha(\xi\eta)}+\bar p_{0,1}(\xi\eta)) \right)^{\frac14},
$$
we define
\begin{equation}\label{varphi_p_01}
\varphi:(\xi,\eta)\mapsto \left(\Theta(\xi\eta)\xi, \Theta^{-1}(\xi\eta)\eta\right).
\end{equation}
Since $|p_{0,1}|_{\omega,\beta_+}<\varepsilon$, we can apply Proposition \ref{prop_transform_product} with $r=\tilde r$, $r'=r_+$, and get
$$(\varphi^{-1}\circ\tilde\tau_1\circ\varphi)(\xi,\eta)=\left(\begin{array}{l}
e^{\frac{\rm i}2\alpha_+(\xi\eta)}\eta+p_+(\xi,\eta) \\
e^{-\frac{\rm i}2\alpha_+(\xi\eta)}\xi+q_+(\xi,\eta)
\end{array}
\right),$$
with $\alpha_+\in{\CA}^\R_{2\beta_+,\tilde r}(\CO_\delta)$ satisfying
$$\alpha_+(\xi\eta)-\alpha(\xi\eta)=-{\rm i}(e^{-\frac{\rm i}2\alpha(\xi\eta)}p_{0,1}(\xi\eta)-e^{\frac{\rm i}2\alpha(\xi\eta)}\bar p_{0,1}(\xi\eta)).$$
By Lemma \ref{lemma_real_res}, we obtain \re{error_alpha}.

With $\psi:=\phi\circ\varphi$, which satisfies $\psi\circ\rho=\rho\circ\psi$,
let ${\cL U}=\left(\begin{array}{c}
                                     u \\
                                     v
                                   \end{array}
 \right):=\psi-{\rm Id}$.
In view of the definitions of $\phi={\rm Id}+\hat{\cL U}$ in Lemma \ref{lem_cohomo} and $\varphi$ given in (\ref{varphi_p_01}), we have that
$$(\psi-{\rm Id})(\xi,\eta)= \left(\begin{array}{c}
                                     (\Theta-1) \xi \\
                                     (\Theta^{-1}-1)\eta
                                   \end{array}
 \right)+\left(\begin{array}{c}
                                     {\hat u}(\Theta \xi,\Theta^{-1}\eta) \\
                                     {\hat v}(\Theta \xi,\Theta^{-1}\eta)
                                   \end{array}
 \right).$$
Since $\|\hat u\|_{\CO_\delta,\tilde\beta,  r^{(7)}}$, $\displaystyle \|\hat v\|_{\CO_\delta,\tilde\beta, r^{(7)}} < \frac{\varepsilon^{\frac{49}{50}}}{20}$, by Lemma \ref{lemma_h_Lambda}, we have
\begin{equation}\label{hat_uv_Lambda}
\|{\hat u}(\Theta \xi,\Theta^{-1}\eta)\|_{\CO_\delta,\beta_+,r_+}, \  \|{\hat v}(\Theta \xi,\Theta^{-1}\eta)\|_{\CO_\delta,\beta_+,r_+}<\frac{\varepsilon^{\frac{49}{50}}}{20},
\end{equation}
Moreover, by Lemma \ref{lemma_Lambda_k}, we obtain that
\begin{equation}\label{Lambda-1_Lambda_inver_-1}
\|\Theta-1\|_{\CO_\delta,2\beta_+,\tilde r}, \ \|\Theta^{-1}-1\|_{\CO_\delta,2\beta_+,\tilde r}<\frac34\|p_{0,1}\|_{\CO_\delta,\beta, r}<\frac{3\varepsilon}{40r}.
\end{equation}
Combining \re{hat_uv_Lambda} and \re{Lambda-1_Lambda_inver_-1}, we obtain $\|u\|_{\CO_\delta,\beta_+,r_+}, \ \|v\|_{\CO_\delta,\beta_+,r_+}<\frac{\varepsilon^{\frac{49}{50}}}2$.

It remains to prove \re{esti_p_+q_+} and \re{esti_Lp_+q_+}. By \re{p_+}, \re{q_+} in Proposition \ref{prop_transform_product}, we obtain that $p_+$, $q_+\in\CA_{\beta_+, r_+}(\CO_\delta)$, and
\begin{eqnarray*}
\|p_+\|_{\CO_\delta,\beta_+, r_+}&<&\left(1+\frac34\|p_{0,1}\|_{\CO,\beta,r}\right)\|\tilde p\|_{\CO_\delta,2\beta_+,\tilde r}+\|p_{0,1}\|^2_{\CO,2\beta,r},\\
&<&\left(1+\frac{3\varepsilon}{40r}\right)\left(\frac{\varepsilon^{\frac{61}{32}}}{3}+22(K+1)\delta^{-1}\|e^{\frac{\rm i}2\alpha}\eta q+e^{-\frac{\rm i}2\alpha}\xi p\|_{\CO,\tilde\beta, r^{(7)}}\right)+\varepsilon^2\\
&<&\frac{\varepsilon^{\frac{61}{32}}}{2}+24(K+1)\delta^{-1}\|e^{\frac{\rm i}2\alpha}\eta q+e^{-\frac{\rm i}2\alpha}\xi p\|_{\CO,\tilde\beta, r^{(7)}},\\
\|q_+\|_{\omega,\beta_+, r_+}&<&\frac{\varepsilon^{\frac{61}{32}}}{2}+24(K+1)\delta^{-1}\|e^{\frac{\rm i}2\alpha}\eta q+e^{-\frac{\rm i}2\alpha}\xi p\|_{\CO,\tilde\beta, r^{(7)}}.
\end{eqnarray*}
Moreover, by \re{crossing_theta_+} in Proposition \ref{prop_transform_product} and (\ref{Ltilde_p_tilde_q}), (\ref{rough_tilde_pq}), we have
\begin{eqnarray*}
& &\|e^{-\frac{\rm i}2\alpha_+}\xi p_++e^{-\frac{\rm i}2\alpha_+}\eta q_+\|_{\CO_\delta,\beta_+, r_+}\nonumber\\
&<&\|e^{-\frac{\rm i}2\alpha} \xi \tilde p+ e^{\frac{\rm i}2\alpha} \eta  \tilde q\|_{\CO_\delta,2\beta_+,\tilde r}+\|p_{0,1}\|_{\CO,\beta,r}(\|\tilde p\|_{\CO_\delta,2\beta_+,\tilde r}+\|\tilde q\|_{\CO_\delta,2\beta_+,\tilde r})+\|p_{0,1}\|_{\CO,\beta,r}^2\\
&<& \frac{\varepsilon^{\frac{61}{32}}}{2}+\frac{\varepsilon}{10r}\cdot\frac{\varepsilon^{\frac{31}{32}}}{8}+\frac{\varepsilon^2}{100r^2}<\varepsilon^{\frac{61}{32}}.
\end{eqnarray*}

It remains to show that
\begin{equation}\label{tran_crown}
\psi(\xi,\eta)\in {\cL C}^{r}_{\omega,\beta} \quad {\rm for} \ (\xi,\eta)\in{\cL C}^{r_+}_{\omega,\beta_+}, \quad  \omega\in{\cL O}_{\delta}(r_{+},\beta_{+}).
\end{equation}
Recalling the definition of the set given in (\ref{cwb}), we have
$${\cL C}^{r_+}_{\omega,\beta_+}:=\left\{(\xi,\eta)\in \C^2: |\xi\eta -\omega|<\beta_+,  \;\ |\xi|, |\eta|<r_+ \right\}.$$
Since $\|u\|_{\CO_\delta,\beta_+,r_+}$, $\displaystyle \|v\|_{\CO_\delta,\beta_+,r_+}<\frac{\varepsilon^{\frac{49}{50}}}2$,
we have that
\begin{eqnarray*}
|(\xi+u(\xi,\eta))(\eta+v(\xi,\eta)) -\omega|
&<&\beta_++\|\eta u+\xi v+uv\|_{\omega,\beta_+,r_+}\\
&<&\beta_++\varepsilon^{\frac{49}{50}}+\varepsilon^{\frac{49}{25}}<\beta,\\
|\xi+u(\xi,\eta)|, \ |\eta+v(\xi,\eta)|
&<&r_++\varepsilon^{\frac{49}{50}}<r.
\end{eqnarray*}
(\ref{tran_crown}) is shown.

\appendix

\section{Proof of Lemma \ref{lem_norm}}\label{proof_lemma_norm}

Let $\varsigma:=\max\{\|f_1-f_2\|_{\CO,\beta'',r''},\ \|g_1-g_2\|_{\CO,\beta'',r''}\}$,
which is, by (\ref{f_g_m}), smaller than $\frac{\beta'^2}{16}$.
Let $\omega\in {\CO}(r'', \beta'')$. In order to estimate its norm, let us first decompose the following expression:
\begin{eqnarray}
& &h(e^{{\rm i}b\alpha(\xi\eta)}\xi+f_1,e^{-{\rm i}b\alpha(\xi\eta)}\eta+g_1)-h(e^{{\rm i}b\alpha(\xi\eta)}\xi+f_2,e^{-{\rm i}b\alpha(\xi\eta)}\eta+g_2) \nonumber \\
&=& \sum_{l\geq0}\left(h_{l,0}(\xi\eta+e^{-{\rm i}b\alpha}\eta f_1+e^{{\rm i}b\alpha}\xi g_1+f_1g_1)-h_{l,0}(\xi\eta+e^{-{\rm i}b\alpha}\eta f_2+e^{{\rm i}b\alpha}\xi g_2+f_2g_2)\right)\label{error0_}\\
& & \  \  \  \  \  \   \cdot \, \left(e^{{\rm i}b \alpha}\xi+f_2\right)^l\nonumber \\
& &+\, \sum_{l\geq1}  h_{l,0}(\xi\eta+e^{-{\rm i}b\alpha}\eta f_1+e^{{\rm i}b\alpha}\xi g_1+f_1g_1) \, \left((e^{{\rm i}b\alpha}\xi+f_1)^l- (e^{{\rm i}b\alpha}\xi+f_2)^l \right)\label{error2_}\\
&&+ \,  \text{similar expressions involving $h_{0,j}$ instead of $h_{l,0}$.}\label{error-sim}
\end{eqnarray}

\begin{itemize}
\item {\bf Terms in (\ref{error0_})} 
\end{itemize}

Expanding $h_{l,0}$ around $\xi\eta+e^{-{\rm i}b\alpha}\eta f_2+e^{{\rm i}b\alpha}\xi g_2+f_2 g_2$ in (\ref{error0_}), we obtain~:
\begin{eqnarray}
 & & h_{l,0}(\xi\eta+e^{-{\rm i}b\alpha}\eta f_1+e^{{\rm i}b\alpha}\xi g_1+f_1 g_1)-h_{l,0}(\xi\eta+e^{-{\rm i}b\alpha}\eta f_2+e^{{\rm i}b\alpha}\xi g_2+f_2 g_2)  \label{develop_h_l0}\\
 &=& \sum_{k\geq 1}\frac{1}{k!}h^{(k)}_{l,0}(\xi\eta+e^{-{\rm i}b\alpha}\eta f_2+e^{{\rm i}b\alpha}\xi g_2+f_2 g_2)\nonumber\\
  & &  \  \  \  \  \  \  \  \cdot \left((e^{-{\rm i}b\alpha}\eta f_1+e^{{\rm i}b\alpha}\xi g_1+f_1 g_1)-(e^{-{\rm i}b\alpha}\eta f_2+e^{{\rm i}b\alpha}\xi g_2+f_2 g_2)\right)^k.\nonumber
\end{eqnarray}
Combining \nrc{lem_eiba} together with Remark \ref{rem_eiba} and (\ref{f_g_m}), we obtain, for $(\xi,\eta)\in {\cL C}^{r''}_{\om,\beta''}$ and $m=1,2$~:
\begin{eqnarray*}
|\xi\eta+e^{-{\rm i}b\alpha}\eta f_m+e^{{\rm i}b\alpha}\xi g_m+f_m g_m-\omega|
&\leq&|\xi\eta-\omega|+|e^{-{\rm i}b\alpha}\eta f_m+e^{{\rm i}b\alpha}\xi g_m+f_m g_m|\\
&<&\beta''+2e^{\frac98\tilde\beta}r''\cdot \frac{\beta'^2}{16}+\frac{\beta'^4}{256}\\
&<&\frac{\beta'}{2}+\frac{101}{200}\cdot \frac{\beta'}{16}+\frac{\beta'}{256}<\beta'.
\end{eqnarray*}
By Cauchy's inequality, we have, for all $ (\xi,\eta)\in {\cL C}_{\omega,\beta''}$ and $\tilde k\geq 1$,
\begin{equation}\label{Cauchy_h_l0}
\frac{1}{\tilde k!}\left|h^{(\tilde k)}_{l,0}(\xi\eta)\right|\leq \sup_{|z-\xi\eta|=\frac{\beta'}{2}}|h_{l,0}(z)|\left(\frac{2}{\beta'}\right)^{\tilde k} \leq  \left(\frac{2}{\beta'}\right)^{\tilde k}|h_{l,0}|_{\omega,\beta'}.
\end{equation}
We recall that, for $|z|<1$,
$$
	\sum_{\tilde k\geq k}C_{\tilde k}^k z^{\tilde k-k}
	=\frac{1}{k!}\sum_{\tilde k\geq 0} \frac{d^k}{dz^k}\left(z^{\tilde k}\right)=\frac{1}{k!}\frac{d^k}{dz^k}\left(\frac{1}{1-z}\right)
	=(1-z)^{-(k+1)}.
$$
Hence, developing $h^{(k)}_{l,0}$ around $\xi\eta$, we have, for $k\geq 1$,
\begin{eqnarray}
& &
\frac{1}{k!}\left\|h^{(k)}_{l,0}(\xi\eta+e^{-{\rm i}b\alpha}\eta f_2+e^{{\rm i}b\alpha}\xi g_2+f_2 g_2)\right\|_{\omega,\beta'',r''}\nonumber\\
&\leq& 
\sum_{\tilde k\geq k} \frac{\left|h^{(\tilde k)}_{l,0}\right|_{\omega,\beta''}}{k!(\tilde k-k)!} \|e^{-{\rm i}b\alpha}\eta f_2+e^{{\rm i}b\alpha}\xi g_2+f_2 g_2\|^{\tilde k-k}_{\omega,\beta'',r''}\nonumber\\
&=& \sum_{\tilde k\geq k} \frac{\tilde k!}{(\tilde k-k)!\cdot k!}\cdot \frac{\left|h^{(\tilde k)}_{l,0}\right|_{\omega,\beta''}}{\tilde k!} \|e^{-{\rm i}b\alpha}\eta f_2+e^{{\rm i}b\alpha}\xi g_2+f_2 g_2\|^{\tilde k-k}_{\omega,\beta'',r''}\nonumber\\
&\leq & |h_{l,0}|_{\omega,\beta'}\sum_{\tilde k\geq k}C_{\tilde k}^k \left(\frac{2}{\beta'}\right)^{\tilde k} \left(\frac{\beta'^2}{16}\right)^{\tilde k-k}
= \left(\frac{2}{\beta'}\right)^{k} |h_{l,0}|_{\omega,\beta'}\sum_{\tilde k\geq k}C_{\tilde k}^k \left(\frac{\beta'}{8}\right)^{\tilde k-k}\nonumber\\
&=&\left(\frac{2}{\beta'}\right)^{k} |h_{l,0}|_{\omega,\beta'}\left(1-\frac{\beta'}{8}\right)^{-(k+1)}.\label{hk_k!}
\end{eqnarray}
Recalling that $|\om|<r''^2-\beta''$, by \rl{lem_norm00}, we have, for $-1\leq b\leq 1$,
\begin{align}
	 &\left\|(e^{-{\rm i}b\alpha}\eta f_1+e^{{\rm i}b\alpha}\xi g_1+f_1g_1)-(e^{-{\rm i}b\alpha}\eta f_2+e^{{\rm i}b\alpha}\xi g_2+f_2g_2)\right\|_{\omega,\beta'',r''} \nonumber   \\
	\leq&\left\|e^{-{\rm i}b\alpha}\eta (f_1-f_2)\right\|_{\omega,\beta'',r''}+ \left\|e^{{\rm i}b\alpha}\xi(g_1-g_2)\right\|_{\omega,\beta'',r''} \nonumber  \\
	& +  \|(f_1-f_2)g_1-f_2(g_1-g_2)\|_{\omega,\beta'',r''} \nonumber \\
	< &e^{\frac98\tilde\beta}r''\left(\|f_1-f_2\|_{\omega,\beta'',r''}+\|g_1-g_2\|_{\omega,\beta'',r''}\right) \nonumber \\
	 & +  \, \|f_1-f_2\|_{\omega,\beta'',r''}\|g_1\|_{\omega,\beta'',r''}+\|f_2\|_{\omega,\beta'',r''}\|g_1-g_2\|_{\omega,\beta'',r''}  \nonumber \\
	< & \frac{101}{50} r''\varsigma+\frac{\beta'^2}{8}\varsigma<\frac{13}{25}\varsigma.\label{eiba_fg}
\end{align}
\vspace{-.1cm}
Since $\frac{\varsigma}{\beta'}<\frac{\beta'}8$ and according to \re{tildebeta}, we have $1-\frac{\beta'}{8}>\frac{99}{100}$ and  $1-\frac{26\varsigma}{25\beta'}\left(1-\frac{\beta'}{8}\right)^{-1}>\frac{99}{100}$, so that
 $\frac{26}{25}\left(1-\frac{\beta'}{8}\right)^{-2}\left(1-\frac{26\varsigma}{25\beta'}\left(1-\frac{\beta'}{8}\right)^{-1}\right)^{-1}<\frac{27}{25}.$
Combining (\ref{hk_k!}) and (\ref{eiba_fg}), together with (\ref{esi_coeff}), we obtain
\begin{align}
	& \left\|h_{l,0}(\xi\eta+e^{-{\rm i}b\alpha}\eta f_1+e^{{\rm i}b\alpha}\xi g_1+f_1 g_1)-h_{l,0}(\xi\eta+e^{-{\rm i}b\alpha}\eta f_2+e^{{\rm i}b\alpha}\xi g_2+f_2 g_2)\right\|_{\omega,\beta'',r''}\nonumber\\
	\leq&\sum_{k\geq 1}\frac{1}{k!}\|h^{(k)}_{l,0}(\xi\eta+e^{-{\rm i}b\alpha}\eta f_2+e^{{\rm i}b\alpha}\xi g_2+f_2 g_2)\|_{\omega,\beta'',r''}\nonumber\\
 & \, \cdot\left\|(e^{-{\rm i}b\alpha}\eta f_1+e^{{\rm i}b\alpha}\xi g_1+f_1g_1)-(e^{-{\rm i}b\alpha}\eta f_2+e^{{\rm i}b\alpha}\xi g_2+f_2g_2)\right\|_{\omega,\beta'',r''}^k \nonumber\\
	<& \, \|h\|_{\omega,\beta',r'} r'^{-l} \cdot  \sum_{k\geq 1}\left(\frac{26\varsigma}{25\beta'}\right)^k \left(1-\frac{\beta'}{8}\right)^{-(k+1)} \nonumber\\
	=& \, \|h\|_{\omega,\beta',r'} r'^{-l}\frac{26\varsigma}{25\beta'}  \left(1-\frac{\beta'}{8}\right)^{-2}\left(1-\frac{26\varsigma}{25\beta'}\left(1-\frac{\beta'}{8}\right)^{-1}\right)^{-1}  \nonumber\\
	<& \, \frac{27\varsigma}{25\beta'}\|h\|_{\CO,\beta',r'} r'^{-l}.\label{h_lj}
\end{align}
Hence, according to Lemma \ref{lem_norm00}, for $l\geq 0$,
\begin{eqnarray*}
& &\left\|\left(h_{l,0}(\xi\eta+e^{-{\rm i}b\alpha}\eta f_1+e^{{\rm i}b\alpha}\xi g_1+f_1 g_1)-h_{l,0}(\xi\eta+e^{-{\rm i}b\alpha}\eta f_2+e^{{\rm i}b\alpha}\xi g_2+f_2 g_2)\right)\right.\\
& & \  \  \  \  \  \  \  \  \cdot \, \left.(e^{{\rm i}b\alpha} \xi+f_2)^l\right\|_{\omega,\beta'',r''}\\
&<&\frac{27\varsigma}{25\beta'}\|h\|_{\omega,\beta',r'} r'^{-l}\left(e^{\frac98 \tilde\beta} r''+\frac{\beta'^2}{16}\right)^{l}.
\end{eqnarray*}
On the other hand, (\ref{small_beta02}) implies $2r''\tilde\beta <\frac{r'-r''}{8}$. Indeed,
	since $0<r''<r'<\frac14$ and $0<\beta<1$, then $8\beta^{\frac{1}{2}}<(r'-r'')r''$ implies
\begin{equation}\label{tildebeta_r'_r''}
2r''\tilde\beta=2r''\cdot 16\beta^{\frac54}<8\beta^{\frac54}<8\beta<8\beta^{\frac{1}{2}}\cdot \frac{(r'-r'')r''}{8}<\frac{r'-r''}{8}.
\end{equation}
Therefore, according to \re{exptbeta}, we have
	\begin{equation*}
	r'-e^{\frac98\tilde\beta} r''-\frac{\beta'^2}{16}>r'-r''-2\tilde\beta r''>r'-r''-\frac{r'-r''}{8}=\frac{7}{8}(r'-r'').
	\end{equation*}
	As a consequence, we have
$$\sum_{k\geq 0}r'^{-k}\left(e^{\frac98 \tilde\beta} r''+\frac{\beta'^2}{16}\right)^{k} =\frac{r'}{r'-e^{\frac98\tilde\beta} r''-\frac{\beta'^2}{16}}<\frac{8r'}{7(r'-r'')}.$$
Thus, under (\ref{small_beta02}), the $\|\cdot\|_{\CO,\beta'',r''}-$norm of (\ref{error0_}) is bounded by
\begin{eqnarray}
\frac{27\varsigma}{25\beta'}\|h\|_{\CO,\beta',r'}\sum_{l\geq 0}r'^{-l}\left(e^{\frac98 \tilde\beta} r''+\frac{\beta'^2}{16}\right)^{l}
&<& \frac{27}{25}\cdot\frac{8r'}{7(r'-r'')}\frac{\varsigma}{\beta'}\|h\|_{\CO,\beta',r'}\nonumber\\
&<& \frac{5r'}{4(r'- r'')}\frac{\varsigma}{\beta'}\|h\|_{\CO,\beta',r'}.\label{part0}
\end{eqnarray}

To show $(\ref{error0_})\in {\CA}_{\beta'', r''}\left({\CO}\right)$ provided that $f_1,f_2,g_1,g_2\in {\CA}_{\beta'', r''}\left({\CO}\right)$, it remains to verify the analyticity on ${\CO}(r'',\beta'') $ for the coefficients $(\ref{error0_})_{l,j}$, $l,j\geq 0$, $lj=0$.
According to \re{develop_h_l0}, \re{hk_k!} and \re{eiba_fg}, for $\tilde l\geq 0$,
$$h_{\tilde l,0}(\xi\eta+e^{-{\rm i}b\alpha}\eta f_1+e^{{\rm i}b\alpha}\xi g_1+f_1g_1)-h_{\tilde l,0}(\xi\eta+e^{-{\rm i}b\alpha}\eta f_2+e^{{\rm i}b\alpha}\xi g_2+f_2g_2)
\in {\CA}_{\beta'',r''}\left({\CO}\right),$$
and, by (\ref{esi_coeff_uniform}) and (\ref{h_lj}), for $\tilde l\geq 0$, for $l,j\geq 0$ with $lj=0$,
\begin{eqnarray*}
& &\left|\left(h_{\tilde l,0}(\xi\eta+e^{-{\rm i}b\alpha}\eta f_1+e^{{\rm i}b\alpha}\xi g_1+f_1g_1)\right.\right.\\
& &\left.\left. \  \  \  - \, h_{\tilde l,0}(\xi\eta+e^{-{\rm i}b\alpha}\eta f_2+e^{{\rm i}b\alpha}\xi g_2+f_2g_2) \right)_{l,j}\right|_{{\CO}(r'',\beta'')}
\leq\frac{27\varsigma}{25\beta'}\frac{\|h\|_{\CO,\beta',r'}}{ r'^{\tilde l}r''^{l+j}}.
\end{eqnarray*}
Note that, for $l\geq 0$,
\begin{align*}
(\ref{error0_})_{l,0}= \,&\sum_{\tilde l\geq 0}\sum_{0\leq k \leq l}\left((e^{{\rm i}b\alpha}\xi+f_2)^{\tilde l}\right)_{k,0}\\
 & \  \  \ \cdot \left(h_{\tilde l,0}(\xi\eta+e^{-{\rm i}b\alpha}\eta f_1+e^{{\rm i}b\alpha}\xi g_1+f_1g_1)-h_{\tilde l,0}(\xi\eta+e^{-{\rm i}b\alpha}\eta f_2+e^{{\rm i}b\alpha}\xi g_2+f_2g_2)\right)_{l-k,0} \\
 & + \, \sum_{\tilde l\geq 0} \sum_{k \geq 1}\left((e^{{\rm i}b\alpha}\xi+f_2)^{\tilde l}\right)_{l+k,0} (\xi\eta)^k\\
 & \  \  \   \cdot
\left(h_{\tilde l,0}(\xi\eta+e^{-{\rm i}b\alpha}\eta f_1+e^{{\rm i}b\alpha}\xi g_1+f_1g_1)-h_{\tilde l,0}(\xi\eta+e^{-{\rm i}b\alpha}\eta f_2+e^{{\rm i}b\alpha}\xi g_2+f_2g_2)\right)_{0,k},\end{align*}
and for $j\geq 1$,
\begin{align*}(\ref{error0_})_{0,j}=&\sum_{\tilde l\geq 0}\sum_{0\leq k \leq j}\left((e^{{\rm i}b\alpha}\xi+f_2)^{\tilde l}\right)_{0,k}\\
&  \  \   \cdot \left(h_{\tilde l,0}(\xi\eta+e^{-{\rm i}b\alpha}\eta f_1+e^{{\rm i}b\alpha}\xi g_1+f_1g_1)-h_{\tilde l,0}(\xi\eta+e^{-{\rm i}b\alpha}\eta f_2+e^{{\rm i}b\alpha}\xi g_2+f_2g_2)\right)_{0,j-k}\\
&  + \, \sum_{\tilde l\geq 0} \sum_{k \geq 1}\left((e^{{\rm i}b\alpha}\xi+f_2)^{\tilde l}\right)_{0,j+k} (\xi\eta)^k\\
&  \  \  \  \cdot
\left(h_{\tilde l,0}(\xi\eta+e^{-{\rm i}b\alpha}\eta f_1+e^{{\rm i}b\alpha}\xi g_1+f_1g_1)-h_{\tilde l,0}(\xi\eta+e^{-{\rm i}b\alpha}\eta f_2+e^{{\rm i}b\alpha}\xi g_2+f_2g_2)\right)_{k,0},\end{align*}
where, by Lemma \ref{lem_eiba} and  by (\ref{esi_coeff_uniform}), for $l,j\geq 0$ with $lj=0$,
$$\left|\left((e^{{\rm i}b\alpha}\xi +f_2)^{\tilde l}\right)_{l,j}\right|_{{\CO}(r'',\beta'')}\leq \left(e^{\frac98\tilde\beta}r''+\frac{\beta'^2}{16}\right)^{\tilde l}r''^{-(l+j)}.$$
Then we see that, for $\omega\in {\CO}(r'',\beta'') $,
\begin{eqnarray*}
& &\sum_{\tilde l\geq 0}\sum_{0\leq k \leq l}\left|\left((e^{{\rm i}b\alpha}\xi+f_2)^{\tilde l}\right)_{k,0}(\omega)\right.\\
& & \  \  \  \cdot \left. \left(h_{\tilde l,0}(\xi\eta+e^{-{\rm i}b\alpha}\eta f_1+e^{{\rm i}b\alpha}\xi g_1+f_1g_1)-h_{\tilde l,0}(\xi\eta+e^{-{\rm i}b\alpha}\eta f_2+e^{{\rm i}b\alpha}\xi g_2+f_2g_2)\right)_{l-k,0}(\omega)\right| \\
&\leq&\frac{27(l+1)\varsigma\|h\|_{\CO,\beta',r'}}{25\beta' r''^{l} }\sum_{\tilde l\geq 0}\left(e^{\frac98\tilde\beta}\frac{r''}{r'}+\frac{\beta'^2}{16r'}\right)^{\tilde l},
\end{eqnarray*}
\begin{eqnarray*}
& &\sum_{\tilde l\geq 0} \sum_{k \geq 1}\left|\left((e^{{\rm i}b\alpha}\xi+f_2)^{\tilde l}\right)_{l+k,0}(\omega)\cdot \omega^k\right.\\
& & \  \  \   \cdot \left.
\left(h_{\tilde l,0}(\xi\eta+e^{-{\rm i}b\alpha}\eta f_1+e^{{\rm i}b\alpha}\xi g_1+f_1g_1)-h_{\tilde l,0}(\xi\eta+e^{-{\rm i}b\alpha}\eta f_2+e^{{\rm i}b\alpha}\xi g_2+f_2g_2)\right)_{0,k}(\omega)\right|\\
&\leq&\sum_{\tilde l\geq 0}\sum_{k\geq 1} \left(e^{\frac98\tilde\beta}r''+\frac{\beta'^2}{16}\right)^{\tilde l}r''^{-(l+k)}(r''^2-\beta'')^k\cdot \frac{27\varsigma}{25\beta'}\frac{\|h\|_{\CO,\beta',r'}}{ r'^{\tilde l}r''^{k}}\\
&\leq& \frac{27\varsigma\|h\|_{\CO,\beta',r'}}{25\beta' r''^{l} } \sum_{k\geq 1}\left(1-\frac{\beta''}{r''^2} \right)^k \sum_{\tilde l\geq 0}\left(e^{\frac98\tilde\beta}\frac{r''}{r'}+\frac{\beta'^2}{16r'}\right)^{\tilde l},
\end{eqnarray*}
which implies the analyticity of $(\ref{error0_})_{l,0}$ under (\ref{small_beta02}), and it is similar for that of $(\ref{error0_})_{0,j}$.
Hence, $(\ref{error0_})\in {\CA}_{\beta'', r''}\left({\CO}\right)$.


\begin{itemize}
\item {\bf Terms in (\ref{error2_})}
\end{itemize}

Note that (\ref{esi_coeff}) and (\ref{Cauchy_h_l0}) imply that, for $\omega\in {\CO}(r'',\beta'')$, for $l\geq 0$,
\begin{gather}
	\|h_{l,0}(\xi\eta+e^{-{\rm i}b\alpha}\eta f_1+e^{{\rm i}b\alpha}\xi g_1+f_1g_1) \|_{\omega,\beta'',r''}\leq\sum_{k\geq 0}\frac{|h^{(k)}_{l,0}|_{\omega,\beta''}}{k!}\|e^{-{\rm i}b\alpha}\eta f_1+e^{{\rm i}b\alpha}\xi g_1+f_1g_1 \|^k_{\omega,\beta'',r''}\nonumber\\
	\leq|h_{l,0}|_{\omega,\beta'}\sum_{k\geq 0}\left(\frac{2}{\beta'}\right)^{k}\left(2e^{\frac98\tilde\beta}r'' \frac{\beta'^2}{8}+\frac{\beta'^4}{256}\right)^k\leq\frac{101}{100}\frac{\|h\|_{\CO,\beta',r'}}{r'^l}.\label{h_l0_comp}
\end{gather}
In (\ref{error2_}), we have, for $l\geq 1$, \begin{eqnarray}
\|(e^{{\rm i}b\alpha}\xi+f_1)^l- (e^{{\rm i}b\alpha}\xi+f_2)^l\|_{\omega,\beta'', r''}&\leq &\sum_{k=1}^l C_{l}^k \|e^{{\rm i}b\alpha}\xi+f_1\|_{\omega,\beta'', r''}^{l-k} \|f_2-f_1\|_{\omega,\beta'', r''}^k\nonumber\\
&<&\sum_{k=1}^l C_{l}^k \left(e^{\frac98\tilde\beta} r''+\frac{\beta'^2}{16}\right)^{l-k} \varsigma^k\nonumber\\
&=&\left(e^{\frac98\tilde\beta}r''+\frac{\beta'^2}{16}+\varsigma\right)^l- \left(e^{\frac98\tilde\beta} r''+\frac{\beta'^2}{16}\right)^{l}\nonumber\\
&\leq&l \left(e^{\frac98\tilde\beta}r''+\beta'^2+\varsigma\right)^{l-1} \varsigma.\label{eiba_comp}
\end{eqnarray}
Furthermore, by (\ref{small_beta02}), we have
 $\beta^{\frac54}<\beta^{\frac12}<\frac{r'-r''}{32}$. Recalling that $\varsigma<\frac{\beta'^2}{8}\leq\frac{(16\beta^\frac54)^2}{8}$ and using \re{tildebeta_r'_r''}, we have
 \begin{equation}\label{e98betar}
 1-\frac{e^{\frac98\tilde\beta}r''+\beta'^2+\varsigma}{r'} > 1-\frac{(1+2\tilde\beta) r''}{r'}
 > \frac{r'-r''}{r'}-\frac{r'-r''}{8r'}
 =\frac{7(r'-r'')}{8r'}.
 \end{equation}
Therefore, the $\|\cdot\|_{\CO,\beta'', r''}-$norm of (\ref{error2_}) is bounded by
\begin{gather}
	\sum_{l\geq1} \|h_{l,0}(\xi\eta+e^{-{\rm i}b\alpha}\eta f_1+e^{{\rm i}b\alpha}\xi g_1+f_1g_1)\|_{\omega,\beta'', r''} \cdot \,  \| (e^{{\rm i}b\alpha}\xi+f_1)^l- (e^{{\rm i}b\alpha}\xi+f_2)^l \|_{\CO,\beta'', r''}\nonumber\\
	< \frac{101\varsigma}{100r'}\|h\|_{\CO,\beta',r'}\sum_{l\geq1}l \left(\frac{e^{\frac98\tilde\beta}r''+\beta'^2+\varsigma}{r'}\right)^{l-1}=\frac{101\varsigma}{100r'}\|h\|_{\CO,\beta',r'} \left(1-\frac{e^{\frac98\tilde\beta}r''+\beta'^2+\varsigma}{r'}\right)^{-2} \nonumber\\
	< \frac{101\varsigma}{100r'}\|h\|_{\CO,\beta',r'}\cdot\frac{8^2 r'^2}{7^2(r'-r'')^2} <\frac{7r'\varsigma\|h\|_{\CO,\beta',r'}}{5(r'-r'')^2}.\label{part2}
\end{gather}

For $l\geq 0$, we have that $(\ref{error2_})_{l,0}$ equals to
\begin{align*}
&\sum_{\tilde l\geq 1}\sum_{0\leq k \leq l} \left(h_{\tilde l,0}(\xi\eta+e^{-{\rm i}b\alpha}\eta f_1+e^{{\rm i}b\alpha}\xi g_1+f_1g_1)\right)_{l-k,0}\left( (e^{{\rm i}b\alpha}\xi+f_1)^{\tilde l}-  (e^{{\rm i}b\alpha}\xi+f_2)^{\tilde l}\right)_{k,0}\\
 & + \, \sum_{\tilde l\geq 1}\sum_{k \geq 1} \left(h_{\tilde l,0}(\xi\eta+e^{-{\rm i}b\alpha}\eta f_1+e^{{\rm i}b\alpha}\xi g_1+f_1g_1)\right)_{l+k,0}  \left( (e^{{\rm i}b\alpha}\xi+f_1)^{\tilde l}-  (e^{{\rm i}b\alpha}\xi+f_2)^{\tilde l}\right)_{0,k}(\xi\eta)^k\\
 & + \, \sum_{\tilde l\geq 1} \sum_{k \geq 1} \left(h_{\tilde l,0}(\xi\eta+e^{-{\rm i}b\alpha}\eta f_1+e^{{\rm i}b\alpha}\xi g_1+f_1g_1)\right)_{0,k} \left( (e^{{\rm i}b\alpha}\xi+f_1)^{\tilde l}-  (e^{{\rm i}b\alpha}\xi+f_2)^{\tilde l}\right)_{l+k,0}(\xi\eta)^k,
\end{align*}
and for $j\geq 1$, $(\ref{error2_})_{0,j}$ equals to
\begin{align*}
&\sum_{\tilde l\geq 1}\sum_{0\leq k \leq j} \left(h_{\tilde l,0}(\xi\eta+e^{-{\rm i}b\alpha}\eta f_1+e^{{\rm i}b\alpha}\xi g_1+f_1g_1)\right)_{0,j-k}  \left( (e^{{\rm i}b\alpha}\xi+f_1)^{\tilde l}-  (e^{{\rm i}b\alpha}\xi+f_2)^{\tilde l}\right)_{0,k}\\
 & + \, \sum_{\tilde l\geq 1}\sum_{k \geq 1} \left(h_{\tilde l,0}(\xi\eta+e^{-{\rm i}b\alpha}\eta f_1+e^{{\rm i}b\alpha}\xi g_1+f_1g_1)\right)_{0,j+k}  \left( (e^{{\rm i}b\alpha}\xi+f_1)^{\tilde l}-  (e^{{\rm i}b\alpha}\xi+f_2)^{\tilde l}\right)_{k,0} (\xi\eta)^k\\
 & + \, \sum_{\tilde l\geq 1} \sum_{k \geq 1} \left(h_{\tilde l,0}(\xi\eta+e^{-{\rm i}b\alpha}\eta f_1+e^{{\rm i}b\alpha}\xi g_1+f_1g_1)\right)_{k,0}  \left( (e^{{\rm i}b\alpha}\xi+f_1)^{\tilde l}-  (e^{{\rm i}b\alpha}\xi+f_2)^{\tilde l}\right)_{0,j+k} (\xi\eta)^k.
\end{align*}
If $f_1,f_2,g_1,g_2\in {\CA}_{\beta'', r''}\left({\CO}\right)$, then by (\ref{h_l0_comp}) -- (\ref{e98betar}) and (\ref{esi_coeff}), we see the analyticity of $(\ref{error2_})_{l,0}$, since for $\omega\in {\CO}(r'',\beta'') $,
\begin{align*}
 &\left|\sum_{\tilde l\geq 1}\sum_{0\leq k \leq l} \left(h_{\tilde l,0}(\xi\eta+e^{-{\rm i}b\alpha}\eta f_1+e^{{\rm i}b\alpha}\xi g_1+f_1g_1)\right)_{l-k,0}(\omega)  \left( (e^{{\rm i}b\alpha}\xi+f_1)^{\tilde l}-  (e^{{\rm i}b\alpha}\xi+f_2)^{\tilde l}\right)_{k,0}(\omega)\right|\\
\leq&\sum_{\tilde l\geq 1}\sum_{0\leq k \leq l}\frac{101}{100}\frac{\|h\|_{\CO,\beta',r'}}{r'^{\tilde l}r''^{l-k}}\cdot \frac{{\tilde l} (e^{\frac98\tilde\beta}r''+\beta'^2+\varsigma)^{{\tilde l}-1} \varsigma}{r''^{k}}\\
=&\frac{101}{100}\frac{ (l+1) \|h\|_{\CO,\beta',r'}  \varsigma }{r' r''^{l}} \sum_{\tilde l\geq 1} {\tilde l} \left(\frac{e^{\frac98\tilde\beta}r''+\beta'^2+\varsigma}{r'}\right)^{\tilde l-1},
\end{align*}
\begin{align*}
 &\left|\sum_{\tilde l\geq 1}\sum_{k \geq 1} \left(h_{\tilde l,0}(\xi\eta+e^{-{\rm i}b\alpha}\eta f_1+e^{{\rm i}b\alpha}\xi g_1+f_1g_1)\right)_{l+k,0} (\omega) \left( (e^{{\rm i}b\alpha}\xi+f_1)^{\tilde l}-  (e^{{\rm i}b\alpha}\xi+f_2)^{\tilde l}\right)_{0,k}(\omega)\cdot \omega^k\right|\\
\leq&\sum_{\tilde l\geq 1}\sum_{k \geq 1}  \frac{101}{100}\frac{\|h\|_{\CO,\beta',r'}}{r'^{\tilde l}r''^{l+k}}\cdot \frac{{\tilde l} (e^{\frac98\tilde\beta}r''+\beta'^2+\varsigma)^{{\tilde l}-1} \varsigma}{r''^{k}}\cdot (r''^2-\beta'')^k\\
=&\frac{101}{100}\frac{  \|h\|_{\CO,\beta',r'}  \varsigma }{r' r''^{l}} \sum_{\tilde l\geq 1} {\tilde l} \left(\frac{e^{\frac98\tilde\beta}r''+\beta'^2+\varsigma}{r'}\right)^{\tilde l-1}\sum_{k \geq 1} \frac{(r''^2-\beta'')^k}{r''^{2k}},
\end{align*}
\begin{align*}
 &\left|\sum_{\tilde l\geq 1}\sum_{k \geq 1} \left(h_{\tilde l,0}(\xi\eta+e^{-{\rm i}b\alpha}\eta f_1+e^{{\rm i}b\alpha}\xi g_1+f_1g_1)\right)_{0,k} (\omega)    \left( (e^{{\rm i}b\alpha}\xi+f_1)^{\tilde l}-  (e^{{\rm i}b\alpha}\xi+f_2)^{\tilde l}\right)_{l+k,0}(\omega)\cdot \omega^k\right|\\
\leq&\sum_{\tilde l\geq 1}\sum_{k \geq 1}  \frac{101}{100}\frac{\|h\|_{\CO,\beta',r'}}{r'^{\tilde l}r''^{k}}\cdot \frac{{\tilde l} (e^{\frac98\tilde\beta}r''+\beta'^2+\varsigma)^{{\tilde l}-1} \varsigma}{r''^{l+k}}\cdot (r''^2-\beta'')^k\\
=&\frac{101}{100}\frac{  \|h\|_{\CO,\beta',r'}  \varsigma }{r' r''^{l}} \sum_{\tilde l\geq 1} {\tilde l} \left(\frac{e^{\frac98\tilde\beta}r''+\beta'^2+\varsigma}{r'}\right)^{\tilde l-1}\sum_{k \geq 1} \frac{(r''^2-\beta'')^k}{r''^{2k}}.
\end{align*}
The proof for $(\ref{error2_})_{0,j}$ is similar, hence $(\ref{error2_})\in{\CA}_{\beta'',r''}\left({\CO}\right)$.


\smallskip

Combining (\ref{part0}),  (\ref{part2}) and similar estimates obtained for expressions \re{error-sim}, we obtain
\begin{eqnarray*}
& &\|h(e^{{\rm i}b\alpha}\xi+f_1,e^{-{\rm i}b\alpha}\eta+g_1)-h(e^{{\rm i}b\alpha}\xi+f_2,e^{-{\rm i}b\alpha}\eta+g_2)\|_{\CO,\beta'', r''}\\
&<&\frac{1}{r'- r''} \left(\frac{14r'}{5(r'- r'')}+\frac{5r'}{2\beta'}\right)\varsigma\|h\|_{\CO,\beta',r'}  < \, \frac{3r'}{(r'-r'')\beta'}\varsigma\|h\|_{\CO,\beta',r'},
\end{eqnarray*}
since (\ref{small_beta02}) implies that $\beta'\leq\tilde\beta<\frac{r'-r''}{64}$.
This finishes the proof of Lemma \ref{lem_norm}.\qed

\def\cprime{$'$}

\end{document}